\documentclass[reqno,12pt]{amsart}

\usepackage[utf8]{inputenc}
\usepackage[T1]{fontenc}
\usepackage{lmodern}

\usepackage[
backend=biber,   
style=numeric,
giveninits=true,
isbn=false,
url=false, 
sorting=anyt,
citestyle=numeric-comp
]{biblatex}
\usepackage{hyperref}
\usepackage{graphicx}
\usepackage[labelfont={rm}, justification=centering]{subfig}

\usepackage{stmaryrd}           
\usepackage{mathtools}          
\usepackage{amssymb}
\usepackage{bm}                 
\usepackage{accents}            
\usepackage{xfrac}              

\usepackage[foot]{amsaddr}         
\usepackage[svgnames]{xcolor}
\usepackage{todonotes}

\usepackage{textcomp}          

\usepackage{array}             
\usepackage{booktabs}          

\usepackage{algorithm}         
\usepackage{algpseudocode}     

\usepackage{lipsum}

\definecolor{Mosss}{RGB}{0,192,81}
\definecolor{Reddd}{RGB}{240,40,40}
\definecolor{Mocca}{RGB}{148,82,0}
\newcommand{\RI}  [1]{#1}
\newcommand{\RII} [1]{#1}
\newcommand{\RIII}[1]{#1}
\newcommand{\RIV} [1]{#1}
\newcommand{\JS}  [1]{#1}

\calclayout

\newcommand{\doi}[1]{\textsc{doi}: \href{http://dx.doi.org/#1}{\nolinkurl{#1}}}

\DeclareMathAlphabet{\mathbsl}{OT1}{cmr}{bx}{sl}  

\newcommand {\mathsc}[1]{\text{{\rmfamily\scshape #1}}}

\newcommand {\EX}  {\mathsc{ex}}                       
\newcommand {\IM}  {\mathsc{im}}                       

\newcommand {\fc}  {f_\mathrm{c}}
\newcommand {\fd}  {f_\mathrm{d}}

\newcommand {\aim} {a^{\IM}} 
\newcommand {\bim} {b^{\IM}} 
\newcommand {\aex} {a^{\EX}} 
\newcommand {\bex} {b^{\EX}} 

\newcommand {\V}  [1] {\bm{#1}}                        
\newcommand {\T}  [1] {\bm{#1}}                        
\newcommand {\NM} [1] {\underaccent{\bar}{#1}}         

\renewcommand {\d}            {\partial}
\newcommand   {\D}            {\mathrm{d}\mspace{1.0mu}}

\newcommand   {\imag}     [1] {\operatorname{Im}(#1)}
\newcommand   {\real}     [1] {\operatorname{Re}(#1)}
\newcommand   {\transpose}[1] {#1^{\textsc{t}}}

\newcommand   {\trun}         {t_{\mathrm{run}}}
\newcommand   {\dt}           {\Delta t}

\newcommand   {\Clmb}         {\mathit{CFL}_{\lambda}}
\newcommand   {\CClmb}        {\mathit{CFL}_{\lambda}^{\ast}}

\newcommand   {\Emin}         {\varepsilon_{v,\min}}

\begin{document}

\title%
{
Assessment of high-order IMEX methods for incompressible flow
}
\author{Montadhar Guesmi$^\dag$, Martina Grotteschi$^\ddag$, J\"org Stiller$^\ddag$}
\address{$^\dag$
  TU Dresden, Institute of Process Engineering and Environmental Technology, 
  01062 Dresden, Germany}
\address{$^\ddag$
  TU Dresden, Institute of Fluid Mechanics, 01062 Dresden, Germany}
\email{\{montadhar.guesmi,martina.grotteschi,joerg.stiller\}@tu-dresden.de}

\begin{abstract}
This paper investigates the competitiveness of semi-implicit Runge-Kutta (RK) and spectral deferred correction (SDC) time-integration methods up to order six for incompressible Navier-Stokes problems in conjunction with a high-order discontinuous Galerkin method for space discretization.
It is proposed to harness the implicit and explicit RK parts as a partitioned scheme, which provides a natural basis for the underlying projection scheme and yields a straight-forward approach for accommodating nonlinear viscosity. 
Numerical experiments on laminar flow, variable viscosity and transition to turbulence are carried out to assess accuracy, convergence and computational efficiency. 
Although the methods of order 3 or higher are susceptible to order reduction due to time-dependent boundary conditions, two third-order RK methods are identified that perform well in all test cases and clearly surpass all second-order schemes including the popular extrapolated backward difference method.
The considered SDC methods are more accurate than the RK methods, but become competitive only for relative errors smaller than ca $10^{-5}$.
\end{abstract}


\maketitle


\section{Introduction}
\label{sec:intro}

\RIII{%
The present study is concerned with the adaptation and application of higher-order accurate time integration methods to the incompressible Navier-Stokes equations with possibly time-dependent boundary conditions and variable viscosity.}
High-order discretization methods enjoy increasing popularity 
\RIII{in computational fluid dynamics}
because of their low dispersion errors and superior convergence properties\RIII{, e.g.,}
\cite{%
SE_Beck2014a,%
SE_Schaal2015a,%
TI_Tavelli2016a,%
SE_Winters2018a,%
SE_Fehn2019a,%
SE_Bassi2020a}.
\RIII{%
However, while orders up to 16 are common in space, the use of low-order time integration schemes is still widespread.}
As the resolution of simulations increases, this approach requires extremely small time steps to match the temporal and spatial accuracy.
High-order time-integration methods promise higher accuracy with larger steps and lower cost, which shall be assessed in the present paper.

While explicit methods offer the simplest approach, they become infeasible with growing spatial order due to the worsening stiffness of the \RIII{viscous} terms 
\cite{SE_Karniadakis2005a,SE_Canuto2011a}.
This leaves the choice between implicit and semi-implicit methods, depending on the treatment of the convective terms.
Stable implicit methods exist for virtually any order and possess the advantage of choosing the time step only according to accuracy or physical criteria
\cite{TI_Hairer1996a}.
Several authors used Diagonally Implicit Runge-Kutta (DIRK) methods for simulation
of compressible flows, e.g. 
\cite{%
  TI_Jothiprasad2003a,%
  TI_Persson2010a,%
  TI_Uranga2010a,%
  TI_Pazner2017a,%
  TI_Bassi2015a,%
  TI_Pan2021a}.
The application of \JS{fully} implicit Runge-Kutta (IRK) methods to incompressible flows is more complicated due to the quasi-static nature of the continuity equation, but was investigated in %
\RII{%
\cite{%
  TI_John2006a,%
  TI_Sanderse2013a,%
  TI_Bassi2015a,%
  TI_Noventa2016a}}.
In most of these studies, the order of convergence reached 3 to 5, although methods up to order 10 were considered in \cite{TI_Southworth2021a}.
\RII{%
\textcite{TI_Sanderse2013a} developed energy conserving IRK methods for incompressible Navier-Stokes problems based on implicit collocation schemes.}
As an alternative to IRK, variational methods such as Discontinuous Galerkin methods in time
and space-time element methods are gaining interest
\cite{TI_Tavelli2016a,TI_Ahmed2017a}. 
%
\JS{%
They use a polynomial expansion in time, which makes them inherently implicit and capable of reaching any order.}
A common drawback of implicit methods is the need to solve nonlinear systems in each step.
\JS{%
These systems are linearized e.g. using fixed point or Newton methods, or already in the context of the RK method, as in Rosenbrock methods
\cite{TI_John2006a,TI_Bassi2015a,TI_Noventa2016a}.}
Unfortunately, the resulting linear systems require computation and storage of the exact or approximate Jacobian of the nonlinear operator.
Moreover, they are typically non-symmetric and often ill-conditioned.
Their fast solution depends on an efficient preconditioner, which may be hard to achieve.
These difficulties diminish the benefits of implicit methods considerably.

Semi-implicit methods combine an implicit scheme for diffusion with an explicit one for convection and, hence, are often called IMEX methods.
The explicit treatment of convection terms leads to a Courant-Friedrichs-Lewy (CFL) condition, which depends on the type and the order of the spatial discretization.
For methods based on polynomial elements, the admissible time step scales approximately as ${\Delta t \sim \sfrac{\Delta x}{v P^2}}$ where $\Delta x$ is the element length, $v$ the fluid velocity and $P$ the polynomial degree, see e.g.\ 
\cite{SE_Karniadakis2005a}.
The implicit discretization of diffusive terms leads to symmetric linear systems, for which very efficient solution methods exist, e.g.\ 
\cite{%
  SE_Lottes2005a,%
  SE_Janssen2011a,%
  SE_Stiller2016b,%
  SE_Stiller2017a}.
This offers an enormous advantage over fully implicit methods, especially in high fidelity simulations, where the time step is limited by accuracy rather than stability constraints.
Due to their simplicity, IMEX variants of linear multistep methods are popular, especially for incompressible flow simulation
\cite{%
  TI_Karniadakis1991a,%
  TI_Leriche2006a,%
  TI_Klein2015a,%
  TI_Fehn2017a}.
Corresponding methods can be constructed for any order, however, increasing the latter beyond two results in the loss of A-stability for the implicit part and worsens the stability restrictions for the explicit part 
\cite{TI_Hairer1996a}.
As an alternative, high-order IMEX Runge-Kutta (RK) methods have been proposed, e.g.~%
\RIII{%
\cite{%
  TI_Ascher1997a,%
  TI_Kennedy2003a,%
  TI_Pareschi2005a,%
  TI_Boscarino2009a,%
  TI_Boscarino2013a,%
  TI_Cavaglieri2015a,%
  TI_Boscarino2017a,%
  TI_Kennedy2019a}.~}%
These methods have far better stability properties and a higher accuracy than the best IMEX multistep schemes of the same order.
%
\RII{%
IMEX RK methods 
were applied to
compressible
\cite{TI_Kanevsky2007a,TI_Persson2011a}
as well as incompressible flow problems
\cite{TI_Le1991a,TI_Nikitin2006a,TI_Colomes2016a}.
The second-order scheme of \textcite{TI_Le1991a} is probably the first IMEX RK method specifically designed for incompressible Navier-Stokes problems and has been widely used for the simulation of turbulent flows.
\textcite{TI_Nikitin2006a} derived a third-order IMEX method by adding an implicit perturbation to an explicit RK scheme.
\textcite{TI_Colomes2016a} proposed segregated IMEX RK methods which, unlike earlier methods, are not affected by splitting errors.
This is achieved by incorporating the pressure terms in to the explicit RK part.
A possible disadvantage of these methods is the lack a mechanism to enforce continuity, which can be a source of instability in the simulation of turbulent flows.}

Unfortunately, the growing number of order conditions makes it difficult to construct IMEX RK methods of orders greater than 4. 
To the knowledge of the authors only one method of order 5 has been published
\cite{TI_Kennedy2003a}.
In contrast, the spectral deferred correction (SDC) method can reach arbitrary convergence rates by incrementally increasing the order using sweeps of lower order schemes 
\cite{TI_Dutt2000a,TI_Ong2020a}.
\textcite{TI_Minion2003b} developed a semi-implicit spectral deferred correction (SISDC) method based on an IMEX Euler scheme.
This method has stability properties that are comparable with IMEX RK methods, but extends easily to higher order.
It was adopted for incompressible flows in periodic domains in \cite{TI_Minion2004a} and later extended to time-dependent boundary conditions
\cite{TI_Minion2018a}.
In \cite{TI_Stiller2020a} the SISDC method was generalized to variable properties and demonstrated to reach convergence rates up to order 12 for three-dimensional vortical flows with nonlinear (turbulent) viscosity.
\RIV{%
A further alternative is the use of ADER methods, which attain an order of two in space and time and can cope with compressible as well as incompressible flow including variable diffusivity and density \cite{TI_Busto2018a,TI_Bermudez2020a}.
}

\RIII{%
A common issue of high-order time-integration methods is their susceptibility to order reduction, which means that in certain situations the observed convergence rate is lower than the theoretical order of convergence.
Possible reasons for this behavior include 
the stiffness of the problem,
the presence of variables with no evolution equation as in differential-algebraic systems,
time-dependent boundary conditions, and
lacking smoothness of the source terms or the solution itself,
see, e.g., \cite{TI_Kennedy2016a,TI_Minion2018a} and the references cited therein.
The role of stiffness in order reduction was analyzed theoretically in
\cite{TI_Boscarino2007a} for IMEX RK methods and in 
\cite{TI_Boscarino2018a} for integral deferred correction methods, which can be regarded as
as a variant of the SDC method.
For incompressible Navier-Stokes problems, \textcite{TI_Minion2018a} observed an order reduction of their SDC method due to time-dependent boundary conditions.
In \cite{TI_Stiller2020a} it was shown that this issue is related to the divergence error of the approximate velocity and can be significantly reduced by adding a grad-div stabilization to the numerical scheme.
}

The obvious shortcomings of low-order IMEX multistep methods lead to the question whether high-order RK and SDC methods can surpass and replace them as a suitable complement to higher-order spatial discretization methods.
Facing this challenge, the present study investigates the convergence properties and computational efficiency of selected IMEX RK methods as stand-alone integrators and as predictors in SISDC methods for incompressible Navier-Stokes problems with constant or variable viscosity.
To cope with the quasi-static incompressibility constraint and achieve a consistent treatment of nonlinear viscosity, it is proposed to harness the explicit and implicit RK components as a partitioned RK method.
\RII{In contrast to earlier work, the proposed method first applies the explicit scheme to all but the source terms and then improves the result using a projection method based on the DIRK part.}
Spatial discretization is accomplished using a discontinuous Galerkin method of order 16.
The convergence behavior of IMEX RK and SISDC methods up to order 6 is studied using various test cases featuring time-dependent boundary conditions, nonlinear viscosity and turbulence.
For each case an error-cost analysis is performed to assess the computational efficiency.
A comparison with the widespread IMEX BDF2 method demonstrates the superiority of the proposed high-order time-integration methods, even for direct simulations of turbulence with marginal spatial resolution.

The remainder of the paper is organized as follows: 
Section2 describes the time-integration methods using the example of a one-dimensional convection-diffusion equation with variable diffusivity.
Section~3 extends these methods to incompressible Navier-Stokes problems,
Section~4 presents the numerical experiments including a discussion of the results and
Section~5 concludes the paper.


\section{Time-integration methods}
\label{sec:time-integration}


\subsection{Model problem}
\label{sec:time-integration:model-problem}

For explaining the IMEX approach and devising the semi-implicit treatment of nonlinear viscosity, the one-dimensional convection-diffusion equation is considered as a model problem:
\begin{equation}
  \label{eq:conv-diff:pde}
  \d_t u + \d_x (v u) = \d_x (\nu \d_x u)
  \,.
\end{equation}
Here, the solution $u(x,t)$ is a function in space and time,
while the convection velocity $v$ and diffusivity $\nu$ may
depend on $u$.
The simplest case with constant velocity and diffusivity possesses the spatially periodic solution
${u = \hat u(t) e^{\mathrm i k x}}$ 
with wave number $k$.
Its amplitude $\hat u$ satisfies the ordinary differential equation (ODE) 
\begin{equation}
  \label{eq:conv-diff:ode}
  \D_t \hat u 
    = \lambda_{\mathrm c} \hat u + \lambda_{\mathrm d} \hat u
    = \lambda \hat u
  \,, 
\end{equation}
where
${\lambda_{\mathrm c} = -\mathrm i vk}$ and 
${\lambda_{\mathrm d} = -\nu k^2}$
represent the effect of convection and diffusion, respectively.
The associated time scales
${t_{\mathrm c} = |\lambda_{\mathrm c}|^{-1} = \sfrac{1}{v k}}$
and
${t_{\mathrm d} = |\lambda_{\mathrm d}|^{-1} = \sfrac{1}{\nu k^2}}$ 
motivate a semi-implicit treatment as in the case of flow problems, 
since the $t_{\mathrm d}$ decreases much faster with growing wave number
than $t_{\mathrm c}$.
Therefore, this problem serves as a model for characterizing the stability and accuracy of the investigated IMEX methods in Sec.~\ref{sec:time-integration:stability}.

A more general situation arises when a spatial discretization method is applied to the original problem \eqref{eq:conv-diff:pde}. 
In this case, the semi-discrete solution takes the form of a time dependent coefficient vector, $\NM u(t)$ and is governed by the ODE system
\begin{equation}
  \label{eq:conv-diff:ode system}
  \D_t \NM u = \NM C(\NM v) \NM u + \NM D(\NM \nu) \NM u
  \,,
\end{equation}
where 
$\NM C$ and $\NM D$ are the discrete convection and diffusion operators,
and
$\NM v$ and $\NM \nu$ the coefficient vectors of velocity and diffusivity,
both of which may depend on $\NM u$.
As the time integration methods act individually on each component of \eqref{eq:conv-diff:ode system}, their description is based on the corresponding scalar problem
\begin{equation}
  \label{eq:model}
  \d_t u = \fc(u) + \fd(\nu,u) = f(u)
  \,,
\end{equation}
where
$\fc$ corresponds to ${\NM C(\NM v) \NM u}$ and
$\fd$ to ${\NM D(\NM \nu) \NM u}$, respectively. 
Note that the explicit dependence of $\fc$ on $v$ has been dropped for convenience.


\subsection{IMEX methods for constant diffusivity}
\label{sec:time-integration:imex-methods}

As a starting point, the time-integration methods are introduced for constant diffusivity.
This case allows for a clean IMEX approach, where 
the diffusive part $\fd$ is treated implicitly and
the convective part $\fc$ explicitly. 
Using equidistant steps ${\Delta t}$ yields the times ${t^n = n \Delta t}$ for which the approximate solution ${u^n \approx u(t^n)}$ is sought.
The IMEX backward difference formula of order 2 (IMEX BDF2) serves as a reference and is given by \cite{TI_Frank1997a}
\begin{equation}
  \label{eq:bdf2:constant-nu}
    \frac{\gamma_0 u^{n+1} - \alpha_0 u^{n} - \alpha_1 u^{n-1}}{\Delta t}
  = \beta_0 \fc(u^{n}) + \beta_1 \fc(u^{n-1}) + \fd(\nu,u^{n+1})
  \,,
\end{equation}
where
${\alpha_0 =  2}$,
${\alpha_1 = -\JS{1}/2}$,
${\beta_0  =  2}$,
${\beta_1  = -1}$ and
${\gamma_0 =  \JS{3}/2}$.
Since the method is not self-starting, it is initialized using an IMEX Euler step, 
i.e.
${\alpha_0 = \beta_0 = \gamma_0 = 1}$ and
${\alpha_1 = \beta_1 = 0}$.

\RIII{%
The IMEX RK methods investigated in this study are combinations of DIRK and explicit Runge-Kutta (ERK) schemes with an identical number of stages $s$, possibly different coefficients
${\NM b^\IM}$ 
and 
${\NM b^\EX}$, 
but equal time nodes 
${\NM c^\IM = \NM c^\EX = \NM c}$.
Noting that the latter assumption implies 
${\aim_{1,1} = c_1 = 0}$, the Butcher tableaux can be written in the form:}
\renewcommand{\arraystretch}{1.5}%
\begin{equation}
  \label{eq:imex-rk:butcher}
  \begin{array}{c|ccccc}
    0      & 0          &            &              &            \\
    c_2    & \aim_{2,1} & \aim_{2,2} &              &            \\
    \vdots & \vdots     & \!\ddots   & \!\ddots     &            \\
    c_{s}  & \aim_{s,1} & \cdots     & \aim_{s,s-1} & \aim_{s,s} \\\hline
           & \bim_1     & \cdots     & \bim_{s-1}   & \bim_s     \\
  \end{array}
  \qquad
  \begin{array}{c|ccccc}
    0      & 0          &            &              &        \\
    c_2    & \aex_{2,1} & 0          &              &        \\
    \vdots & \vdots     & \!\ddots   & \ddots       &        \\
    c_{s}  & \aex_{s,1} & \cdots     & \aex_{s,s-1} & 0      \\\hline
           & \bex_1     & \cdots     & \bex_{s-1}   & \bex_s \\
  \end{array}
\end{equation}%
\renewcommand{\arraystretch}{1.0}%
\RIII{%
IMEX RK methods with this structure are referred to as CK methods \cite{TI_Boscarino2007a}.
More general methods with ${\NM c^\IM \ne \NM c^\EX}$ may achieve favorable properties
(see e.g.\cite{TI_Pareschi2005a}), but are not considered here because the synchronization of the time nodes is crucial for the splitting methods used in this study.~}%

To reach a certain order \JS{of convergence}, the DIRK and ERK methods must satisfy the corresponding order and coupling conditions \cite{TI_Kennedy2003a}.
Table~\ref{tab:rk properties} provides an overview of the considered methods, including the number of stages, the theoretical order of convergence and further relations between the coefficients.
The Butcher tableaux are given in Appendix~\ref{sec:butcher}. 
\JS{%
Additionally, the IMEX RK methods selected for this study satisfy ${c_s = 1}$ and possess
an L-stable DIRK part.}

\begin{table}[h]
  \caption{Properties of IMEX Runge-Kutta methods.}
  \label{tab:rk properties}
  \begin{tabular}{l@{\qquad}c@{\qquad}c@{\qquad}c@{\qquad}c@{\qquad}c}
    \toprule
    method & stages & 
    order & 
    $a^\IM_{s,j}=b^\IM_j$ &
    $a^\EX_{s,j}=b^\EX_j$ &
    $b^\EX_i=b^\IM_i$  \\\midrule
    RK-TR   & 3 & 2 & $\checkmark$ & $\checkmark$ &              \\
    RK-CB2  & 3 & 2 & $\checkmark$ &              & $\checkmark$ \\
    RK-CB3c & 4 & 3 &              &              & $\checkmark$ \\
    RK-CB3e & 4 & 3 & $\checkmark$ &              & $\checkmark$ \\
    RK-CB4  & 6 & 4 & $\checkmark$ &              & $\checkmark$ \\
    RK-ARS3 & 5 & 3 & $\checkmark$ & $\checkmark$ &              \\ 
    \bottomrule
  \end{tabular}
\end{table}
RK-TR is a 3-stage version of the IMEX trapezoidal rule with better stability properties than the more common 2-stage version.
RK-CB2, RK-CB3c, RK-CB3e and RK-CB4 were developed by \textcite{TI_Cavaglieri2015a}, who named them IMEXRKCB2, IMEXRKCB3c, IMEXRKCB3e and IMEXRKCB4, respectively.
The assumption of identical assembly coefficients $\NM b^\EX = \NM b^\IM$ reduced the number of coupling conditions and thus allowed to optimize these methods for stability and accuracy.
Moreover, the DIRK part of RK-CB4 attains stage order 2, which renders it less susceptible to order reduction when applied to stiff problems.
RK-ARS3 was adopted from Ascher, Ruth and Spiteri \cite{TI_Ascher1997a}, who dubbed it (4,4,3), because the implicit and explicit part can be both reduced to 4 stages and the resulting IMEX method has order 3.
\RIII{%
It should be noted that all methods except RK-CB3c satisfy ${\aim_{s,j}=\bim_j}$ and thus have a stiffly accurate (SA) DIRK part \cite{TI_Hairer1996a}.
The explicit parts of RK-TR and RK-ARS3 fulfill the condition ${\aex_{s,j}=\bex_j}$ which along with ${c_s = 1}$ yields the first-same-as-last (FSAL) property \cite{TI_Hairer1993a}. 
Moreover, the combination of SA and FSAL renders these schemes globally stiffly accurate
(GSA) \cite{TI_Boscarino2013a}. 
This property is important in PDE applications such as Navier-Stokes problems because it preserves the boundary conditions enforced in the last stage, which is not guaranteed with non-GSA schemes.
}

\JS{%
Finally, the application to the model problem~\eqref{eq:model} with constant diffusivity is described in Algorithm~\ref{alg:imex-rk:model:constant}.}

\begin{algorithm}[t]
\caption{IMEX Runge-Kutta method for model problem with constant diffusivity.}
\label{alg:imex-rk:model:constant}
\begin{algorithmic}[1]
\Procedure{ImexRK}{$u$}
   \State  $u_1 \gets u$
   \For{$i=2,s$}
      \State 
         ${u_i \gets u + \Delta t
                         \sum_{j=1}^{i-1} 
                            \big[ \aex_{i,j} \fc(u_j) 
                                + \aim_{i,j} \fd(\nu,u_j) 
                            \big]
                       + \Delta t\,\aim_{i,i} f_{\mathrm d}(\nu,u_i)
          }$
   \EndFor
   \State 
     ${u \gets u + \Delta t \sum_{i=1}^{s} 
                               \big[ \bex_{i} \fc(u_i)
                                   + \bim_{i} \fd(\nu,u_i)
                               \big]
     }$
\EndProcedure
\end{algorithmic}
\end{algorithm}


\subsection{Extension to variable diffusivity}
\label{sec:time-integration:variable-nu}

The presence of a variable diffusivity ${\nu(u)}$ renders the diffusion term 
${\fd(\nu,u)}$ nonlinear and, thus, complicates the solution of the implicit 
equations in each time step.
With IMEX BDF2 this difficulty can be circumvented by extrapolating the diffusivity, i.e.
\begin{equation}
  \label{eq:bdf2:variable-nu}
    \frac{\gamma_0 u^{n+1} - \alpha_0 u^{n} - \alpha_1 u^{n-1}}{\Delta t}
  = \beta_0 \fc(u^{n}) + \beta_1 \fc(u^{n-1}) 
  + \fd(\beta_0\nu^n\! + \beta_1\nu^{n-1},u^{n+1})
  \,,
\end{equation}
where ${\nu^n = \nu(u^n)}$ etc.
This approach preserves second order accuracy and leads to discrete equations that are linear in $u^{n+1}$ and, hence, easier to solve as in the fully implicit case.
As a possible side effect, diffusive stability restrictions are to be expected.

Unfortunately, the semi-implicit treatment of nonlinear diffusion is not trivial with
Runge-Kutta methods.
Rosenbrock-type methods can handle variable coefficients, but are expensive due to the introduction of Jacobians and the more complicated structure of linear systems.
\textcite{TI_Boscarino2016a} used partitioned RK methods to construct semi-implicit schemes without the need of nonlinear iterations.
The IMEX RK method for variable diffusivity presented in 
Algorithm~\ref{alg:semi-rk:model} was developed independently, 
but can be seen as a special case of this approach. 
At the begin of stage $i$, it applies the ERK component to compute the preliminary solution $v_i$ (line 6), which is used to evaluate the diffusivity $\nu_i$ (line 7).
The computation of the stage solution $u_i$ then resembles a fully implicit treatment of diffusion, but is actually semi-implicit, since $\nu_i$ is already known (line 8).
As a consequence, each RK stage has a complexity comparable to the IMEX BDF2 method for variable diffusivity stated in \eqref{eq:bdf2:variable-nu}, while the full scheme still retains the accuracy of the underlying IMEX RK method
\cite{TI_Boscarino2016a}.

\begin{algorithm}[t]
\caption{Semi-implicit partitioned Runge-Kutta method for variable diffusivity.}
\label{alg:semi-rk:model}
\begin{algorithmic}[1]
\Procedure{ImexRK}{$u$}
   \State  $u_1 \gets u$
   \State  $v_1 \gets u$
   \State  $\nu_1 \gets \nu(u)$
   \For{$i=2,s$}
      \State  ${v_i \gets u + \Delta t 
                              \sum_{j=1}^{i-1}
                                \aex_{i,j} 
                                \big[ f_{\mathrm c}(u_j)
                                    + f_{\mathrm d}(\nu_j, u_j)
                                \big]}$
      \State  $\nu_i \gets \nu(v_i)$
      \State  ${u_i \gets u + \Delta t 
                              \sum_{j=1}^{i-1}
                                \aex_{i,j} f_{\mathrm c}(u_j)
                            + \Delta t 
                              \sum_{j=1}^{i} 
                                \aim_{i,j} f_{\mathrm d}(\nu_j, u_j)
                                \big]}$
   \EndFor
   \State 
     ${u \gets u + \Delta t 
                   \sum_{i=1}^{s} 
                     \big[ \bex_i f_{\mathrm c}(u_i)
                         + \bim_i f_{\mathrm d}(\nu_i, u_i)
                     \big]}$
\EndProcedure
\end{algorithmic}
\end{algorithm}


\subsection{Spectral deferred correction}
\label{sec:time-integration:sdc}

The SDC method divides the interval ${(t^{n},t^{n+1})}$ into subintervals 
${\{(t_{m-1}, t_m)\}_{m=1}^M}$ 
with intermediate times
${t^{n} \le t_{0} < \dots t_{m-1} < t_m < \dots t_M \le t^{n+1}}$.
Following \cite{TI_Dutt2000a} the Gauss-Lobatto-Legendre (GLL) points are chosen, such that
${t_0 = t^n}$
and 
${t_M = t^{n+1}}$.
A predictor sweep through all subintervals provides the initial approximations
${u_m^0 \approx u(t_m)}$, which are joined in ${\NM u^0}$.
The predictor can be any self-starting time-integration method \cite{TI_Layton2007a}.
In the present study, IMEX Euler and RK methods serve for this purpose.
Starting from ${u_m^0}$, the approximate solution is gradually improved by successive correction sweeps.
Given the $k$-th approximation, ${\NM u^k = \{u_m^k\}}$, the corrector seeks a solution to the error equation \cite{TI_Dutt2000a}
\begin{equation}
  \label{eq:model:error-equation}
  \begin{aligned}
    \delta(t)
      \equiv u(t) - \mathcal I \NM u^k(t)
      = u(t_0) 
    & + \int_{t_0}^{t} \big[ f(\mathcal I \NM u^k(\tau) + \delta(\tau)) 
                           - f(\mathcal I \NM u^k(\tau))
                       \big] \D \tau
    \\
    & + \int_{t_0}^{t} f(\mathcal I \NM u^k(\tau)) \D \tau
      - \mathcal I \NM u^k(t)
      \,,
  \end{aligned}
\end{equation}
where
${\mathcal I \NM u^k(t) = \sum_{m=0}^M u^k_m \ell_m(t)}$ 
is the Lagrange interpolant to $\NM u^k$.
Combining the equations for ${t=t_m}$ and $t_{m+1}$ 
leads to the correction equation
\begin{equation}
  \label{eq:model:correction-equation}
  \begin{aligned}
    u^k_{m+1} + \delta(t_{m+1})
    & = u^k_{m} + \delta(t_{m})
    \\
    & + \int_{t_m}^{t_{m+1}} \big[ f(\mathcal I \NM u^k(\tau) + \delta(\tau)) 
                                 - f(\mathcal I \NM u^k(\tau))
                             \big] \D \tau
    \\
    & + \int_{t_m}^{t_{m+1}} f(\mathcal I \NM u^k(\tau)) \D \tau
    \,.
  \end{aligned}
\end{equation}
As proposed by \textcite{TI_Minion2003b}, 
the IMEX Euler rule is applied to the first integral,
while the second one is approximated using
${f(\mathcal I \NM u^k(\tau)) \approx \mathcal I \NM f^k(\tau)}$
with ${f^k_m = f(u^k_m)}$.
Introducing
${u^{k+1}_m \approx u^{k}_m + \delta(t_m)}$
and exploiting 
${\delta(t_0) = 0}$
finally yields the discrete correction equation
\begin{equation}
  \label{eq:model:corrector}
  \begin{aligned}
    u^{k+1}_{m} 
      = u^{k+1}_{m-1}
      + \Delta t_m [\;\,& \fc(u^{k+1}_{m-1}) + \fd(\nu^{k+1}_{m-1},u^{k+1}_{m}) \\
                   -\, (& \fc(u^{k}_{m-1})   + \fd(\nu^{k}_{m-1},u^{k}_{m}))  \;]
      + \sum_{q=0}^{M} w_{q,m} f(u^k_q)
  \end{aligned}
\end{equation}
with subinterval length ${\Delta t_m = t_{m} - t_{m-1}}$ and weights
${
  w_{q,m} = \int_{t_m}^{t_{m+1}} \!\ell_q(\tau) \D \tau
}$, 
where $\ell_q$ is the Lagrange basis polynomial associated with $t_q$.
Each sweep solves Eq.~\eqref{eq:model:corrector} for ${m = 1}$ to $M$ and
increases the order by one \cite{TI_Dutt2000a}, although stiffness or time-dependent boundary conditions can degrade convergence
\cite{TI_Boscarino2018a,TI_Minion2018a}.
The converged solution satisfies the Lobatto IIIA collocation scheme with $M$ stages and, hence, achieves order $2M$ at the end of the time interval
\cite{TI_Hagstrom2006a}.


\subsection{Linear stability and accuracy}
\label{sec:time-integration:stability}

\begin{figure}
\subfloat[$|R(z)| = 1$]
  {\includegraphics[scale=0.60]{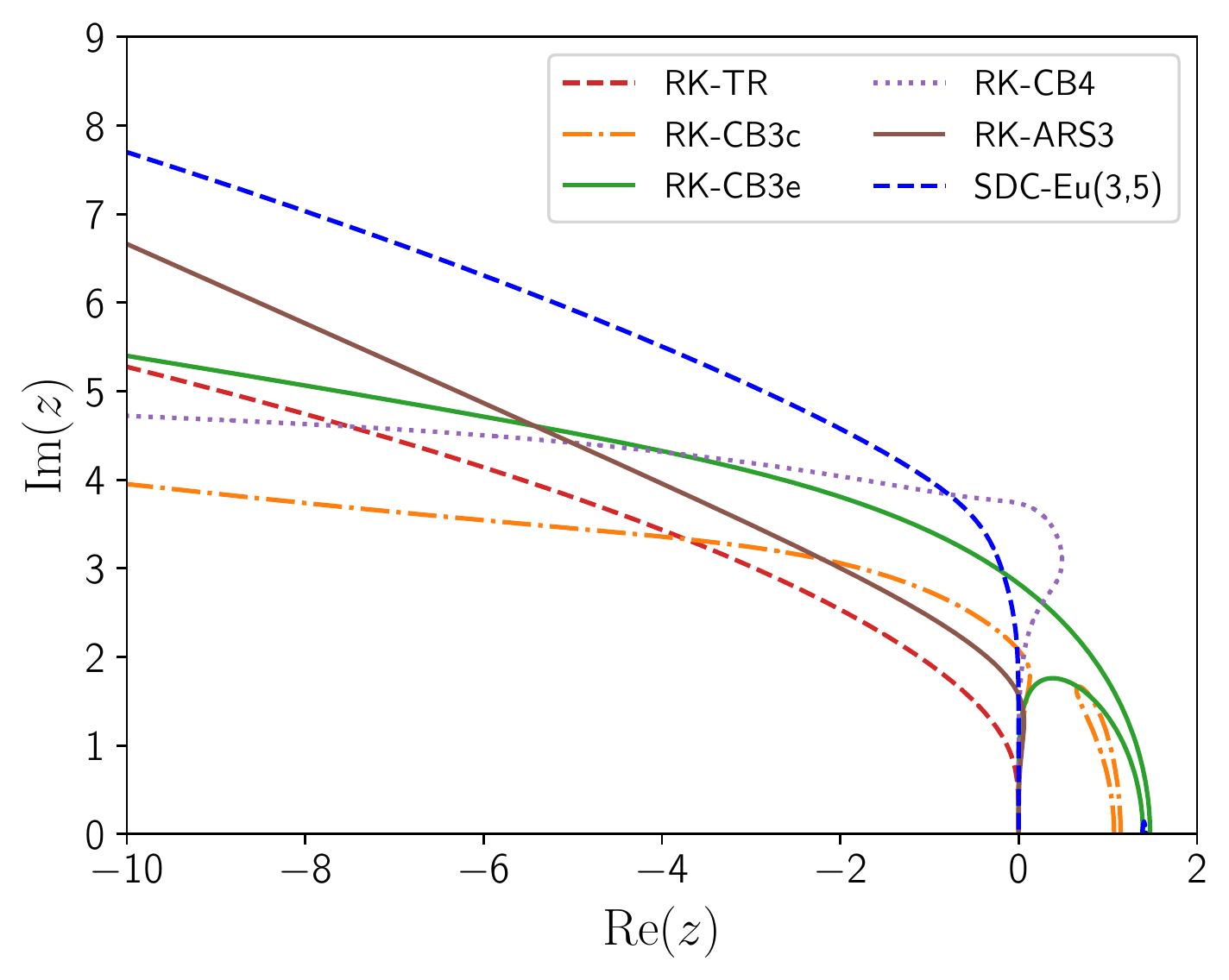}
   \label{fig:model:stability:unscaled}}
\hfill
\subfloat[$|R(z_s)| = 1$]
  {\includegraphics[scale=0.60]{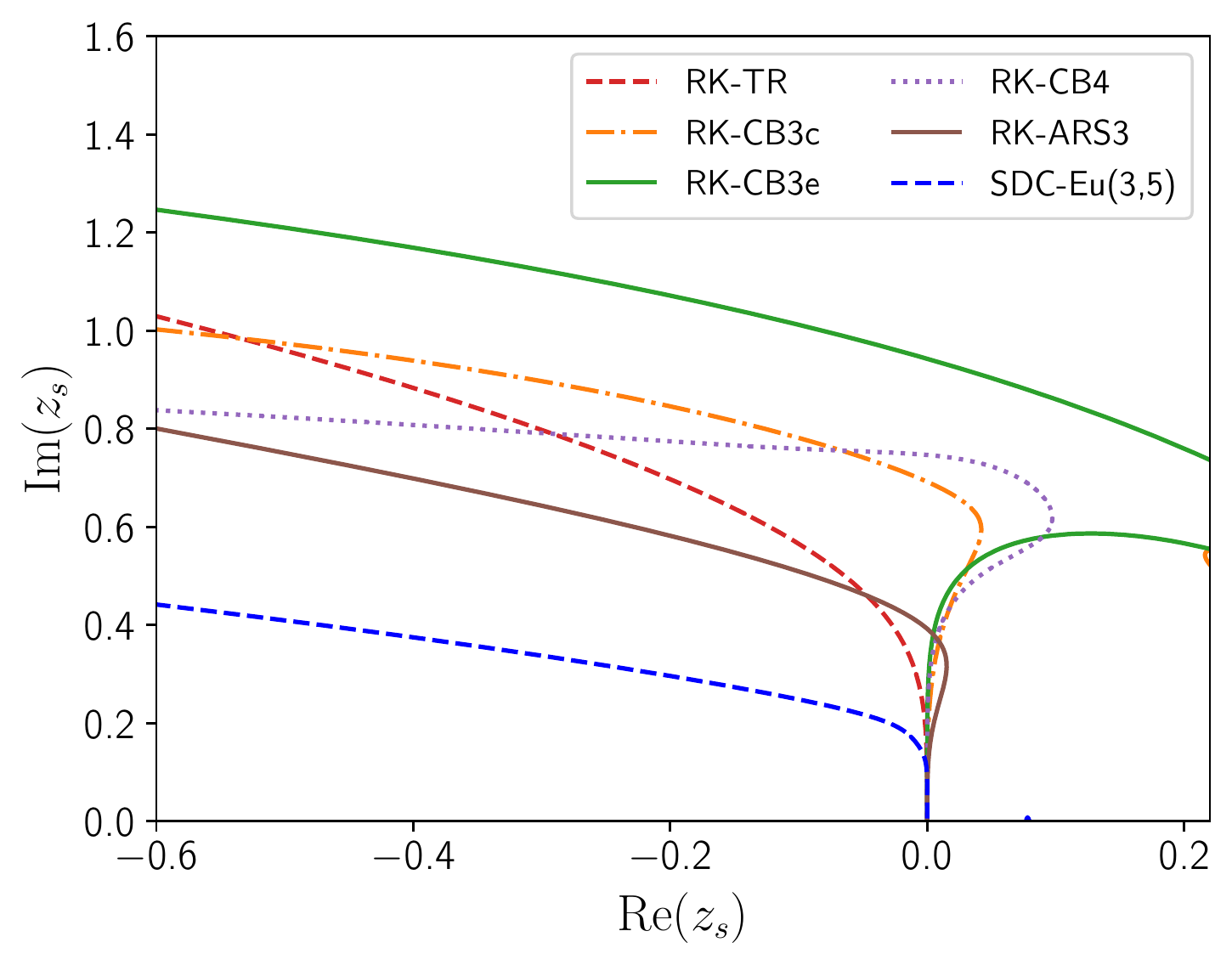}
   \label{fig:model:stability:scaled}}
  \caption{%
  Unscaled and scaled stability domains of selected semi-implicit methods. 
  \label{fig:model:stability}
  }
\end{figure}

For comparing the stability and accuracy of the considered time-integration methods, the  RHS of model problem \eqref{eq:conv-diff:ode} is split into an implicit and an explicit part such that 
\begin{equation}
  \label{eq:model-imex-splitting}
  \D_t \hat u = \lambda^\IM \hat u + \lambda^\EX \hat u
  \,.
\end{equation}
The application of an IMEX one-step method yields the discrete evolution equation
\begin{equation}
  \label{eq:model-imex-one-step}
  \hat u^{n+1} = R(\Delta t\lambda^\IM, \Delta t\lambda^\EX) \hat u^{n}
  \,,
\end{equation}
where ${R}$ is the (unknown) stability function, which depends on two complex parameters,
${z^\IM = \Delta t\lambda^\IM}$ and 
${z^\EX = \Delta t\lambda^\EX}$.
Apart from the general case, the following special choices are possible:
\begin{itemize}
\item
  \makebox[6.5em][l]{implicit:}
  \makebox[4.5em][l]{${\lambda^\IM = \lambda}$,}
  ${\lambda^\EX = 0}$,
\item
  \makebox[6.5em][l]{explicit:}
  \makebox[4.5em][l]{${\lambda^\IM = 0}$,}
  ${\lambda^\EX = \lambda}$,
\item
  \makebox[6.5em][l]{semi-implicit:}
  \makebox[4.5em][l]{${\lambda^\IM = \lambda_{\mathrm d}}$,} 
  ${\lambda^\EX = \lambda_{\mathrm c}}$.
\end{itemize}
Each of these cases yields a specific stability function $R(z)$ that depends only on a single parameter, ${z = \Delta t (\lambda^\IM + \lambda^\EX})$.
Given a one-step IMEX method such as RK or SDC, the corresponding stability function can be evaluated for any $z$ by executing \eqref{eq:model-imex-one-step} once with 
${\Delta t = 1}$, 
${\lambda_{\mathrm c} = \imag z}$,
${\lambda_{\mathrm d} = \real z}$ and
initial condition ${\hat u^0 = 1}$, such that ${R(z) = \hat u^1}$.
Additionally, this procedure provides the error 
${\varepsilon(z) = \hat u^1 - e^z}$.
For a fair comparison of different methods the argument is scaled to the number of nontrivial substeps, i.e.,
${z_s = z / (s-1)}$ for IMEX RK, 
${z_s = z / M(K+1)}$ for SISDC with Euler predictor and
${z_s = z / M(K+s-1)}$ for SISDC with RK predictor.
As $\lambda$ is fixed for a given problem, 
the scaling applies also to the time step such that $\Delta t_s \sim z_s$.
Figure~\ref{fig:model:stability} shows the unscaled and scaled stability domains 
of selected methods for the semi-implicit case 
${\lambda^\IM = \lambda_{\mathrm d}}$ and 
${\lambda^\EX = \lambda_{\mathrm c}}$.
The threshold ${|R(z)| = 1}$ is of interest because the imaginary value $\imag{z}$
corresponds to the eigenvalue-based Courant-Friedrichs-Lewy number, i.e., 
${CFL_{\lambda} = \Delta t |\lambda_{\mathrm c}|}$.
Its critical value depends on the mangnitude of real part 
${|\real z| = \Delta t |\lambda_{\mathrm d}|}$, 
which is sometimes referred to as diffusion number.

\begin{figure}
\subfloat[$|\varepsilon| \le 10^{-5}$]
  {\includegraphics[scale=0.60]{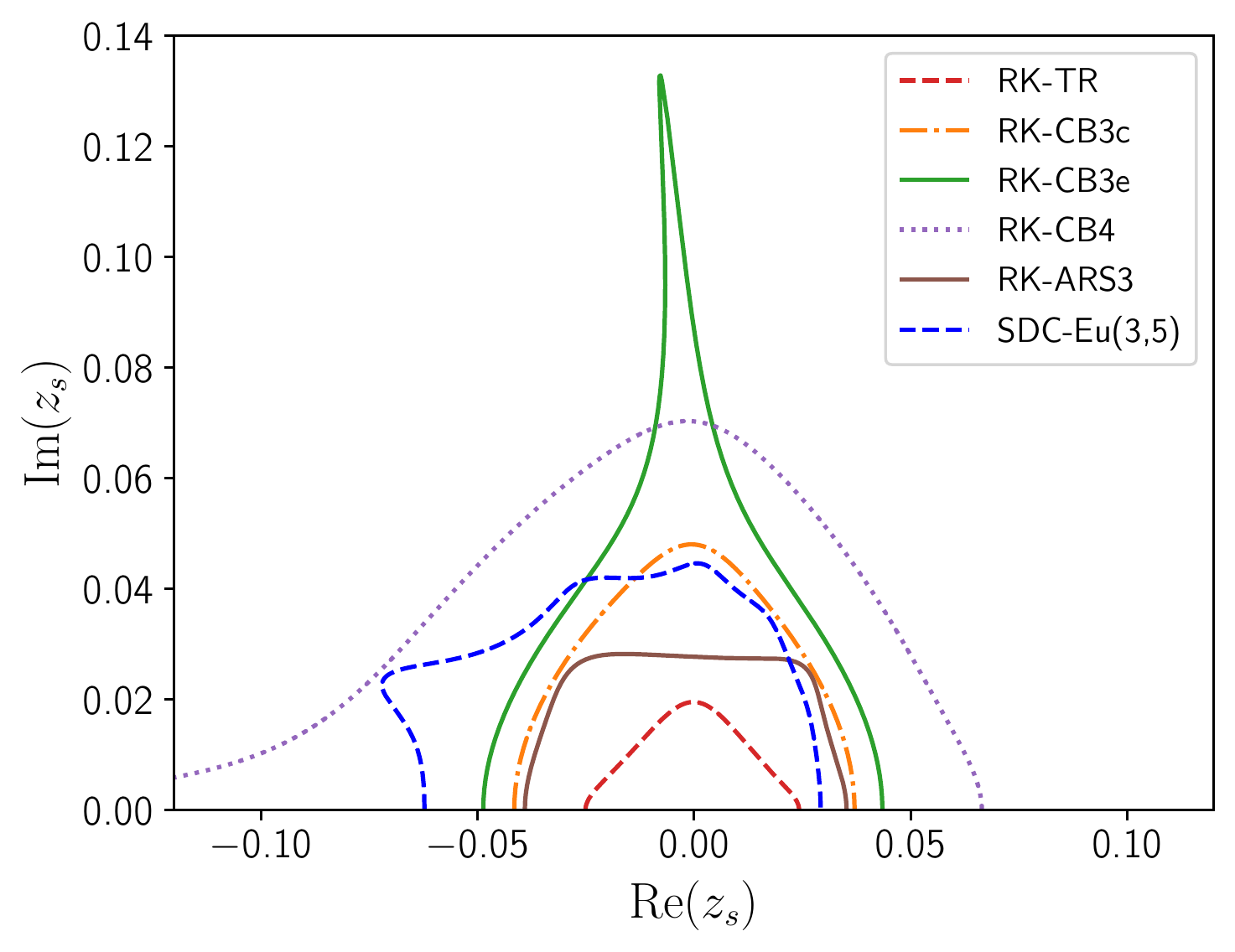}
   \label{fig:model:error-5}}
\hfill
\subfloat[$|\varepsilon| \le 10^{-10}$]
  {\includegraphics[scale=0.60]{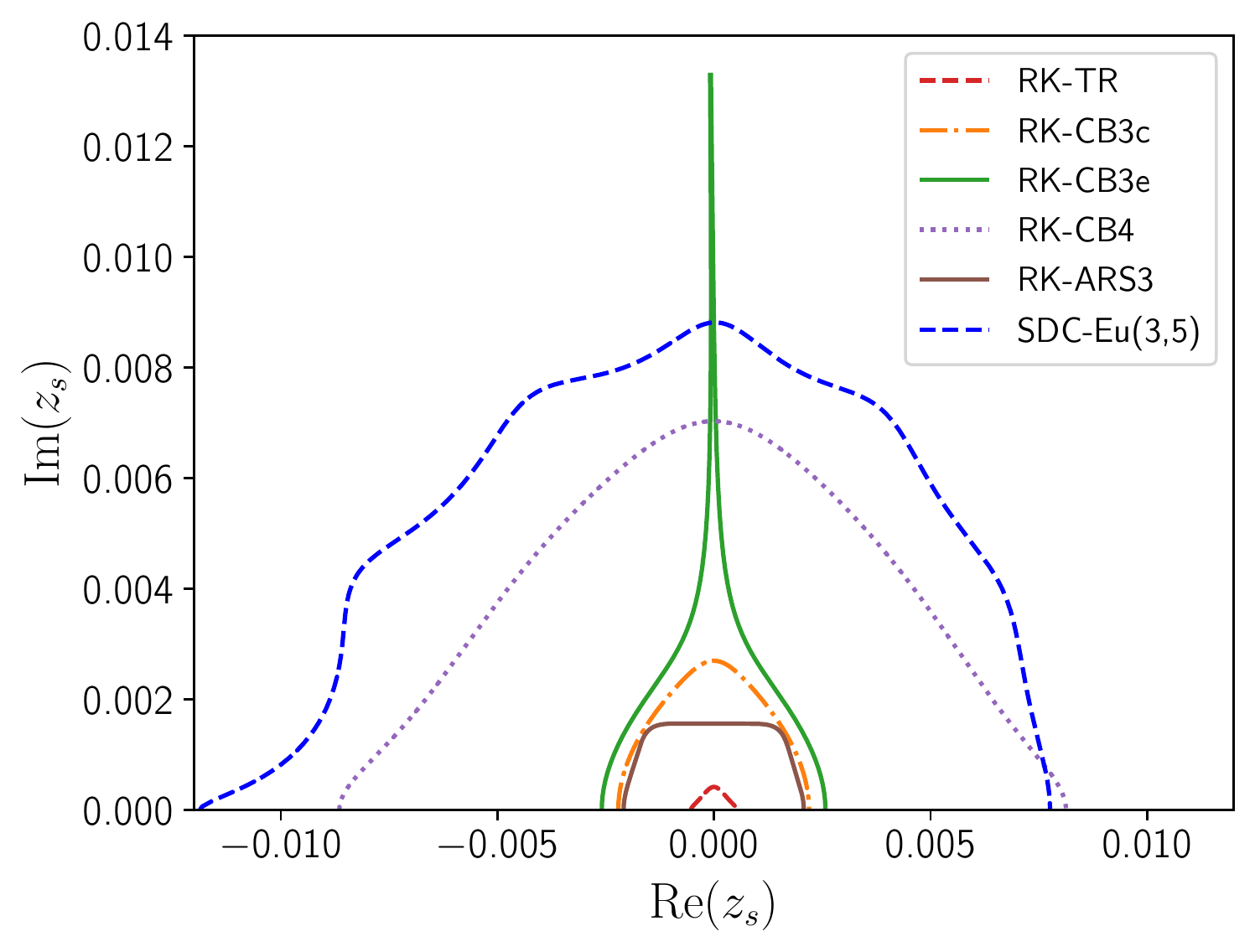}
   \label{fig:model:error-10}}
\caption{%
  Scaled accuracy domains of selected semi-implicit methods. 
  \label{fig:model:error}
  }
\end{figure}

As expected, the sixth-order SDC method admits the largest time steps over a wide 
range of $\real z$ (Fig.~\ref{fig:model:stability:unscaled}).
However, when applying the substep scaling (Fig.~\ref{fig:model:stability:scaled}), 
RK-CB3e allows scaled time steps $\Delta t_s$ roughly 
1.5 times larger than RK-ARS3 and 4 times larger than SDC-Eu(3,5).
SDC-CB3e(3,3) and SDC-ARS3(3,3) resemble the latter and are omitted for clarity. 
All remaining RK methods, including CB2, range between ARS3 and CB3e.
Figure~\ref{fig:model:error} compares the scaled accuracy domains of the same methods for two different thresholds.
For errors smaller than $10^{-5}$, RK-CB3e is clearly the most efficient method, followed by RK-CB4 and then RK-CB3c and SDC-Eu(3,5), with considerably smaller $z_s$ or, respectively, $\Delta t_s$.
With increasing accuracy, the 6-th order SDC methods gain efficiency and become 
competitive with RK-CB3e for ${|\varepsilon| \approx 10^{-10}}$ 
(Fig.~\ref{fig:model:error-10}).
In summary, its superior stability and accuracy render RK-CB3e the most promising method for the considered model problem.
However, this advantage is not necessarily preserved in the application to incompressible flow problems, where nonlinearity, boundary conditions, splitting errors and spatial discretization introduce additional challenges that can affect stability as well as accuracy.

\section{Application to incompressible Navier-Stokes problems}
\label{sec:navier-stokes}


\subsection{Governing equations}
\label{sec:navier-stokes:equations}

The time-integration methods introduced above are applied to incompressible flows with constant density and variable viscosity in a simply connected spatial domain ${\Omega \in \mathbb R^3}$.
The velocity ${\V v(\V x, t)}$ satisfies the momentum (Navier-Stokes) and continuity equations
\begin{gather}
  \label{eq:momentum}
  \d_t \V v + \nabla \cdot \V v \V v + \nabla p
  = \nabla \cdot \T\tau + \V f
  \,,
\\
  \label{eq:continuity}
  \nabla \cdot \V v = 0
\end{gather}
in $\Omega$, where 
$p$ 
is the pressure and
${\T\tau = \nu [\nabla \V v + \transpose{(\nabla \V v)}]}$
the viscous stress tensor, both divided by density; 
$\nu$ 
the kinematic viscosity, which may depend on $\V v$, and
${\V f(\V x, t)}$ a known forcing term.
The flow problem is closed by initial and boundary conditions
\begin{alignat}{2}
\label{eq:ic}
\V v(\V x, 0) &= \V v_0(\V x)
&\quad& \V x \in  \Omega
\,,
\\
\label{eq:bc}
\V v(\V x, t) &= \V v_{\mathrm b}(\V x, t)
&\quad& \V x \in \d\Omega
\,.
\end{alignat}
satisfying the continuity and compatibility constraints,
${\nabla \cdot \V v_0 = 0}$
and
${\int_{\d\Omega} \V n \cdot \V v_{\mathrm b} \D\Gamma = 0}$,
respectively.


\subsection{Time integration}
\label{sec:navier-stokes:ti}

The IMEX methods are adapted to incompressible Navier-Stokes problems by applying a projection-based splitting scheme in each step and with RK at each stage.
This scheme is based on the rotational velocity correction scheme developed in \cite{TI_Guermond2003b} and extended to variable viscosity in \cite{TI_Stiller2020a}.
For stating the methods, the flow equations are rewritten in the form  
\begin{equation}
  \d_t \V v = \V F(\V x, t, \V v, p),
  \quad
  \nabla \cdot \V v = 0
  \,,
\end{equation}
where 
${\V F = \V F_{\mathrm c}}+ \V F_{\mathrm d} + \V F_{\mathrm p} + \V F_{\mathrm s}$
with
${\V F_{\mathrm d}  = \V F_{\mathrm d1} + \V F_{\mathrm d2} + \V F_{\mathrm d3}}$, 
and
\begin{subequations}
\begin{alignat}{9}
& \V F_{\mathrm c}  && = &\,-& \nabla \cdot \V v \V v \,,  &\qquad
& \V F_{\mathrm p}  && = &\,-& \nabla p \,,                &\qquad
& \V F_{\mathrm s}  && = &   & \V f \,,                       \\
& \V F_{\mathrm d1} && = &   & \nabla \cdot \nu \nabla \V v \,,               &\qquad
& \V F_{\mathrm d2} && = &   & \nabla \cdot \nu \transpose{(\nabla \V v)} \,, &\qquad
& \V F_{\mathrm d3} && = &\,-& \nabla (\nu \nabla \cdot \V v) \,.
\end{alignat}
\end{subequations}
The additional viscous contribution $\V F_{d3}$ introduces a divergence penalty
which, for constant viscosity, yields the laplacian form of the viscous term, i.e.,
${\V F_{\mathrm d}  =  \nu \nabla^2 \V v}$.
Semi-discrete quantities write as follows
${\V v^n   \simeq \V v(\V x, t^n)}$ in general,
${\V v_i   \simeq \V v(\V x, t_i)}$ with RK at stage $i$ and
${\V v^k_m \simeq \V v^k(\V x, t_m)}$ with SDC at time $t_m$ after $k$ corrections.
Intermediate results are indicated by one or more primes, e.g. $\V v'$.


\subsubsection{IMEX\,BDF2}

With IMEX BDF2, the splitting scheme is applied once per time step, which yields four substeps: 
\begin{enumerate}

\item 
  Extrapolation
  \begin{equation}
    \label{eq:bdf2:extra}
    \frac{\gamma_0 \V v' - \alpha_0 \V v^n - \alpha_1 \V v^{n-1}}{\Delta t}
       = \beta_0 ( \V F_{\mathrm c}^{n}   
                 + \V F_{\mathrm d}^{n}
                 + \V F_{\mathrm d3}^{n} )
       + \beta_1 ( \V F_{\mathrm c}^{n-1}
                 + \V F_{\mathrm d}^{n-1}
                 + \V F_{\mathrm d3}^{n-1} )
       + \V F_{\mathrm s}^{n+1}
  \end{equation}

\item 
  Projection
  \begin{equation}
    \label{eq:bdf2:projection}
    \frac{\gamma_0 (\V v'' - \V v')}{\Delta t}
      = -\nabla p'', \quad
    \nabla \cdot \V v'' = 0, \quad
    \V n \cdot \V v''|_{\d\Omega} = \V n \cdot \V v_{\textrm b}(t_{n+1})
  \end{equation}

\item 
  Diffusion
  \begin{equation}
    \label{eq:bdf2:diffusion}
      \begin{aligned}
        \frac{\gamma_0 (\V v''' - \V v'')}{\Delta t} 
        & = \nabla \cdot [(\beta_0 \nu^{n} + \beta_1 \nu^{n-1}) \nabla \V v'''] \\
        & - \beta_0 ( \V F_{\mathrm d1}^{n}   + \V F_{\mathrm d3}^{n}   )   
          - \beta_1 ( \V F_{\mathrm d1}^{n-1} + \V F_{\mathrm d3}^{n-1} ),
        \quad
        \V v'''|_{\d\Omega} = \V v_{\textrm b}(t_{n+1})
      \end{aligned}
    \end{equation}

\item 
  Additional projection
  \begin{equation}
    \label{eq:bdf2:additional projection}
    \frac{\gamma_0 (\V v^{n+1} - \V v''')}{\Delta t}
      = -\nabla (p^{n+1} - p''), \quad
    \nabla \cdot \V v^{n+1} = 0, \quad
    \V n \cdot \V v^{n+1}|_{\d\Omega} = \V n \cdot \V v_{\textrm b}(t_{n+1})
  \end{equation}
 
\end{enumerate}
The extrapolation \eqref{eq:bdf2:extra} can be viewed as an explicit step that includes all terms except forcing, which is treated implicitly.
Adding $\V F_{\mathrm d3}$ yields an extra divergence penalty and transforms the diffusion term into rotational form, 
${\V F_{\mathrm d} + \V F_{\mathrm d3} = -\nu\nabla\times(\nabla\times\V v)}$,
if $\nu$ is constant.
The projection step \eqref{eq:bdf2:projection} requires the pressure $p''$, which is obtained by solving the related Poisson problem
\begin{equation}
  \label{eq:bdf2:pressure}
  \nabla^2 p'' = \frac{\gamma_0}{\Delta t} \nabla \cdot \V v',
  \quad
  \d_n p''|_{\d\Omega} = \V n \cdot
                         \frac{\gamma_0}{\Delta t}(\V v' - \V v_{\textrm b}^{n+1})     
  \,.
\end{equation}
Application of $p''$ in the first part of \eqref{eq:bdf2:projection} yields the divergence-free intermediate velocity $\V v''$.
The following step solves the implicit diffusion problem \eqref{eq:bdf2:diffusion}, which removes the explicit approximation of ${\nabla\cdot\nu\nabla\V v}$ and the extra divergence penalty introduced in \eqref{eq:bdf2:extra}.
If $\nu$ is constant, the diffusion step preserves continuity.
However, a variable viscosity can result in a non-solenoidal contribution to the RHS in \eqref{eq:bdf2:diffusion} and thus affect the continuity of the velocity field $\V v'''$.
Therefore, it is generally necessary to conclude the time step with an additional projection according to Eq.~\eqref{eq:bdf2:additional projection}.
If the viscosity is constant, the projection is skipped such that
${\V v^{n+1} = \V v'''}$.

\RII{%
While the splitting into substeps simplifies the computational procedure, it also introduces an additional temporal error, which reduces the convergence rate of the pressure to ${O(\Delta t^{3/2})}$ and leads to a divergence error of the same order \cite{TI_Guermond2003a}.
Especially the latter can affect the stability of the method and thus require smaller time steps than predicted by the theoretical analysis in Sec.~\ref{sec:time-integration:stability}.}


\subsubsection{IMEX\,RK}

Algorithm~\ref{alg:imex-rk:flow} outlines the IMEX Runge-Kutta method using a splitting similar to BDF2 at each stage except the first.
In comparison to Algorithm~\ref{alg:semi-rk:model}, $\V v''_i$ takes the role of $v_i$ and is used to compute the viscosity $\nu_i$.
The auxiliary quantity 
${\psi''_i \approx \Delta t \sum_{j=1}^{i-1} \aim_{i,j} p_j + \aim_{i,i} p''_i}$
combines the pressure contributions of different stages.
This simplifies the handling and removes the need for keeping the stage-level pressures.
Since no initial value is required, the computation of $p^{n+1}$ is dropped as well.
Line 11 of the algorithm presents an incremental form of the assembly which allows to omit terms with vanishing factors ${b^\ast_i - a^\ast_{s,i}}$.
The final projection (line 13--14) removes the residual velocity divergence in case of variable viscosity and is skipped with constant $\nu$.

\RII{%
It should be noted that the proposed algorithm introduces a splitting error, for which no theoretical analysis is available.
The stability and convergence properties of the resulting IMEX RK methods are therefore assessed by means of numerical experiments.
}

\begin{algorithm}[t]
\caption{IMEX\,Runge-Kutta method for incompressible flow with variable viscosity. 
Compact notation is applied were possible, e.g.\
${\V F_{\mathrm c,i} = \V F_{\mathrm c}(\V x, t + c_i \Delta t, \V v_{i})}$.}
\label{alg:imex-rk:flow}
\begin{algorithmic}[1]
\Procedure{ImexRK}{$\V v, t$}
   \State \makebox[1em][l]{$\V v_1$} $\gets \V v$
   \For{$i=2,s$}
      \State \makebox[1em][l]{$\V v'_i$} 
             ${\gets \V v 
               + \Delta t 
                   \big[ \sum_{j=1}^{i-1}
                           \aex_{i,j} ( \V F_{\mathrm c,j}
                                      + \V F_{\mathrm d,j}
                                      + \V F_{\mathrm d3,j} )
                       + \sum_{j=1}^{i} 
                           \aim_{i,j} \V F_{\mathrm s, j}
                     \big]
             }$
      \State \makebox[1em][l]{$\psi''_i$} 
             ${\gets  
                  \mathbf{solve}
                  \big[ 
                    \nabla^2 \psi''_i = \nabla\cdot\V v'_i, 
                    \;\;
                    \d_n \psi''_i|_{\d\Omega} = 
                      \V n\cdot (\V v'_i - \V v_{\textrm b,i})
                  \big]
             }$
      \State \makebox[1em][l]{$\V v''_i$} 
             ${\gets \V v'_i - \nabla \psi''_i}$
      \State \makebox[1em][l]{$\nu_i$} 
             ${\gets \nu(\V x, t+c_i\Delta t, \V v''_i)}$
      \State \makebox[1em][l]{$\V g_i$} 
             ${\gets 
                   \V v''_i 
                   + \Delta t 
                     \sum_{j=1}^{i-1}
                       \big[ (\aim_{i,j} - \aex_{i,j}) \V F_{\mathrm d1,j}
                           - \aex_{i,j}\V F_{\mathrm d3,j}
                       \big]
             }$
      \State \makebox[1em][l]{$\V v_i$} 
             ${\gets
                  \mathbf{solve}
                  \big[ 
                    \V v_i - \Delta t\,\aim_{i,i}\nabla\cdot(\nu_i \nabla\V v_i) 
                      = \V g_i, 
                    \;\;
                    \V v_i|_{\d\Omega} = \V v_{\textrm b,i}
                  \big]
             }$
   \EndFor
   \State \makebox[0.75em][l]{$\V{v}$} 
           ${\gets \V v_s
                 + \Delta t 
                       \sum_{i=1}^{s} 
                          \big[ (\bex_i - \aex_{s,i}) 
                                ( \V F_{\mathrm c,i} 
                                + \V F_{\mathrm d2,i}
                                + \V F_{\mathrm d3,i} )
                              + (\bim_i - \aim_{s,i})
                                ( \V F_{\mathrm d1,i}
                                + \V F_{\mathrm s,i} )
                        \big]
             }$

   \If{final projection}
      \State \makebox[0.75em][l]{$\psi$} 
             ${\gets  
                  \mathbf{solve}
                  \big[ 
                    \nabla^2 \psi = \nabla\cdot\V{v}, 
                    \;\;
                    \d_n \psi|_{\d\Omega} = 0
                  \big]
             }$
      \State \makebox[0.75em][l]{$\V v$} 
             ${\gets \V v - \nabla\psi}$
   \EndIf
\EndProcedure
\end{algorithmic}
\end{algorithm}


\subsubsection{SISDC}

The spectral deferred correction method considered in this study uses either IMEX Euler or a Runge-Kutta method as the predictor.
Each time step is divided into $M$ subintervals such that the intermediate times ${\{t_m\}}$ represent the GLL points, as described in Sec.~\ref{sec:time-integration:sdc}.
Starting from 
${\V v^0_0 = \V v(\V x, t^n)}$,
the predictor sweeps through the subintervals and computes
${\V v^0_m \approx \V v(\V x, t_m)}$
for ${m=1}$ to $M$.
This initial approximation is improved by the IMEX Euler corrector outlined in Algorithm~\ref{alg:imex-sdc:flow}.
It is a streamlined version of the most promising approach devised in
\cite{TI_Stiller2020a}
with simplified pressure handling analogous to IMEX RK.
Except for the extrapolation that contains additional SDC terms (line 17), a corrector step over one subinterval corresponds to an ordinary IMEX Euler step,
performing the projection on lines 18-19, followed by the diffusion step on line 20 and, in case of variable $\nu$, an additional projection on lines 22-23.
After $K$ sweeps, the corrector returns the solution and thus completes the time step.
\RII{%
Similar to IMEX BDF2 and RK, the SISDC method can be affected by splitting errors, which degrade the efficiency of the correction sweeps. 
Nevertheless, optimal accuracy can be recovered by increasing the number of sweeps \cite{TI_Stiller2020a}.
}

\begin{algorithm}[t]
\caption{SISDC corrector for incompressible flow with variable viscosity.}
\label{alg:imex-sdc:flow}
\begin{algorithmic}[1]
\Procedure{SISDC}{}($\NM{\V v, t}$)
   \State Initialize $\{t_m\}$ and $\{\Delta t_m\}$
   \State ${\NM{\V v}^0 \gets \NM{\V v}}$
   \For{$m=0,M$}
      \Comment{initial viscosity and RHS}
      \State $\nu^0_m \gets \nu(\V x, t_m, \V v^0_m)$
      \State\algorithmicif\ \makebox[3.25em][l]{$m < M$} \algorithmicthen\
             \makebox[2em][l]{$\V F^{\EX,0}_m$}
             ${\gets \V F_{\mathrm c} (\V v^0_m)
                   + \V F_{\mathrm d2}(\nu^0_m, \V v^0_m)
                   + \V F_{\mathrm d3}(\nu^0_m, \V v^0_m)}$
      \State\algorithmicif\ \makebox[3.25em][l]{$m > 0$} \algorithmicthen\
                \makebox[2em][l]{$\V F^{\IM,0}_m$}
                ${\gets \V F_{\mathrm d1}(\nu^0_{m-1}, \V v^0_m)}$
   \EndFor
   \For{$k=0,K-1$}
      \Comment{correction sweeps}

      \State \makebox[3em][l]{$\V v^{k+1}_0$}     ${\gets \V v^0_0}$
      \State \makebox[3em][l]{$\nu^{k+1}_0$}      ${\gets \nu^0_0}$
      \State \makebox[3em][l]{$\V F^{\EX,k+1}_0$} ${\gets \V F^{\EX,0}_0}$

      \For{$m=1,M$}
         \Comment{sweep through subintervals}

         \State \makebox[1.5em][l]{$\V F'_{\mathrm d1}$}
                ${\gets \V F_{\mathrm d1}(\nu^{k+1}_{m-1}, \V v^k_m)}$
         \State \makebox[1.5em][l]{$\V F'_{\mathrm d3}$}
                ${\gets \V F_{\mathrm d3}(\nu^{k+1}_{m-1}, \V v^k_m)}$
         \State \makebox[1.5em][l]{$\V S^k_m$} 
                ${\gets\sum_{q=0}^{M} w_{q,m} \V F(\V x, t_m, \V v^{k}_m,0)}$
         \State \makebox[1.5em][l]{$\V v'_m$}
               ${\gets \V v^{k+1}_{m-1}
                     + \Delta t_m \big[ 
                                    \V F^{\EX,k+1}_{m-1}
                                  + \V F'_{\mathrm d1}
                                  + \V F'_{\mathrm d3}
                                  - ( \V F^{\EX,k}_{m-1}
                                    + \V F^{\IM,k}_{m}
                                    )
                                  \big]
                     + \V S^k_m
                 }$

         \State \makebox[1.5em][l]{$\psi''_m$}
                ${\gets  
                    \mathbf{solve}
                    \big[ 
                      \nabla^2 \psi''_m = \nabla\cdot\V v'_m, 
                      \;\;
                      \d_n \psi''_m|_{\d\Omega} = 
                        \V n\cdot (\V v'_m - \V v_{\textrm b,m})
                    \big]
                }$

         \State \makebox[1.5em][l]{$\V v''_m$}
                ${\gets \V v'_m -  \nabla \psi''_m}$

         \State \makebox[1.5em][l]{$\V v'''_m$}
                ${\gets 
                    \mathbf{solve}
                    \big[
                      \V v'''_m 
                        - \Delta t_m \nabla\cdot\nu^{k+1}_{m-1} \nabla\V v'''_m
                        = \V v''_m
                        - \Delta t_m (\V F'_{\mathrm d1} + \V F'_{\mathrm d3}),
                      \,
                      \V v'''_m|_{\d\Omega} = \V v_{\textrm b,m}
                   \big]
                }$

         \If{final projection}

            \State \makebox[1.5em][l]{$\psi_m$}
                   ${\gets  
                       \mathbf{solve}
                       \big[ 
                         \nabla^2 \psi_m = \nabla\cdot\V v'''_m, 
                         \;\;
                         \d_n \psi_m|_{\d\Omega} = 0
                       \big]
                   }$

            \State \makebox[1.5em][l]{$\V v^{k+1}_m$}
                   ${\gets \V v'''_m - \nabla \psi_m}$

         \Else

            \State \makebox[1.5em][l]{$\V v^{k+1}_m$} ${\gets \V v'''_m }$

         \EndIf

         \State \makebox[3em][l]{$\nu^{k+1}_m$}
                ${\gets \nu(\V x, t_m, \V v^{k+1}_m)}$
        \Comment{update viscosity and RHS}
         \State \makebox[3em][l]{$\V F^{\EX,k+1}_m$}
                ${\gets \V F_{\mathrm c} (\V v^{k+1}_m)
                   + \V F_{\mathrm d2}(\nu^{k+1}_m, \V v^{k+1}_m)
                   + \V F_{\mathrm d3}(\nu^{k+1}_m, \V v^{k+1}_m)}$
         \State \makebox[3em][l]{$\V F^{\IM,k+1}_m$}
                ${\gets \V F_{\mathrm d1}(\nu^{k+1}_m, \V v^{k+1}_m)}$

      \EndFor
   \EndFor
   
   \State ${\NM{\V v} \gets \NM{\V v}^K}$

\EndProcedure
\end{algorithmic}
\end{algorithm}


\subsection{Spatial discretization and implementation}
\label{sec:navier-stokes:dg}

For spatial discretization, the domain is decomposed into hexahedral elements on which the solution is approximated with tensor-product Lagrange polynomials based on GLL points.
A polynomial degree of $P$ is assumed for the velocity and ${P-1}$ for the pressure.
The substeps of the time method are discretized using the discontinuous Galerkin spectral-element method (DG-SEM)
with consistent integration and local Lax-Friedrichs fluxes for convection, and
the interior penalty method for pressure and viscous diffusion equations.
For pressure robustness a divergence/mass-flux stabilization is added as proposed in \cite{SE_Joshi2016a,SE_Akbas2018a}.
A detailed description of the DG-SEM is given in \cite{TI_Stiller2020a}.
The method was implemented in the HiSPEET%
\footnote{Available at \url{https://gitlab.hrz.tu-chemnitz.de/hispeet/hispeet.git}.} 
library which makes use of MPI for parallelization and LIBXSMM \cite{Heinecke2016a} for vectorizing the element operators.
It also provides highly efficient Krylov-accelerated Schwarz/multigrid methods, which are used for solving the implicit pressure and diffusion problems
\cite{SE_Stiller2016a,SE_Stiller2016b,SE_Stiller2017a}.


\section{Numerical experiments}
\label{sec:numerical-experiments}

A series of numerical experiments is conducted to assess the accuracy, convergence properties and computational efficiency of the considered IMEX methods. 
The investigated flow problems include: 
1) a traveling 2D Taylor-Green vortex, 
2) a 3D vortex flow with variable viscosity and
3) the nonlinear instability and turbulent decay of a disturbed Taylor-Green vortex.
All problems are periodic, which is exploited to set up a periodic test cases.
For problem 1 and 2 the exact solution is known and used to construct further test cases that include time-dependent boundary conditions.
Thus, the test cases allow a separate investigation of the order reduction due to boundary conditions or variable viscosity, respectively.

The study encompasses all time integration methods introduced in section \ref{sec:navier-stokes}.
For affordability, the SISDC method is confined to Euler, RK-CB3e and RK-ARS3 predictors.
Using ${M=3}$ subintervals and ${K=5}$ corrections with Euler and ${K=3}$ with RK 
yields a theoretical order of 6.
The resulting methods are denoted as 
SDC-Eu(3,5),
SDC-CB3e(3,3) and
SDC-ARS3(3,3), respectively.
In all simulations, the pressure and diffusion problems are solved with a relative tolerance of 10\textsuperscript{\textminus 12} of the RMS residual evaluated at the collocation points.

\JS{%
For comparison with the model problem studied in Sec.~\ref{sec:time-integration:stability}, the CFL number is defined as 

\begin{equation}
  \Clmb = \Delta t v_{\mathrm{ref}} |\lambda_{\mathrm c}|_{\max}
  \,,
\end{equation}
where $v_{\mathrm{ref}}$ is the reference velocity magnitude and
$\lambda_{\mathrm c}$  are the eigenvalues of the one-dimensional element 
convection operator for unit velocity.
In the present case, the eigenvalues result from the GLL collocation differentiation
operator of degree $P$ combined with one-sided Dirichlet conditions, 
see e.g. \textcite[Sec. 7.3.3]{SE_Canuto2011a}.
}


\subsection{Traveling Taylor-Green vortex}
\label{sec:numerical-experiments:TG}

\begin{figure}
\subfloat[\RI{Case 1:} periodic conditions in $y$-direction]
  {\begin{minipage}{0.56\textwidth}
   \raggedleft
   \includegraphics[scale=0.60]{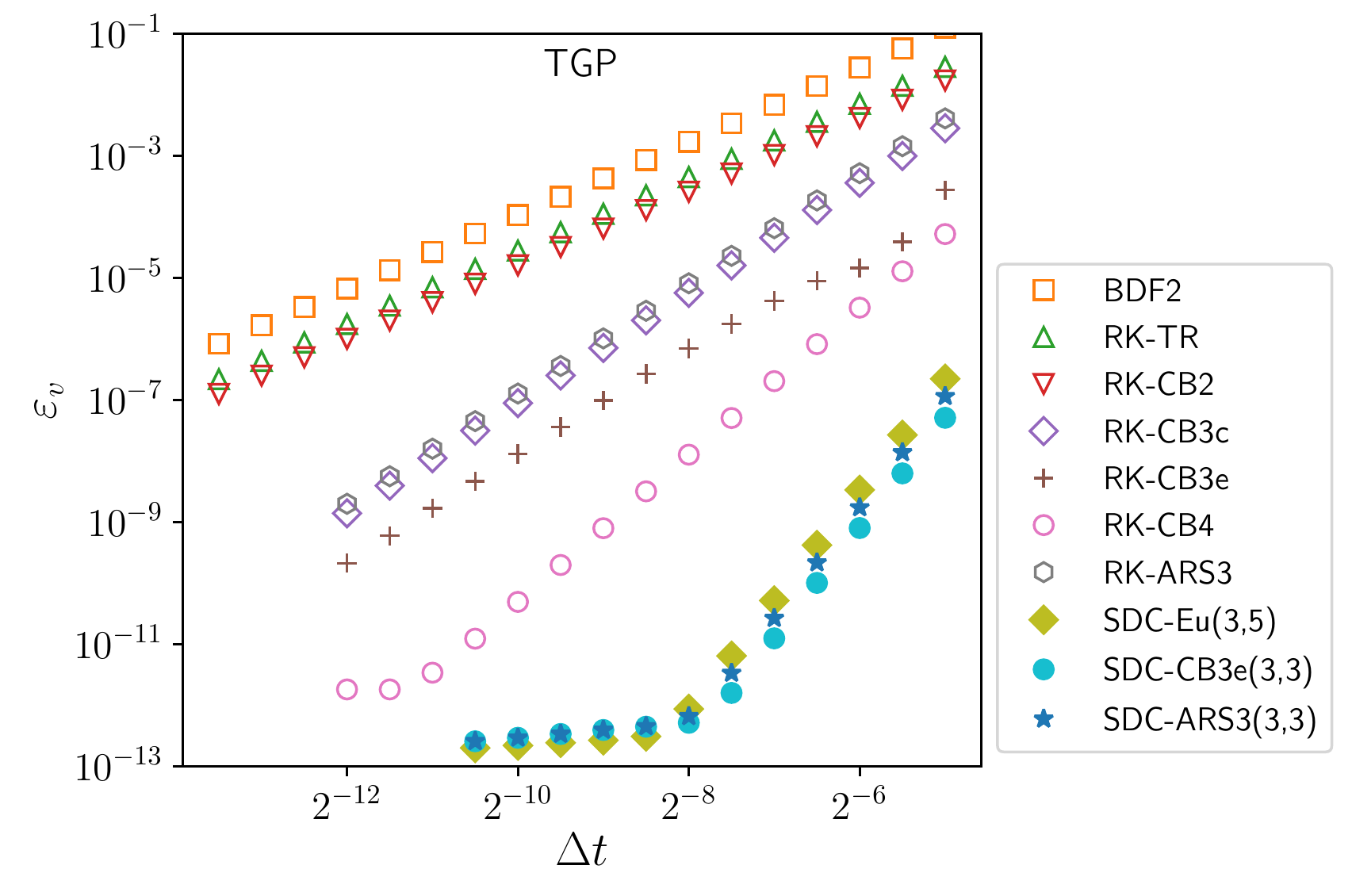} \\
   \includegraphics[scale=0.60]{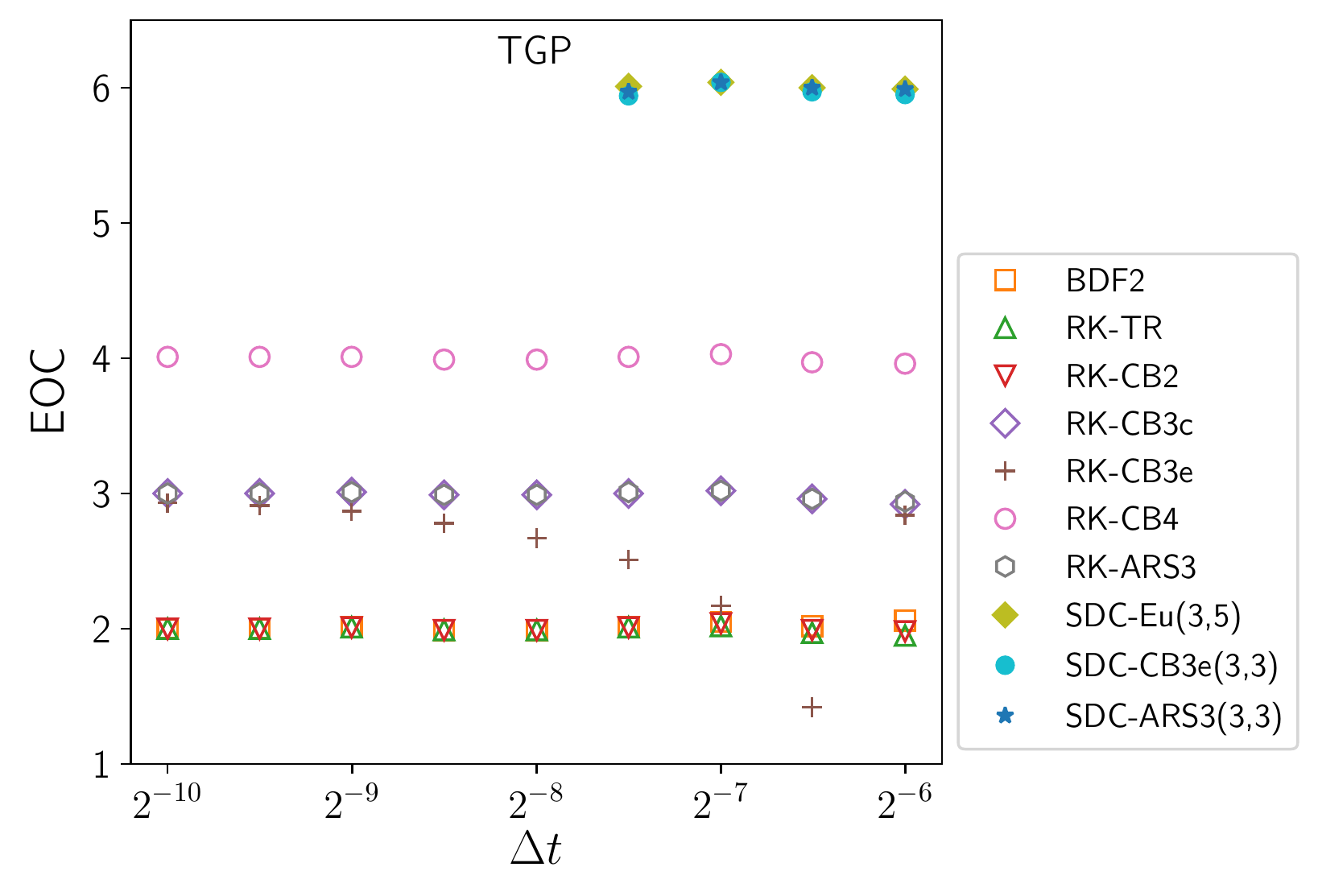}   \\
   \includegraphics[scale=0.60]{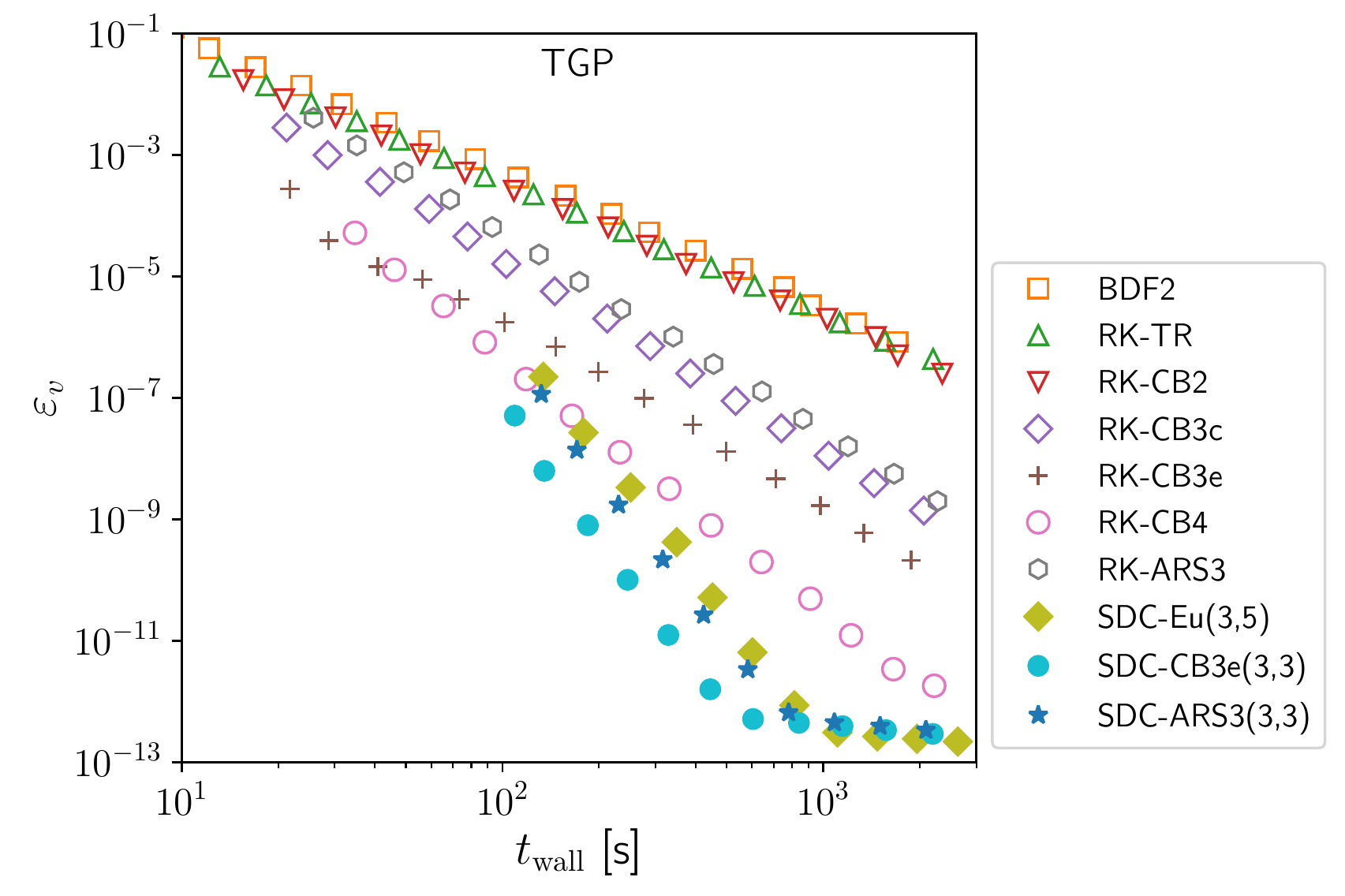}
   \end{minipage}
   \label{fig:TGP}}
   \hspace*{-16mm}
\subfloat[\RI{Case 2:} Dirichlet conditions in $y$-direction\hspace*{-22mm}]
  {\begin{minipage}{0.52\textwidth}
   \raggedleft
   \includegraphics[scale=0.60]{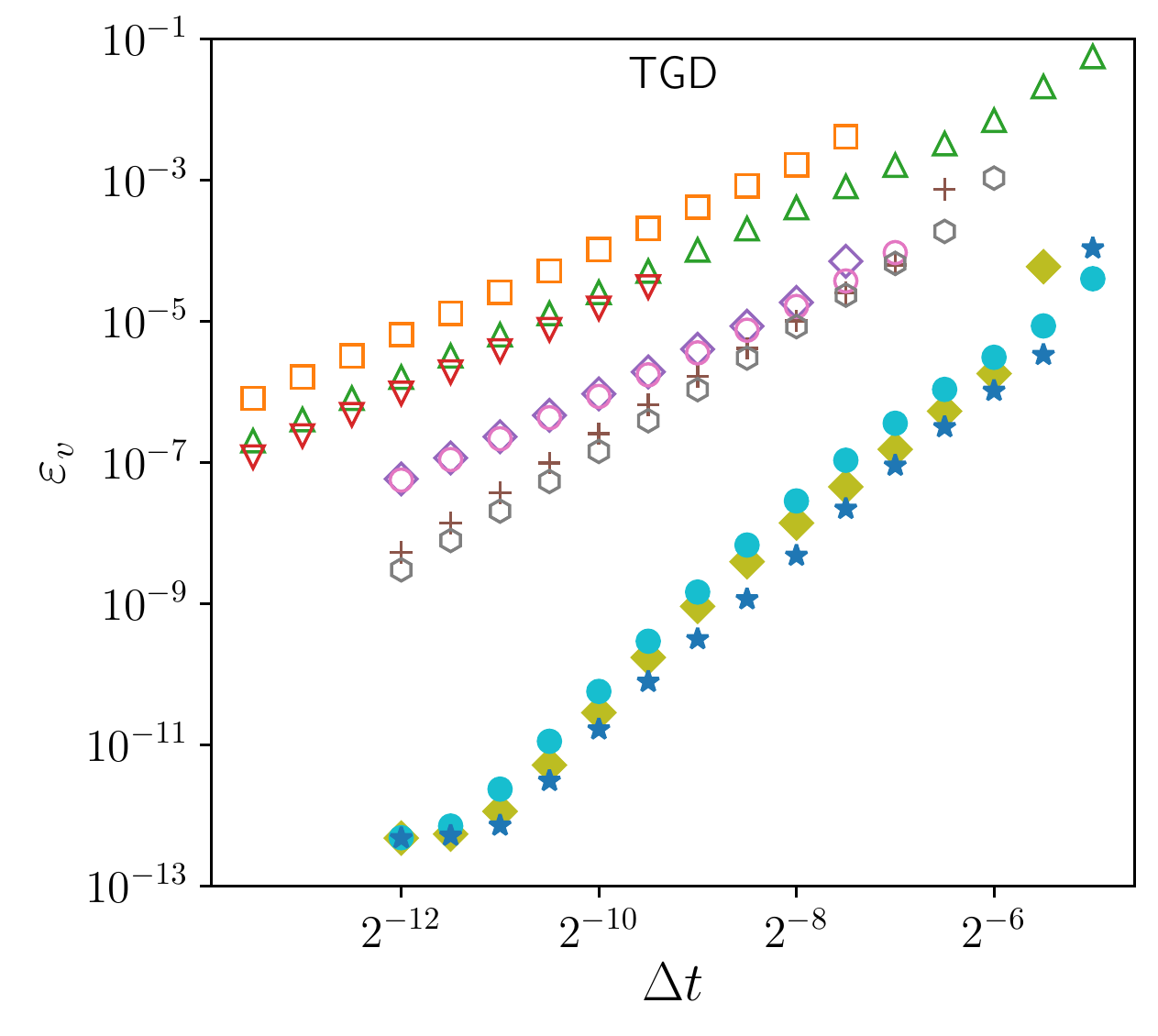} \\
   \includegraphics[scale=0.60]{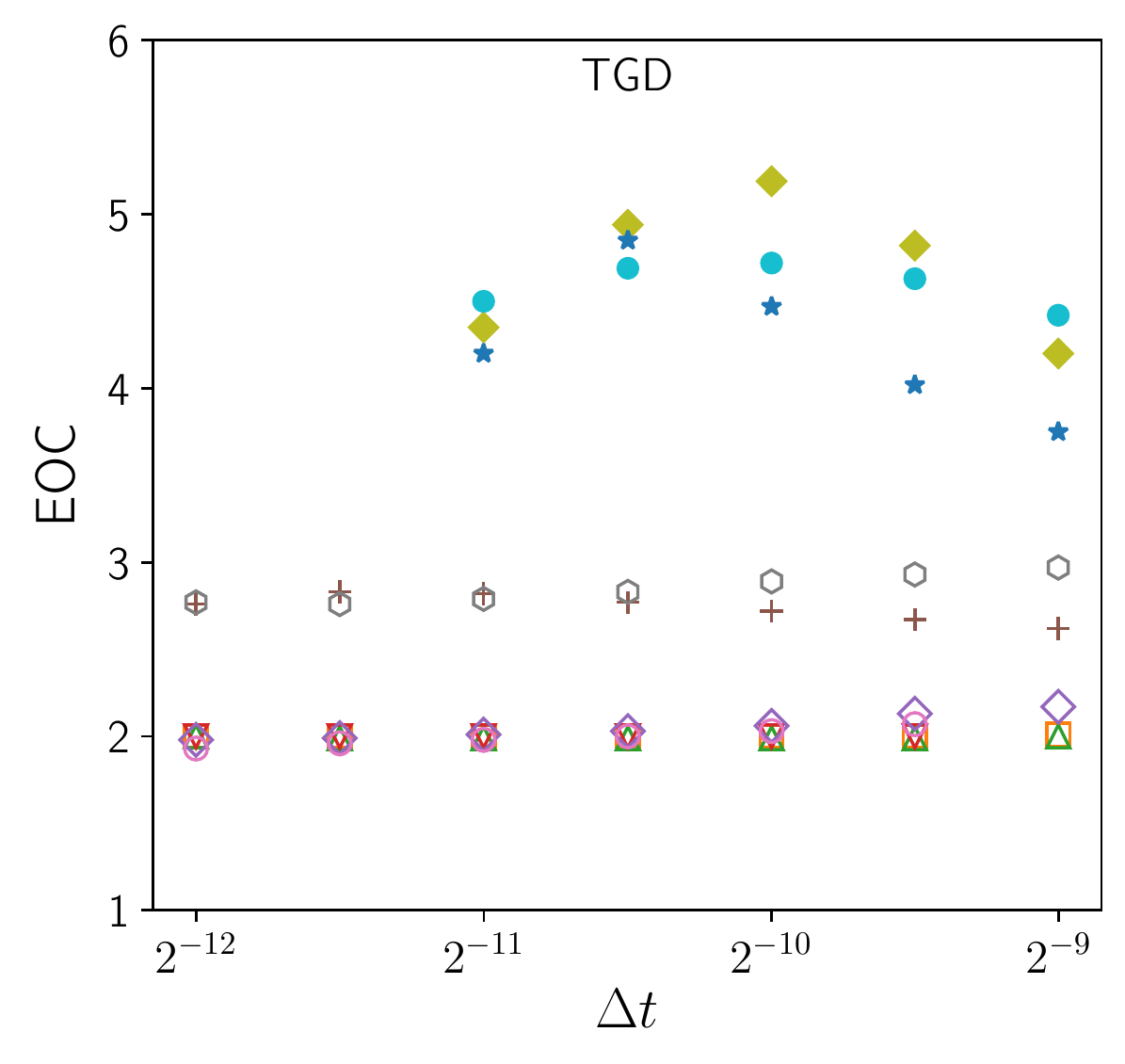}   \\
   \includegraphics[scale=0.60]{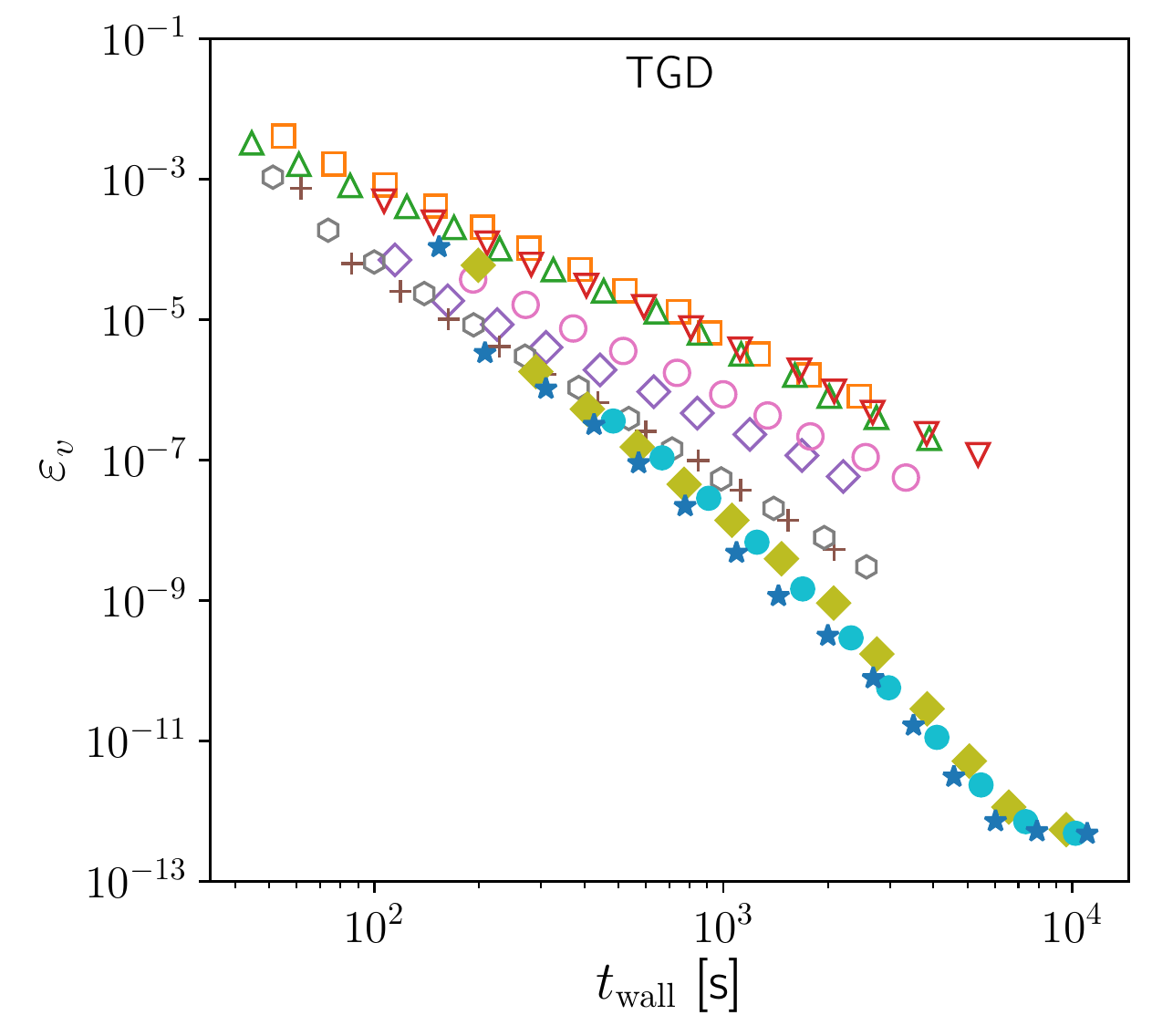}
   \end{minipage}
   \label{fig:TGD}}
\caption{\RI{%
  Results for the Taylor-Green test cases obtained with different time integration methods. 
  The top row shows the RMS velocity error at the final time,  $\varepsilon_v$,
  the middle row the experimental order of convergence, $\mathrm{EOC}$, for subsequent stepsizes and the bottom row the error obtained with a given wall-clock time, $t_{\mathrm{wall}}$.}
  \label{fig:TG}
  }
\end{figure}

The first problem represents a 2D Taylor-Green vortex that travels with a phase speed of one in the $x$ and $y$ directions. 
It was proposed by \textcite{TI_Minion2018a} and possesses the exact solution
\begin{subequations}
\begin{alignat}{2}
  & v_x^{\mathrm{ex}}(x,y,t) 
  &&= 1 + \sin \big( 2\pi ( x                - t )\big) 
          \cos \big( 2\pi ( y - \tfrac{1}{8} - t )\big) 
          \exp( -8\pi^2\nu t ) 
  \,,
  \\
  & v_y^{\mathrm{ex}}(x,y,t) 
  &&= 1 - \cos \big( 2\pi ( x                - t )\big)
          \sin \big( 2\pi ( y - \tfrac{1}{8} - t )\big)
          \exp( -8\pi^2\nu t ) 
  \,,
  \\
  & p^{\mathrm{ex}}(x,y,t)   
  &&=     \tfrac{1}{4}\big[ \cos \big(4\pi (x                - t)\big)
                          + \cos \big(4\pi (y - \tfrac{1}{8} - t)\big) 
                          \big]
          \exp(-16\pi^2\nu t)
  \,.
\end{alignat}
\end{subequations}
There is no forcing and no flow in $z$ direction.
The computational domain 
${\Omega = [-\sfrac{1}{2},\sfrac{1}{2}]^2 \times [-\sfrac{1}{16},\sfrac{1}{16}]}$
is decomposed into $8\times8\times1$ cubic elements of degree ${P=10}$.
Following \cite{TI_Minion2018a} two test cases are defined:
TGP with all directions periodic and kinematic viscosity ${\nu = 0.02}$, and
TGD with Dirichlet conditions in $y$ direction, $x$, $z$ periodic and ${\nu = 0.01}$.
Both cases are initialized with the exact solution and integrated 
with steps ranging from ${\Delta t = 2^{-5}}$ to $2^{-13}$
until reaching the final time of ${T=0.25}$.
The computations are run on 4 cores of an Intel Xeon E5-2680 v3 CPU.
Figure~\ref{fig:TG} presents 
the measured RMS error of velocity at final time $\varepsilon_v$,
the experimental order of convergence ($\mathrm{EOC}$) for successive stepsizes and
the wall-clock time for a given error $t_{\mathrm{wall}}$
for both cases.

In case TGP (Fig.~\ref{fig:TGP}), 
all investigated methods are stable and show a regular convergence behavior, starting from the largest step size, ${\Delta t = 2^{-5}}$.
The $\mathrm{EOC}$ consistently corresponds to the expected rate of convergence with the exception of RK-CB3e, which starts at a lower rate with large steps before reaching the theoretical order with ${\Delta t < 2^{-9}}$.
Nevertheless, RK-CB3e achieves the smallest error of all third-order methods for a given step size.
Among the second-order methods, RK-CB2 is the most accurate, followed by RK-TR, which is still almost twice as accurate as BDF2.
For a given error, this results in runtime savings of ca 75\% with RK-CB2 and 60\% with RK-TR.
Starting from 
${\varepsilon_v \approx 10^{-3}}$,
the RK methods of order ${\ge 3}$ become much more efficient than 
all second-order methods, but are surpassed by sixth-order SDC methods for 
${\varepsilon_v \lesssim 10^{-7}}$.
Its superior accuracy gives RK-CB3 and, hence, SDC-CB3e(3,3) an advantage over competitors and makes them the best methods of order 3 and 6, respectively.

The regular behavior observed in case TGP changes considerably with the imposition of time-dependent boundary conditions in case TGD (Fig.~\ref{fig:TGD}).
Only three of 10 methods succeed with the largest time step, which corresponds to ${\Clmb \approx 11.2}$.
RK-CB2 is the least stable method:
It requires ${\mathit{CFL}_\lambda \lesssim 2}$ which is less than half the critical value of $5.81$ obtained for the model problem. 
All, but the second-order methods are affected by order reduction:
RK-CB3c and RK-CB4 degrade to order two, while
RK-CB3e and RK-ARS3 suffer only a mild reduction to approximately $2.8$.
The $\mathrm{EOC}$ of SDC methods starts around $4$ for the largest time steps and grows up to ca $5$ with decreasing step size, before it stagnates when approaching the minimal error. 
A similar behavior was observed in \cite{TI_Stiller2020a}, where the theoretical order could be recovered by increasing the number of corrections.
However, this option is not considered here because of its higher cost. 
In terms of computational efficiency RK-TR supersedes RK-CB2 as the most efficient second-order method, although the savings are reduced to ca 20\% compared to BDF2.
The third-order Runge-Kutta methods, RK-CB3e and RK-ARS3, gain advantage for errors  smaller then $10^{-3}$, but are outperformed by the SDC methods led by SDC-ARS3(3,3) for 
${\varepsilon_v \lesssim 10^{-6}}$.
%


\subsection{3D Vortex with solution-dependent variable viscosity}
\label{sec:numerical-experiments:VV}

The second flow problem is based on a manufactured three-dimensional solution of the Navier-Stokes equations.
It was proposed by \textcite{Niemann2018a} as a test bench for spatially varying viscosity and generalized to \JS{nonlinear,} velocity-dependent viscosity in \cite{TI_Stiller2020a}.
The exact solution is given by
\begin{subequations}
\label{eq:var-visc-test}
\begin{alignat}{2}
  \label{eq:var-visc-test:vx}
  & v_x^{\mathrm{ex}}
  &&= \big[ \sin( 2\pi(x + t)) + \cos( 2\pi(y + t)) \big] \sin( 2\pi(z + t)) 
  \,,
  \\
  \label{eq:var-visc-test:vy}
  & v_y^{\mathrm{ex}}
  &&= \big[ \cos( 2\pi(x + t)) + \sin( 2\pi(y + t)) \big] \sin( 2\pi(z + t))
  \,,   
  \\
  \label{eq:var-visc-test:vz}
  & v_z^{\mathrm{ex}} 
  &&= \big[ \cos( 2\pi(x + t)) + \cos( 2\pi(y + t)) \big] \cos( 2\pi(z + t))
  \,,
  \\
  \label{eq:var-visc-test:p}
  & p^{\mathrm{ex}}   
  &&= \sin( 2\pi(x + t)) \sin( 2\pi(y + t)) \sin( 2\pi(z + t))
  \,.
\end{alignat}
\end{subequations}
It defines a periodic vortex array with wave length ${l=1}$ and velocity magnitude ${v^{\mathrm{ex}}_{\max} = 2}$, which travels with a phase velocity of 1 in each direction, separately.
The kinematic viscosity is designed to mimic a simple model of unresolved turbulence and takes the form
\begin{equation}
  \nu(\V v) = \nu_0 + \nu_1 \left(\frac{|\V v|}{v^{\mathrm{ex}}_{\max}}\right)^2
\end{equation}
\JS{with constant coefficients $\nu_0$ and $\nu_1$.} 
Rearranging the momentum equations applied to the exact solution yields the
forcing
\begin{equation}
  \label{eq:var-visc-test:rhs}
  \V f(\V x, t)
       = \d_t \V v^{\mathrm{ex}}
       + \nabla \cdot \V v^{\mathrm{ex}} \V v^{\mathrm{ex}}
       - \nabla \cdot \Big[ \nu(\V v^{\mathrm{ex}})
                            \Big( \nabla \V v^{\mathrm{ex}}
                                + \transpose{(\nabla \V v^{\mathrm{ex}})}
                            \Big)
                       \Big]
       + \nabla p^{\mathrm{ex}}
  \, .
\end{equation}
\JS{%
The Reynolds number is defined using the mean viscosity based on $\V v^{\mathrm{ex}}$, i.e.,
${\mathit{Re} = 2 / \bar\nu}$
with
${\bar\nu = \nu_0 + \frac{1}{2}\nu_1}$.

To assess the properties of the time integration method, three different scenarios are considered:
1) convergence properties and efficiency with periodic and Dirichlet conditions at a moderate Reynolds number,
2) impact of viscosity and boundary conditions on stability and accuracy,
3) stability in the convection dominant case.
The computational domain ${\Omega = [ -\sfrac{1}{2}, \sfrac{1}{2} ]^3}$
extends over one period length in every direction,
which allows to define two test cases:
VVP with periodic conditions and
VVD with time-dependent Dirichlet conditions such that
${\V v_{\mathrm b} = \V v^{\mathrm{ex}}}$ 
on $\d\Omega$.
}


\subsubsection{Convergence properties and efficiency at moderate Reynolds number}
\label{sec:numerical-experiments:VV:convergence}

\begin{figure}
\subfloat[\RI{Case VVP:} Periodic conditions]
  {\begin{minipage}{0.57\textwidth}
   \raggedleft
   \includegraphics[scale=0.60]{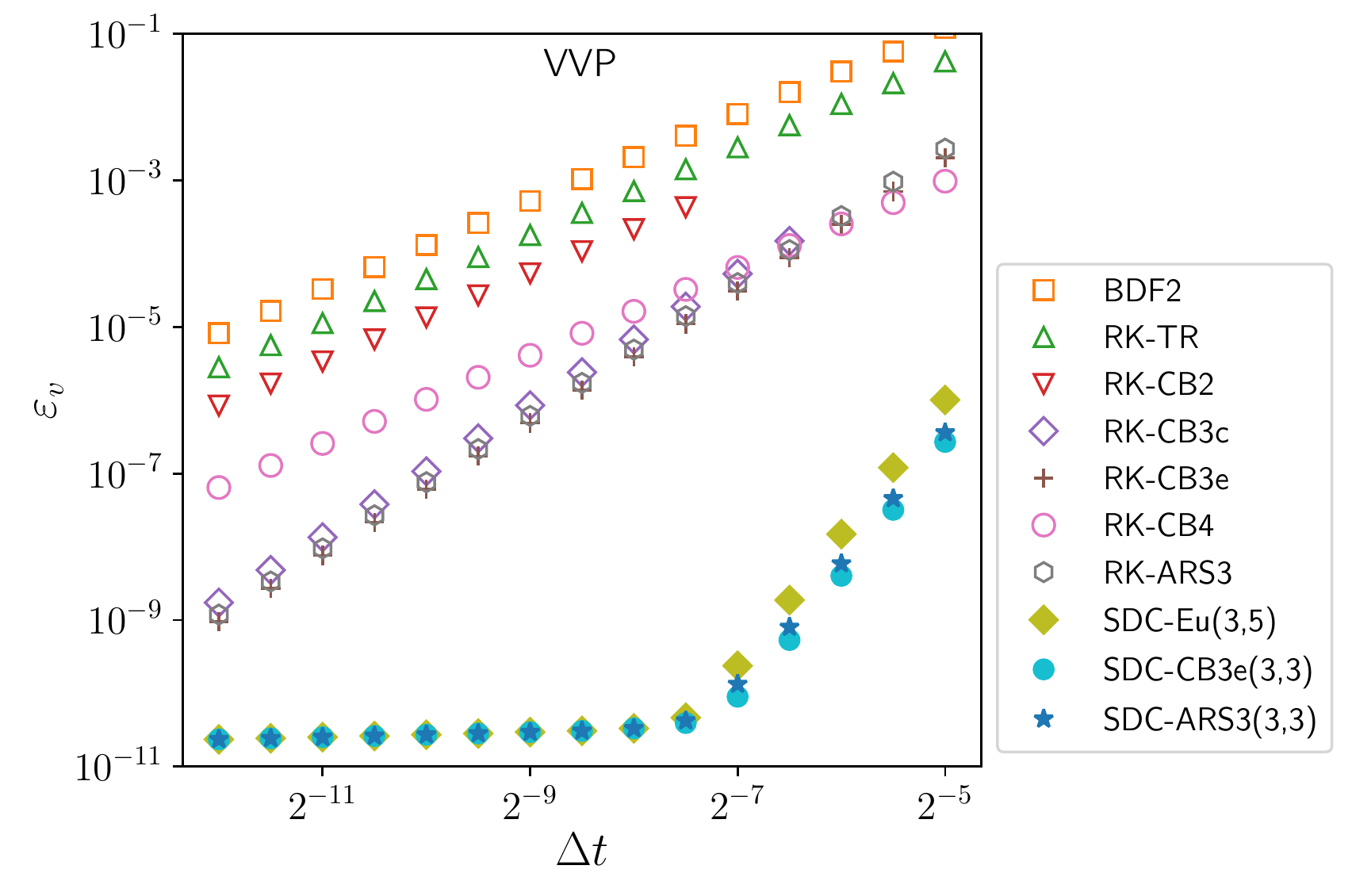} \\
   \includegraphics[scale=0.60]{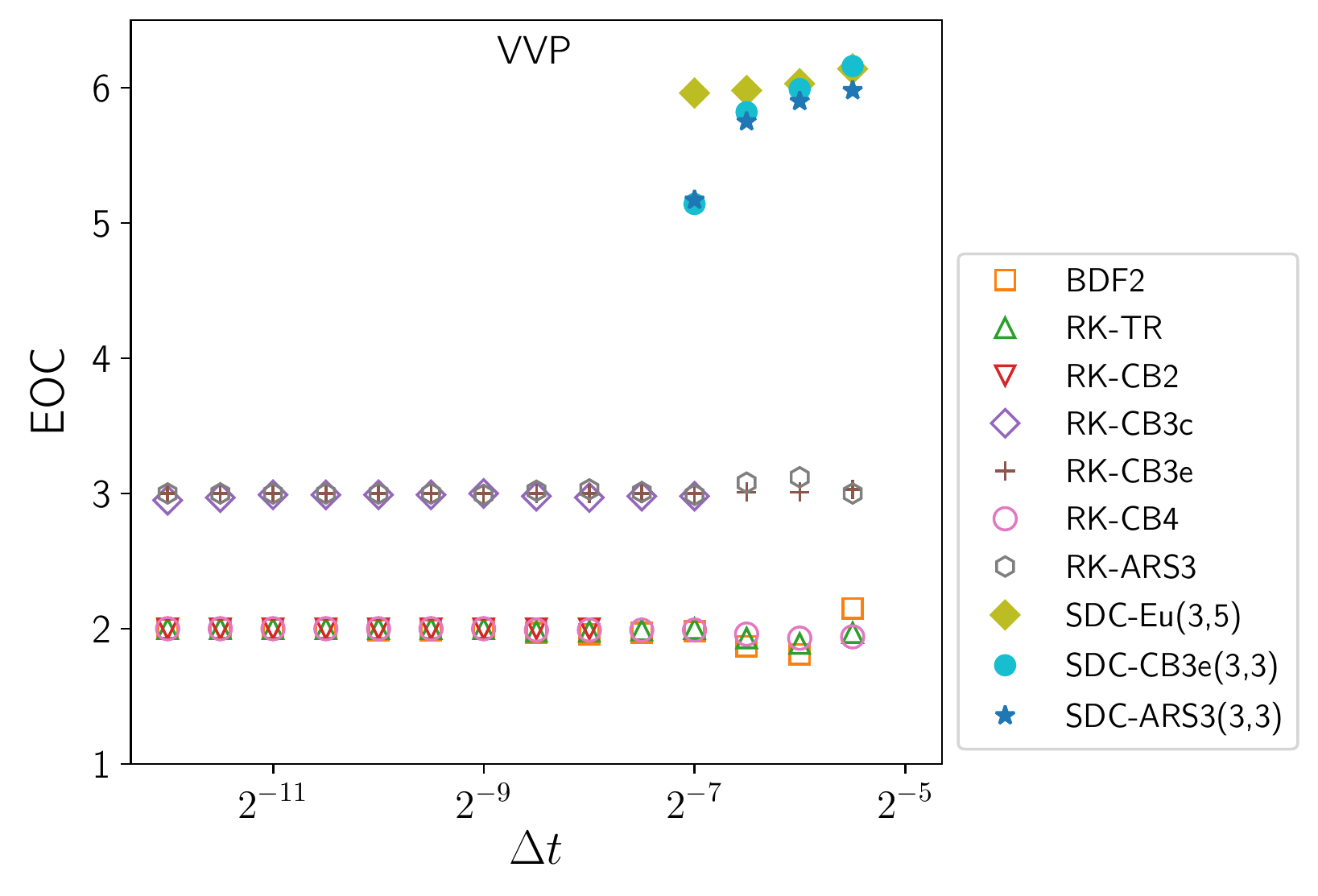}   \\
   \includegraphics[scale=0.60]{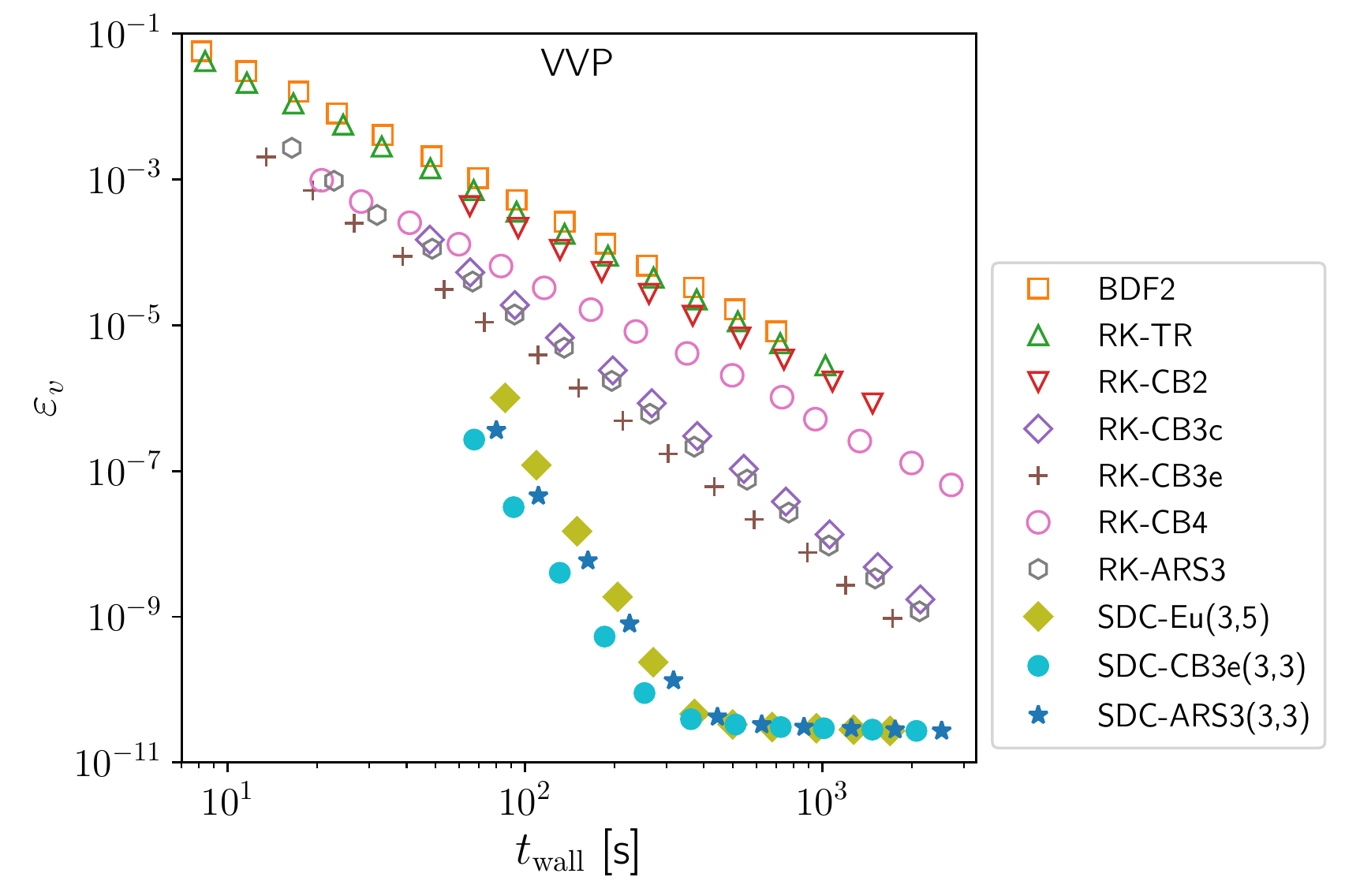}
   \end{minipage}
   \label{fig:VVP}}
\subfloat[\RI{Case VVD:} Dirichlet conditions]
  {\begin{minipage}{0.43\textwidth}
   \raggedleft
   \includegraphics[scale=0.60]{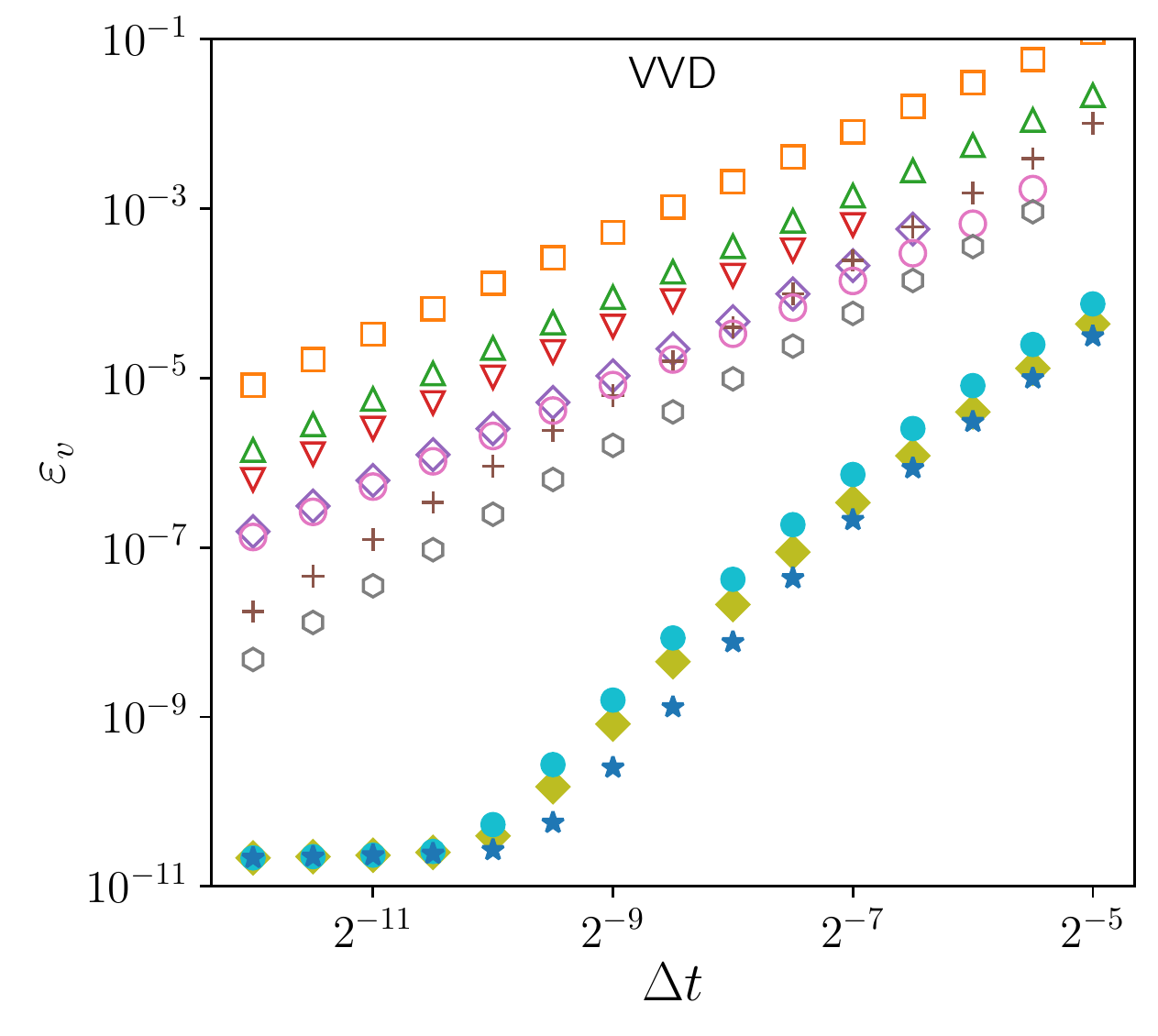} \\
   \includegraphics[scale=0.60]{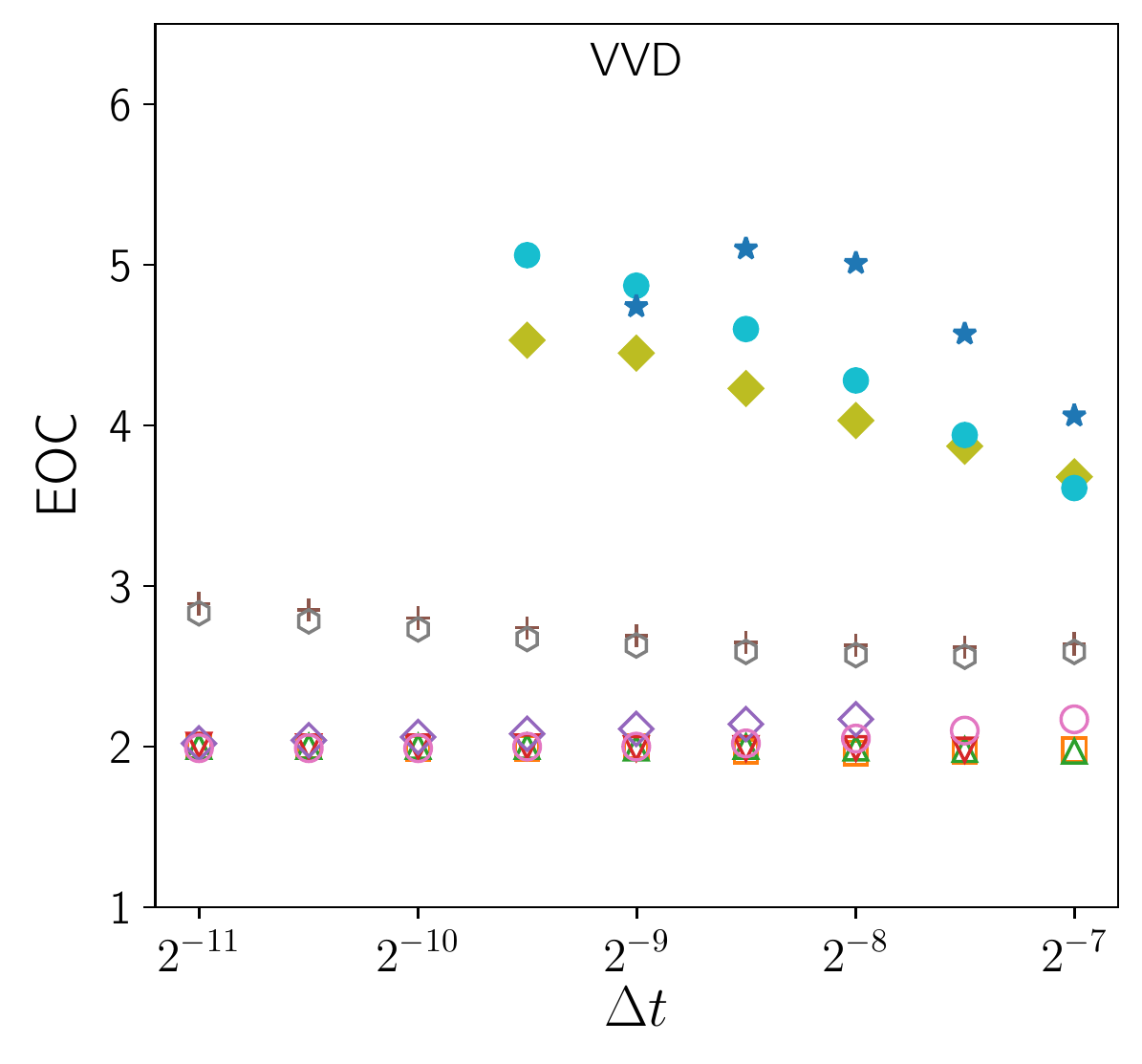}   \\
   \includegraphics[scale=0.60]{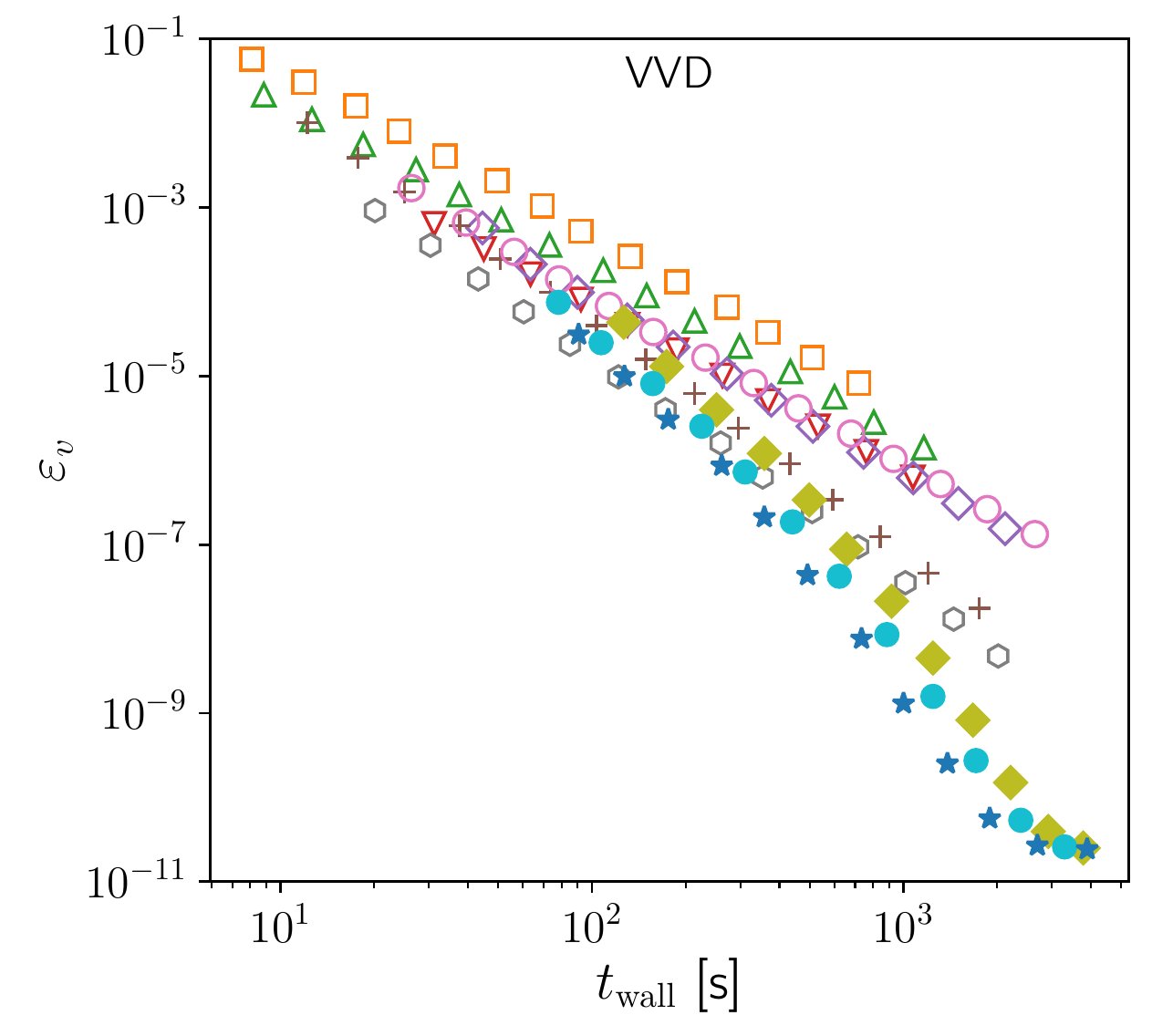}
   \end{minipage}
   \label{fig:VVD}}
\caption{\RI{%
  Results for the 3D vortex array test cases VVP and VVD using a solution-dependent viscosity 
  with ${\nu_0 = \nu_1 = 0.01}$.
  The top row shows the RMS velocity error at the final time,  $\varepsilon_v$,
  the middle row the experimental order of convergence, $\mathrm{EOC}$, for subsequent stepsizes and the bottom row the error obtained with a given wall-clock time, $t_{\mathrm{wall}}$.}
  \label{fig:VV}
  }
\end{figure}

\JS{%
For this scenario, the viscosity coefficients are set to 
${\nu_0 = \nu_1 = 0.01}$, which yields a Reynolds number of
${\mathit{Re} \approx 133}$.} 
The computational domain is decomposed into $2^3$ elements of degree ${P=16}$ and time integration is performed with steps ranging from ${\Delta t = 2^{-5}}$ to $2^{-12}$ until the final time of ${T=0.25}$ is reached.
All simulations are run on 8 cores of an Intel Xeon E5-2680 v3.
Figure~\ref{fig:VV} shows
the error at the final time $\varepsilon_v$,
the experimental order of convergence $\mathrm{EOC}$, and
the wall-clock time for a given error $t_{\mathrm{wall}}$
for both cases.

Under periodic conditions, imposed in case VVP, a similar picture emerges as before in the TGP case: 
all time integration methods achieve their theoretical order except for RK-CB4, which degrades to order 2.
RK-CB2 remains the most accurate second-order method, but fails with
${\Delta t \gtrsim 2^{-7}}$, as does RK-CB3c with $\Delta t \gtrsim 2^{-6}$.
RK-CB3e and RK-ARS3 exhibit nearly the same accuracy, but the former
is more efficient due to the lower number of stages.
The SDC methods outperform all lower order methods in the examined range, led by SDC-CB3e(3,3).
However, extrapolation indicates a break-even between RK-CB3e and SDC-CB3e(3,3) near 
${\varepsilon_v \approx 10^{-4}}$.

The imposition of Dirichlet conditions in case VVD yields an order reduction almost identical to TGD.
This leads to the conclusion that the semi-implicit handling of the nonlinear viscosity does not cause any further degradation.
Similar to case TGD, RK-ARS3 supersedes RK-CB3e as the computationally most efficient method for errors ${\varepsilon_v \gtrsim 10^{-5}}$.
For ${\varepsilon_v = 10^{-3}}$ it saves about 70\% of the time required with BDF2.
A further reduction of the error yields an increasing advantage of SDC methods, in particular SDC-ARS3(3,3).

\RII{ 

\subsubsection{Impact of viscosity and boundary conditions}
\label{sec:numerical-experiments:VV:robustness}

To evaluate the influence of the magnitude and nonlinearity of the viscosity and the imposition of time-dependent Dirichlet boundary conditions, respectively, an additional study is performed in which these parameters are changed independently.
As an example, Figure~\ref{fig:robustness:rk_ars3} shows the behavior of RK-ARS3 for two different magnitudes of viscosity with periodic and Dirichlet conditions, respectively.
To identify the impact of nonlinearity, the viscosity is evaluated using either the computed approximate velocity $\bm{v}$ or the exact solution $\bm{v}^{\mathrm{ex}}$.
The spatial discretization and the final time $T$ are identical to the previous study.
As expected, the fasted convergence and highest stability are observed with moderate viscosity, and periodic boundaries, i.e., 
${\nu_0 = \nu_1 = 0.01}$, 
which corresponds to 
${\mathit{Re} \approx 133}$.
Under these conditions, the theoretical order of $3$ is reproduced and the computations remain stable with time steps up to ${\Delta t_{\max} = 0.05}$, with practically no difference between linear and nonlinear viscosity.
Switching from periodic to Dirichlet conditions leads to a slight reduction in the experimental convergence rate to $2.93$, but has no impact on stability. 
Increasing the viscosity to ${\nu_0 = \nu_1 = 1}$, i.e. ${\mathit{Re} \approx 1.33}$, results in higher errors and affects stability in the case of nonlinear viscosity $\nu(\bm{v})$.
Using the linear viscosity $\nu(\bm{v}^{\mathrm{ex}})$ and periodic conditions, convergence begins at a rate of about $2$ at larger time steps and attains $2.99$ at smaller ones.
Additionally imposed Dirichlet conditions lead to a further reduction of the $\mathrm{EOC}$,
which now reaches only about $2.4$ at smaller time steps.
Regardless of the boundary conditions, the semi-implicit treatment of nonlinear viscosity results in severe restrictions of the time step, which indicates that a fully implicit treatment is advisable in the diffusion dominant case.
On the other hand, the results indicate that the proposed semi-implicit approach does not have a negative impact on accuracy.
Hence, the observed order reduction is likely caused by the stiffness, growing with viscosity, and the splitting used in the RK stages.
Based on the current results it seems impossible to identify the relative importance and possible interaction of these two factors. 
For an in-depth analysis, one could eliminate the splitting error with a coupled method, but this would go beyond the scope of the present work.
Of the other high-order IMEX RK methods, only RK-CB3e remains competitive to RK-ARS3, 
whereas RK-CB3c and RK-CB4 are more severely affected. 
For this reason and because of their poor performance in the previous tests, 
these methods are not included in the following investigations.

Finally, the properties of the SISDC method are studied in the diffusion dominant case with Dirichlet conditions.
Figure~\ref{fig:robustness:sdc_eu} shows the results obtained using an Euler-based method
with ${M=3}$ subintervals and a varying number $K$ of correction sweeps for a nonlinear 
viscosity with coefficients ${\nu_0 = \nu_1 = 1}$.
As expected, increasing the number of sweeps $K$ improves accuracy and stability.
However, for given $K$, the method shows an order reduction similar to the one observed with RK methods.
With ${K=5}$ only ${\mathrm{EOC} \approx 3.3}$ is achieved, which corresponds to a loss $2.7$ compared to the theoretical order.
A further increase in the number of sweeps reduces the error, while the convergence rate behaves irregularly and even decreases in between.
Nevertheless, the results for ${K=72}$ indicate that the optimal order can be restored using a sufficient number of corrections.
This leads to the conclusion that the SISDC methods can cope with time-dependent boundary conditions and high viscosity, but at the expense of increasing costs.
Whether the latter can be reduced by eliminating the splitting error remains an open question.

\begin{figure}
\begingroup
\captionsetup[subfigure]{width=70mm}
\subfloat[RK-ARS3 with different boundary conditions and linear vs nonlinear viscosities]
  {\includegraphics[scale=0.60]{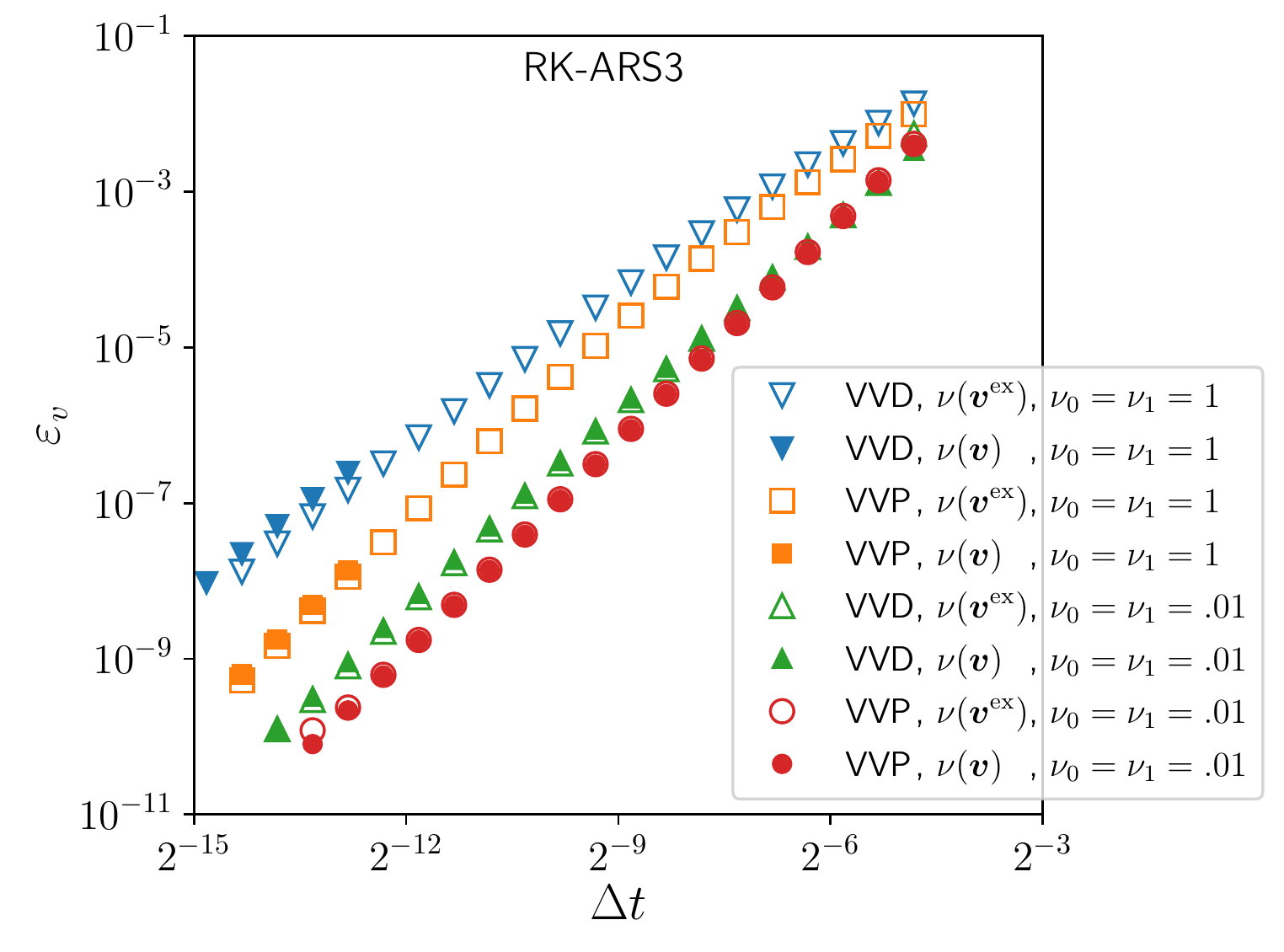}
   \label{fig:robustness:rk_ars3}}
\endgroup
\begingroup
\captionsetup[subfigure]{width=80mm}
\subfloat[SDC-Eu(3,$K$) with varying number of correction sweeps $K$ indicated by symbols]
  {\includegraphics[scale=0.60]{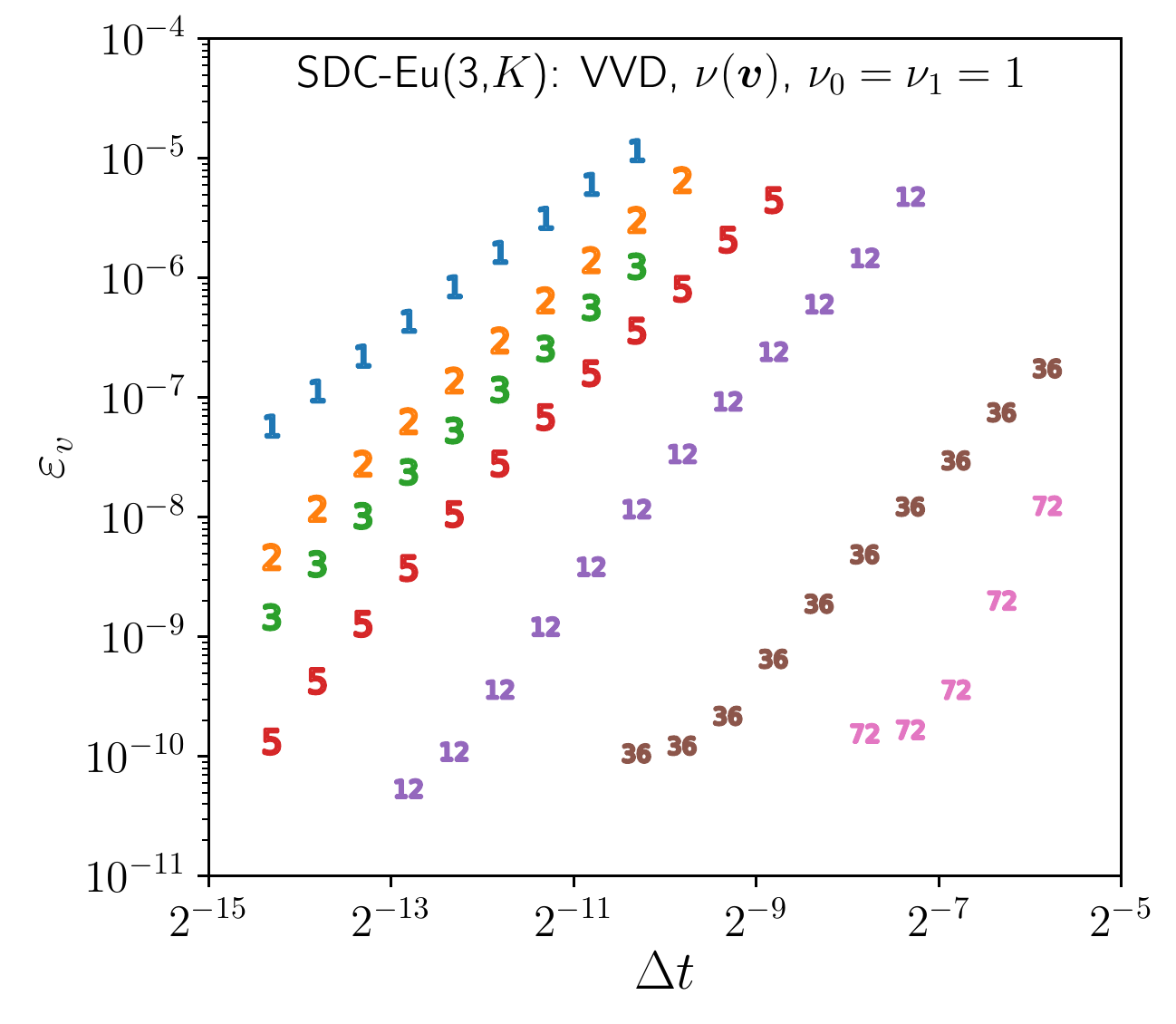}
   \label{fig:robustness:sdc_eu}}
\endgroup
\caption{\RII{Influence of boundary conditions and viscosity on the convergence behavior 
of time integration methods.}
  \label{fig:robustness}
  }
\end{figure}

}
\RII{ 

\subsubsection{Stability in the convection dominant case}
\label{sec:numerical-experiments:VV:stability}

To investigate the stability of the time integration methods in the convection dominant case, the viscosity coefficients are reduced to
  ${\nu_0 = \nu_1 = 10^{-3}}$, 
  which corresponds to a Reynolds number of 
  ${\mathit{Re} \approx 1333}$.
In order to determine a possible influence of the spatial approximation, the domain 
  ${\Omega = [-\sfrac{1}{2},\sfrac{1}{2}]^3}$
  is discretized in three ways:
  1) $6^3$ elements of degree ${P\!=\!6}$,
  2) $3^3$ elements of degree 11 and
  3) $2^3$ elements of degree 16.
Dirichlet conditions are specified at all boundaries.
The numerical tests are run until reaching the final time ${T = 10}$, 
  which corresponds to 10 convective units in terms of the phase velocity 
  of the vortex array.
A test is considered unstable when 
  exceeding a velocity magnitude of ${v^{\mathrm{ex}}_{\max}/2}$ or
  detecting a NaN in the numerical solution.
Based on this criterion the time step is adapted via bisection until
  a sufficiently accurate approximation of the critical 
  time step $\Delta t^\ast$ is reached.
The resulting stability threshold is converted into a CFL number and denoted as $\CClmb$.
Table ~\ref{tab:cfl} summarizes the critical CFL numbers for the time integration methods considered, with the exception of RK-TR, which proved unstable even for very small time steps.
Further, the table compiles the minimal errors achieved with 
${\Delta t \le \Delta t^{\ast}}$
and, in the last column, the theoretical threshold obtained from the model problem with comparable diffusivity.
For IMEX BDF2 no theoretical value is available.
Instead, the observed $\CClmb$ can be compared to the stability limit of the corresponding fully explicit method, which is ca $0.6$ for a real (diffusion) part of 
${\real{z} = -0.1}$, see e.g. \textcite[Ch.\,II,\,3.3]{TI_Hundsdorfer2003a}.
Taking this as a reference, the BDF2 method remains stable up to time steps of approximately 60 percent the size of the explicit threshold.
Similarly, RK-CB3e, RK-ARS3 and SDC-Eu(3,5) reach between 57 and 61 percent of the critical values obtained for the model problem.
Taking into account possible effects of nonlinearity, multidimensionality, spatial discretization and temporal splitting, these observations are in line with expectations.
RK-CB3e, SDC-CB3e and SDC-ARS3 reach only between 26 and 36 percent of the predicted threshold.
The reason for this behavior cannot be clearly identified from the present study.
However, it is worth noting that in the case of RK-CB3e, the lack of the GSA property prevents the accurate fulfillment of the boundary conditions, which can introduce additional disturbances.
With the SDC methods, the transition from a high-order RK predictor to the first-order the Euler method can yield a mismatch in the correction terms and thus result in further perturbations.
Finally it is noted, that BDF2 becomes unstable at errors of about two order of magnitudes larger than RK-CB3e or RK-ARS3, and 3 to 5 orders larger than SDC-Eu(3,5).
This indicates that the higher-order IMEX methods can take advantage only if sufficiently low errors are targeted, while the robustness of IMEX BDF2 makes it more suitable with relaxed accuracy requirements.

\begin{table}
\RII{%
  \caption{\RII{Critical CFL numbers and minimal errors for ${\mathit{Re} \approx 1333}$. 
  For reference, ${z_{\mathrm{im}}^{\ast}}$ gives the critical imaginary value for 
  the model problem studied in Sec.~\ref{sec:time-integration:stability} 
  assuming a real part of ${\real{z} = -10^{-1}}$.}
  \label{tab:cfl}}
  \begin{tabular}{lccccccccccc} \toprule
                 && \multicolumn{2}{c}{$P=6$} && \multicolumn{2}{c}{$P=11$} && \multicolumn{2}{c}{$P=16$} &&  \\
                    \cmidrule{3-4}               \cmidrule{6-7}                \cmidrule{9-10}
   Method        && $\CClmb$ & $\Emin$ && $\CClmb$ & $\Emin$ && $\CClmb$ & $\Emin$ && $z_{\mathrm{im}}^{\ast}$ \\ \midrule
   BDF2          && $0.34$ & $4.7\!\times\!10^{-3}$  && $0.37$ & $3.0\!\times\!10^{-3}$  && $0.36$ & $1.9\!\times\!10^{-3}$ &&        \\
   RK-CB2        && $0.47$ & $1.5\!\times\!10^{-3}$  && $0.66$ & $1.2\!\times\!10^{-3}$  && $0.68$ & $1.0\!\times\!10^{-3}$ && $1.18$ \\
   RK-CB3e       && $0.63$ & $5.8\!\times\!10^{-5}$  && $0.92$ & $5.8\!\times\!10^{-5}$  && $1.03$ & $6.3\!\times\!10^{-5}$ && $2.90$ \\
   RK-ARS3       && $0.70$ & $1.7\!\times\!10^{-5}$  && $1.00$ & $1.4\!\times\!10^{-5}$  && $1.06$ & $1.8\!\times\!10^{-5}$ && $1.73$ \\
   SDC-Eu(3,5)   && $1.52$ & $8.5\!\times\!10^{-6}$  && $1.58$ & $2.5\!\times\!10^{-8}$  && $1.52$ & $4.4\!\times\!10^{-8}$ && $2.68$ \\
   SDC-CB3e(3,3) && $0.97$ & $4.2\!\times\!10^{-6}$  && $1.27$ & $4.7\!\times\!10^{-4}$  && $1.08$ & $7.1\!\times\!10^{-5}$ && $4.20$ \\
   SDC-ARS3(3,3) && $0.80$ & $1.7\!\times\!10^{-5}$  && $1.04$ & $1.4\!\times\!10^{-5}$  && $1.08$ & $1.8\!\times\!10^{-5}$ && $2.99$ \\
  \bottomrule
  \end{tabular}
}
\end{table}

}


\subsection{Transition and turbulent decay of a Taylor-Green vortex}
\label{sec:numerical-experiments:TTG}

The last flow problem involves the nonlinear transition and turbulent decay of a disturbed Taylor-Green vortex \cite{Taylor1937a}.
Following \textcite{Brachet1991a}, the initial conditions are given by
\begin{equation}
\label{eq:transition:v0}
  \V v_0 = \begin{bmatrix}
             \phantom{-} \cos(x) \sin(y) \sin(z)  \\
                       - \sin(x) \cos(y) \sin(z)  \\
                         0
           \end{bmatrix}
  \,.
\end{equation}
This test case has been used in numerous studies for investigating transition mechanisms as well as for validating numerical methods, e.g.\
\cite{%
  Gassner2013a,%
  CartonDeWiart2014a,%
  Piatkowski2018a,%
  SE_Fehn2018a,%
  SE_Fehn2019a}.
The evolution of the vortex depends on the Reynolds number 
${\mathit{Re} = \nu^{-1}}$,
which is set to $1600$ for this study.
The simulations are run in the periodic domain ${\Omega = [-\pi,\pi]^3}$ from 
${t = 0}$ to ${T=20}$.
For spatial discretization, the domain is divided into $32^3$ cubic elements of degree ${P=15}$.
This yields approximately half the resolution of the reference study of 
Fehn et al.~\cite{SE_Fehn2018a}, who used $128^3$ elements of degree $7$.
Because of the expected higher spatial errors and the high computational cost, only the following methods are considered for time integration: BDF2, RK-ARS3 and RK-CB3e.
Depending on the method, the time step varies from $0.0002$ to $0.0024$ which corresponds to 
${\Clmb = 0.04}$  to $0.50$
\RII{based on the maximum initial velocity ${v_{0,\max} = 1}$,} 
respectively.
The computations are run until on 512 Intel Xeon E5-2680 v3 core until reaching the final time of 
${T=20}$.
\RI{%
Table~\ref{tab:DNS:overview} gives an overview of the simulations including the consumed runtime.}

\begin{table}\RI{%
  \newcommand{\tstep}{\bar t_{\mathrm{step}}}
  \newcommand{\phh}{\phantom{0}}
  \caption{\RI{%
   Overview of direct numerical simulations of the disturbed Taylor-Green vortex.
   The CFL number $\Clmb$ is based on the maximum initial velocity, 
   $\tstep$ is the average runtime required for a single time step.}
   }
  \label{tab:DNS:overview}
  \begin{tabular}{l@{\hspace{2em}}c@{\hspace{2em}}c@{\hspace{2em}}r@{\hspace{2em}}r}
  \toprule
  Method    & $\dt$  & $\Clmb$ & $\trun$ in h & $\tstep$ in s \\ \midrule
  BDF2      & 0.0002 &  0.040  &  118.4\phh   &    4.26\phh  \\
  BDF2      & 0.0012 &  0.250  &   16.7\phh   &    3.61\phh  \\
  RK-CB3e   & 0.0024 &  0.500  &   34.2\phh   &   14.77\phh  \\
  RK-ARS3   & 0.0008 &  0.167  &  126.3\phh   &   18.19\phh  \\
  RK-ARS3   & 0.0016 &  0.333  &   51.1\phh   &   14.71\phh  \\
  \bottomrule
  \end{tabular}
}
\end{table}

The accuracy of the results is assessed on the basis of the viscous dissipation rate $\epsilon$, the energy spectrum $E$ and selected snapshots of the instantaneous velocity field.
Since the first two represent statistical quantities, they can evaluated using the homogeneity of the flow:
The dissipation rate defined as the volume average of the local dissipation,
\begin{equation}
  \epsilon(t) 
   = \frac{\nu}{8\pi^3} \int_{\Omega} \nabla \V v : \nabla \V v \,\D \Omega
  \,.
\end{equation}
The energy spectrum $E(k,t)$ identifies the contributions of wave-number intervals ${[k-\frac{1}{2}, k+\frac{1}{2}]}$ for ${k=0,1,2,\ldots}$ such that
\begin{equation}
  \sum_{k \ge 0} E(k,t) 
  = \frac{1}{2}\int_{\Omega} |\V v(\V x, t)|^2 \,\D \Omega
  \,.
\end{equation}
It is computed by interpolating the velocity field to a regular grid,
performing a discrete 3D Fourier transform and 
summing up the contributions of wave vectors within each interval.

\begin{figure}
\subfloat[Evolution of the dissipation rate]
  {\begin{minipage}{0.51\textwidth}
   \raggedright
   \includegraphics[scale=0.68]{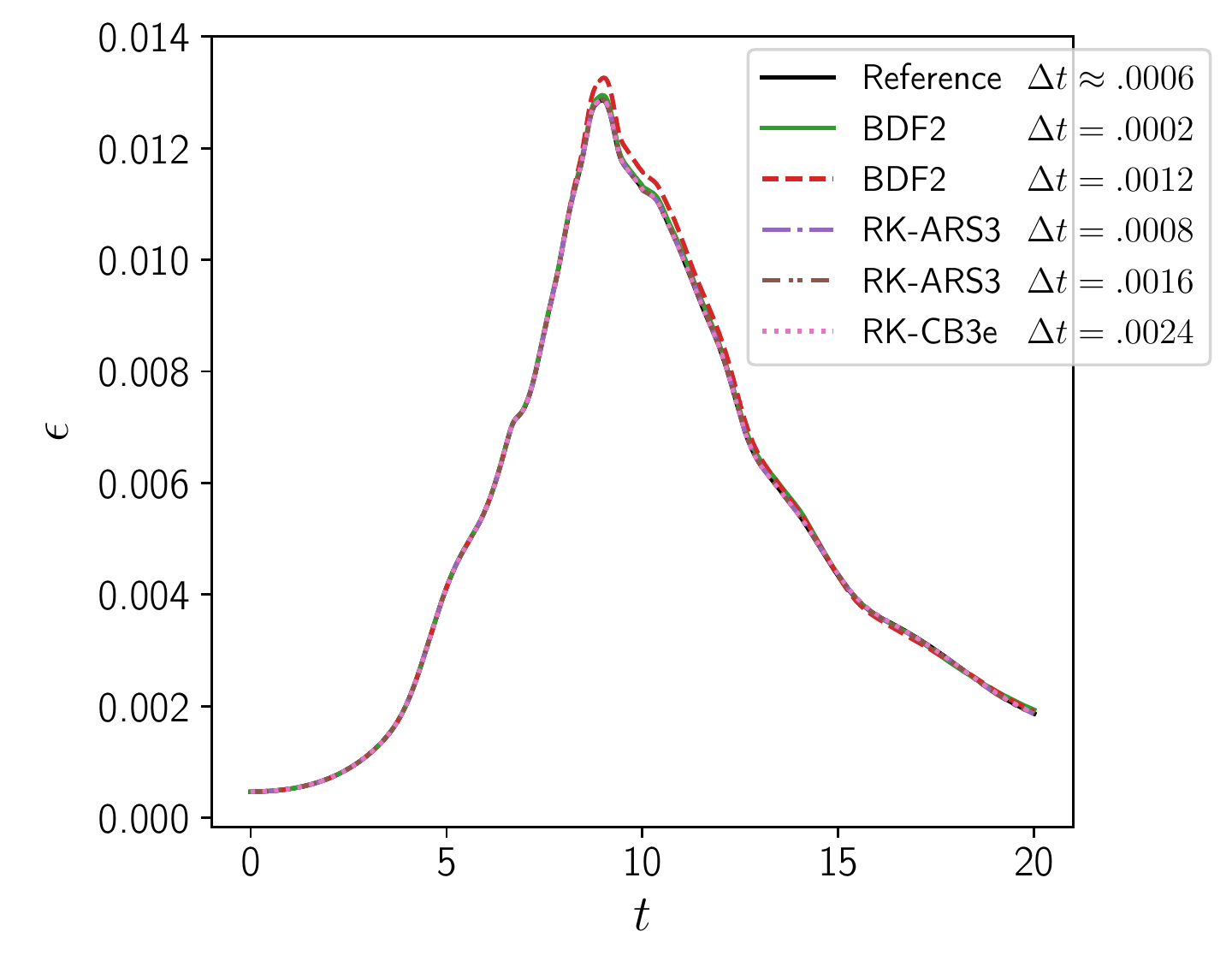}
   \end{minipage}
   \label{fig:DNS:statistics:dissipation}}
\subfloat[Energy spectra at ${t=8}$]
  {\begin{minipage}{0.49\textwidth}
   \raggedleft
   \includegraphics[scale=0.68]{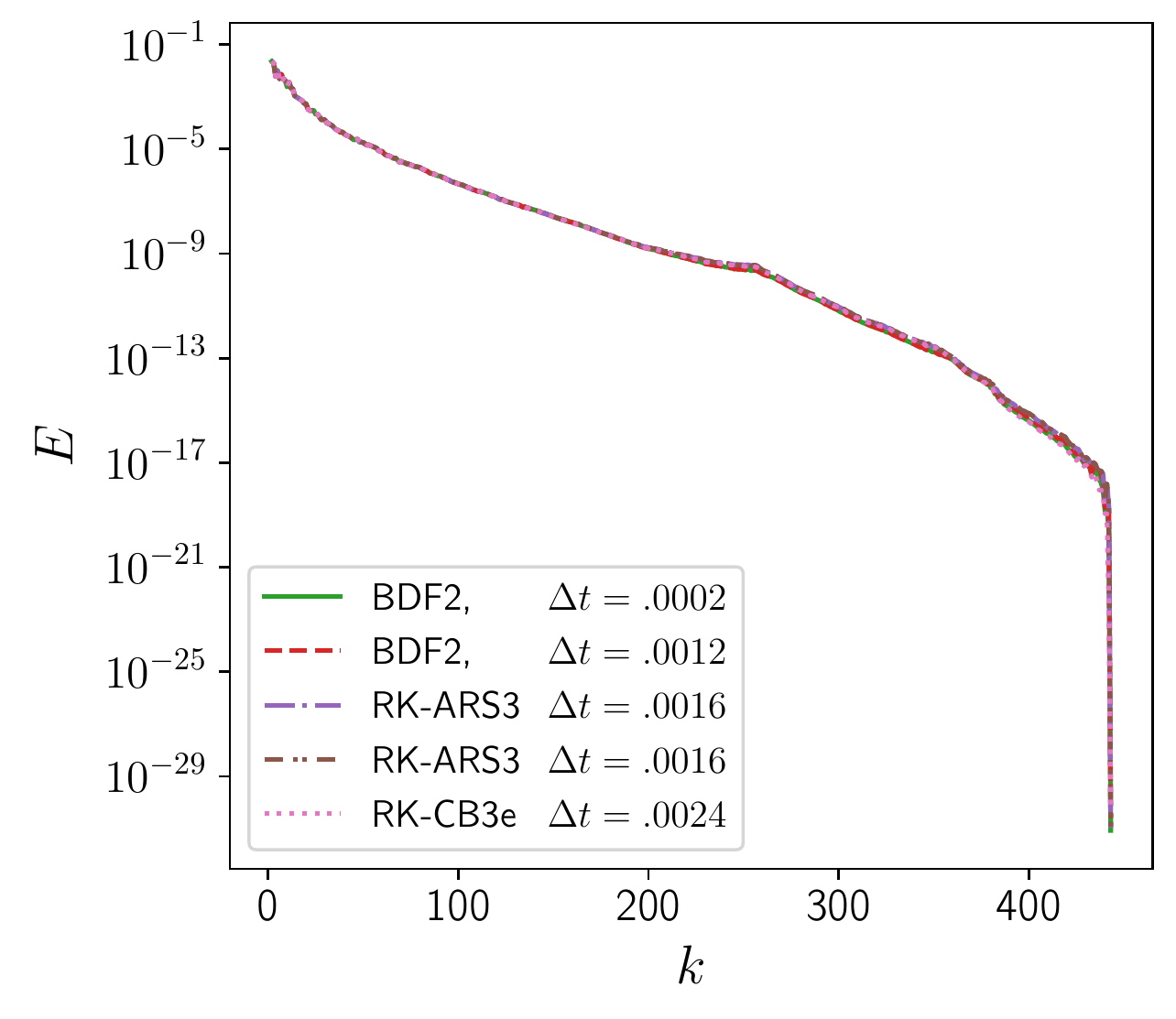}
   \end{minipage}
   \label{fig:DNS:statistics:spectrum}}
\caption{%
  Dissipation rate (a) and energy spectra (b) 
  for the decaying Taylor-Green vortex.
  The reference data for $\epsilon$ stems from \cite{SE_Fehn2018a}
  and was gratefully provided by Martin Kronbichler.
  \label{fig:DNS:statistics}
  }
\end{figure}

\begin{figure}
\begingroup
\captionsetup[subfigure]{width=60mm}
\subfloat[Top view of ${\lambda_2}$ contours for RK-ARS3 with ${\Delta t = 0.008}$]
  {\begin{minipage}{0.37\textwidth}
   \centering
   \includegraphics[height=60mm]{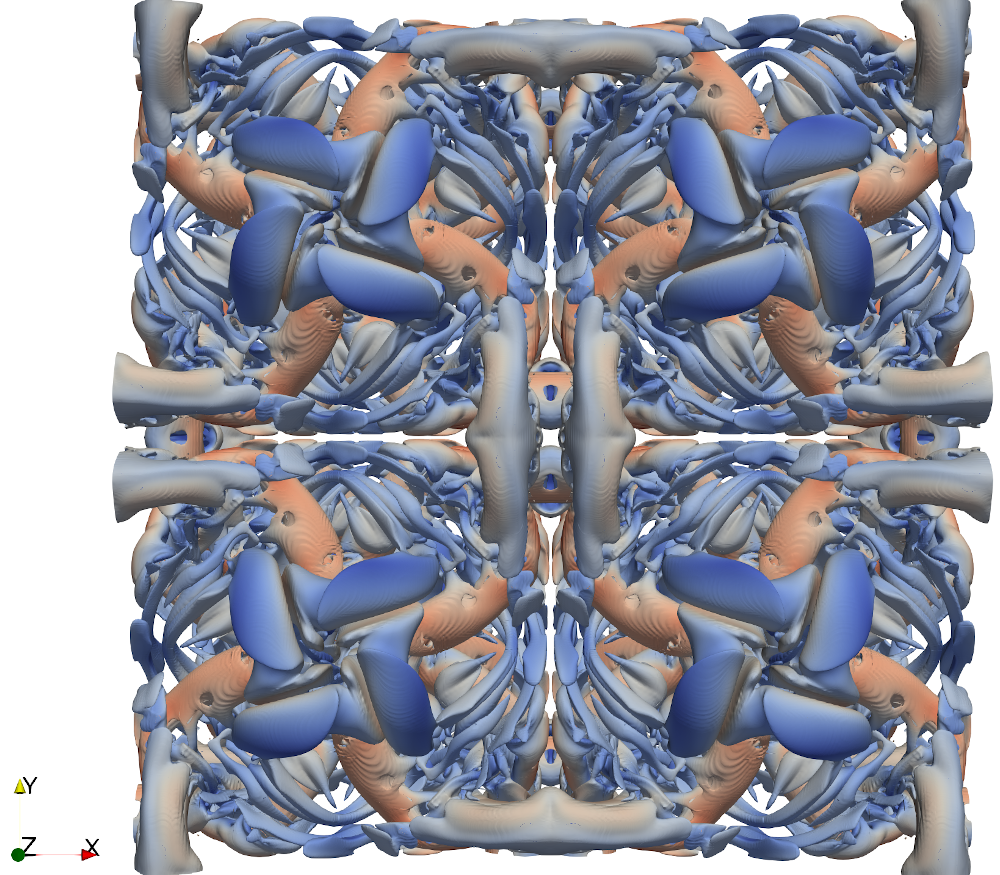}
   \end{minipage}
   \label{fig:DNS:lambda2:topview}}
\hfill
\subfloat[Contour lines for RK-ARS3 with enlarged region marked by red circle]
  {\begin{minipage}{0.37\textwidth}
   \centering
   \includegraphics[height=60mm]{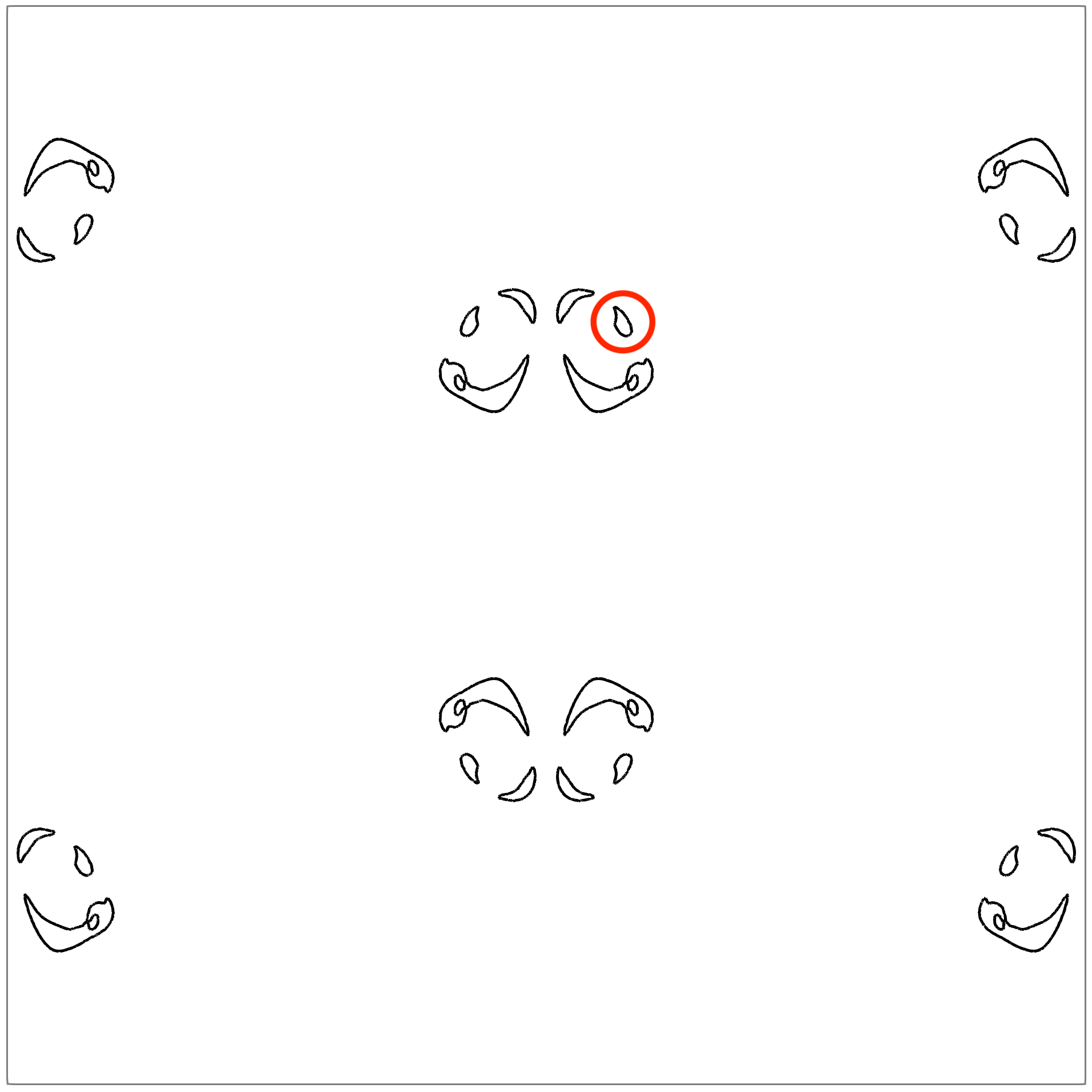}
   \end{minipage}
   \label{fig:DNS:lambda2:slice}}
\endgroup
\hfill
\begingroup
\captionsetup[subfigure]{width=45mm}
\subfloat[Comparison of methods in the marked region]
  {\begin{minipage}{0.23\textwidth}
   \centering
   \includegraphics[height=60mm]{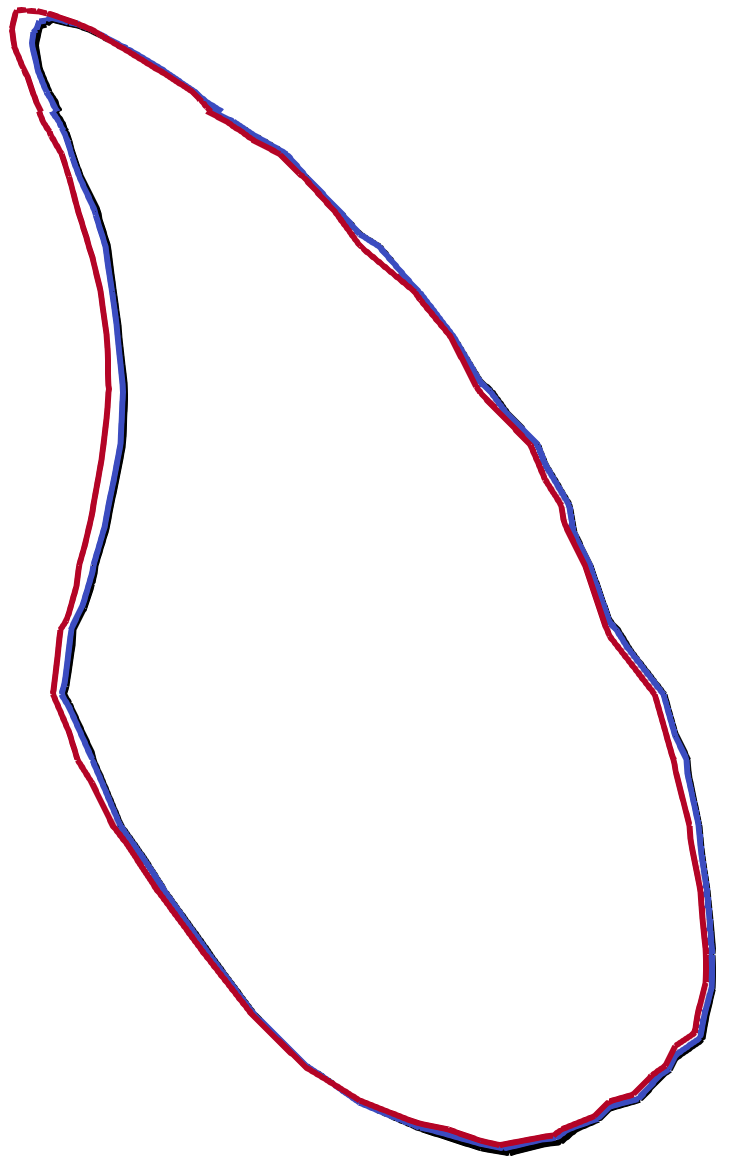}
   \end{minipage}
   \label{fig:DNS:lambda2:detail}}
\endgroup
\caption{%
  Visualization of ${\lambda_2 = -1.5}$ contours for the decaying Taylor-Green vortex 
  at ${t=8}$.
  The left panel (a) shows a top view the contours colored with the velocity magnitude 
  ranging from 0 (blue) to 1.5 (red)
  and the middle (b) depicts the contour lines in the plane ${x = 0}$.
  Both figures were obtained using RK-ARS3 with ${\Delta t = 0.008}$.
  The right panel (c) compares the contours obtained with different methods in the marked region.
  The black line corresponds to RK-ARS3 with ${\Delta t = 0008}$, 
  the virtually coinciding blue one to RK-CB3e with ${\Delta t = 0.024}$ and
  the red to BDF2 with ${\Delta t = 0.002}$.
  \label{fig:DNS:lambda2}
  }
\end{figure}

Figure~\ref{fig:DNS:statistics:dissipation} shows the dissipation rates obtained with the present simulations in comparison with the reference data from \cite{SE_Fehn2018a}.
Except for BDF2 with ${\Delta t = 0.00012}$, operating near the stability limit, all agree excellently with the latter.
At the peak, a maximum deviation of $10^{-5}$ is observed, which corresponds to a relative error less than $10^{-3}$.
Figure~\ref{fig:DNS:statistics:spectrum} reveals a similar agreement between the energy spectra at ${t=8}$, although no reference data is available.
The inspection of the peak values at ${k=2}$ yields a maximum relative deviation of about ${4\times10^{-4}}$.
Altogether the statistical results indicate relative errors between $10^{-3}$ and $10^{-4}$, but do not allow a clear distinction between the methods.
Moreover, the error level is in a range, where high-order time integration just starts to gain advantage.
In order to identify possible advantages of the Runge-Kutta methods, instantaneous velocity fields were compared with RK-ARS3 using 
${\Delta t = 0.0008}$ which, based on the previous test results, 
is assumed to be the most accurate method.
Figures~\ref{fig:DNS:lambda2:topview} and \ref{fig:DNS:lambda2:slice}
show a top view and a slice of the ${\lambda_2 = -1.5}$ contours 
for this method at time ${t = 8}$, respectively.
On this scale, the corresponding results obtained with the other methods
or time steps, respectively, can be hardly distinguished from the 
reference case.
Therefore, a detailed comparison is made in Fig.~\ref{fig:DNS:lambda2:detail}
for a single vortex filament.
It reveals an excellent agreement between the reference and RK-CB3e using a 
$\Delta t$ three times as large.
The contours obtained with RK-ARS3 using ${\Delta t = 0.0016}$ are not included, because they would coincide with these two cases.
In contrast, the contours of BDF2 using ${\Delta t = 0.0002}$ show a clear deviation from the reference, especially in the upper left part of the filament.
To allow for a quantitative assessment, the RMS deviation of the instantaneous velocity fields from the reference was computed at ${t=2,4,6,8}$ and compiled in Tab.~\ref{tab:DNS:velocity-deviation}.
Especially at larger times, ${t = 6}$ and $8$, the deviation of RK-ARS3 and RK-CB3e is 10 times smaller than that of BDF2, despite using 8 or 12 times the step size.
These observations confirm that temporal precision matters even in situations, where the spatial error is considerable.
In the present case of a marginally resolved direct numerical simulation, 
RK-CB3e produces significantly more accurate results than IMEX-BDF2 at 
\JS{less than one third} of the cost.

\begin{table}
  \caption{RMS deviation of instantaneous velocity from RK-ARS3 with ${\Delta t = 0.0008}$}.
  \label{tab:DNS:velocity-deviation}
  \begin{tabular}{l@{\hspace{2em}}c@{\hspace{2em}}c@{\hspace{2em}}c}
  \toprule
       & BDF2    
       & RK-ARS3
       & RK-CB3e \\
   $t$ & ${\Delta t = 0.0002}$ 
       & ${\Delta t = 0.0016}$
       & ${\Delta t = 0.0024}$ \\ \midrule
    2  & 1.02e\textminus08  & 1.40e\textminus10  & 1.40e\textminus10  \\
    4  & 1.67e\textminus05  & 4.39e\textminus06  & 5.44e\textminus06  \\
    6  & 3.36e\textminus04  & 2.89e\textminus05  & 3.79e\textminus05  \\
    8  & 2.43e\textminus03  & 2.03e\textminus04  & 2.27e\textminus04  \\
  \bottomrule
  \end{tabular}
\end{table}


\section{Conclusions}
\label{sec:conclusions}

The present study investigates the suitability of high-order implicit-explicit (IMEX) time integration methods for the simulation of incompressible flows.
It covers IMEX Runge-Kutta (RK) methods with convergence orders from 2 to 4 and semi-implicit spectral deferred correction (SDC) methods of order 6 and, for comparison, the IMEX difference formula of order 2 (BDF2).
The IMEX RK methods are adapted for use with a projection-based splitting scheme and extended to variable viscosity by operating as a partitioned Runge-Kutta method.
The SDC method is adopted from \cite{TI_Stiller2020a} and generalized for using IMEX RK predictors.
Three flow problems are selected to serve as test bench for examining the performance of these methods:
a 2D traveling Taylor-Green vortex,
a 3D manufactured vortex array with solution-dependent variable viscosity, and
the transition and turbulent decay of a Taylor-Green vortex.
Only in case of the 2D Taylor-Green vortex with periodic boundaries do all methods achieve their theoretical order.
With time-dependent Dirichlet conditions, the fourth-order and one of three third-order Runge-Kutta methods drop to order 2.
Using a fixed number of corrections, the SDC methods suffer a similar order reduction and reach experimental orders of 3.5 up to 5 instead of 6.
Except for the fourth-order RK, all methods cope with the nonlinear variable viscosity without any further degradation, which confirms the proposed semi-implicit approach.

As expected, SDC are methods by far the most accurate for a given size of the time step.
In terms of wall time per error, the RK and SDC methods are always more economic than BDF2.
The third-order RK-CB3e and RK-ARS3 are the best choice for relative errors down to approximately $10^{-5}$.
For errors below that they are surpassed by SDC methods, especially when using an RK predictor.
In case of the decaying Taylor-Green vortex, direct numerical simulations are performed with marginal spatial resolution.
Nonetheless, IMEX BDF2 and third-order RK methods achieved an excellent agreement with reference data from \cite{SE_Fehn2018a}.
However, a thorough comparison shows that RK-CB3e and RK-ARS3 are significantly more accurate than BDF2, even when using a 12-fold larger time step.
Translated into runtime, this corresponds to at least a saving of a factor of 4.

In summary, the study supports the use of high-order IMEX methods for incompressible Navier-Stokes problems.
Even with moderate accuracy requirements, IMEX Runge-Kutta methods present a much better choice than second-order schemes like BDF2.
For high-fidelity simulations, higher-order SDC methods offer an attractive alternative.
However, further research should be aimed at increasing their robustness and competitiveness, possibly by using multilevel techniques or Runge-Kutta correctors as proposed e.g. in \cite{TI_Speck2015a} and \cite{TI_Christlieb2011a}, respectively.


\appendix
\pagebreak

\section{Butcher tableaux of the IMEX Runge-Kutta methods}
\label{sec:butcher}

\renewcommand{\arraystretch}{1.25}

\subsection{RK-TR}
 
\begin{equation}
  \label{eq:butcher:RK-TR}
  \begin{array}{c|ccc}
    0 & 0           &   &             \\
    1 & 0           & 1 &             \\
    1 & \frac{1}{2} & 0 & \frac{1}{2} \\ \hline
      & \frac{1}{2} & 0 & \frac{1}{2} \\
  \end{array}
  \qquad
  \begin{array}{c|ccc}
    0 & 0           &             &   \\
    1 & 1           & 0           &   \\
    1 & \frac{1}{2} & \frac{1}{2} & 0 \\ \hline
      & \frac{1}{2} & \frac{1}{2} & 0 \\
  \end{array}
\end{equation}

\subsection{RK-CB2 \cite{TI_Cavaglieri2015a}}
 
\begin{equation}
  \label{eq:butcher:RK-CB2}
   \begin{array}{c|cccc}
     0           & 0 &               &             \\
     \frac{2}{5} & 0 &  \frac{2}{5}  &             \\
     1           & 0 &   \frac{5}{6} & \frac{1}{6} \\ \hline
        & 0  &   \frac{5}{6}  &  \frac{1}{6}
    \end{array}
    \qquad
  \begin{array}{c|cccc}
     0           & 0		        &              &      \\
     \frac{2}{5} & \frac{2}{5}  &   0          &      \\
     1           & 0  		      &   1          &  0   \\ \hline
                 & 0            & \frac{5}{6}  &  \frac{1}{6}
    \end{array} 
\end{equation}

\subsection{RK-CB3c \cite{TI_Cavaglieri2015a}}

\begin{subequations} 
\begin{align}
  \label{eq:butcher:RK-CB3c:im}
  &  
  \begin{array}{c|ccccc}
      0                                  & 0 &  &  &  \\
     \frac{3375509829940}{4525919076317} & 0 &  \frac{3375509829940}{4525919076317}  &   &  \\
     \frac{272778623835}{1039454778728}  & 0 & -\frac{11712383888607531889907}{32694570495602105556248} & \frac{566138307881}{912153721139} & \\
      1                                  & 0 &  \frac{673488652607}{2334033219546}                      & \frac{493801219040}{853653026979} & \frac{184814777513}{1389668723319} \\ \hline 
                                         & 0 &  \frac{673488652607}{2334033219546}                      & \frac{493801219040}{853653026979} & \frac{184814777513}{1389668723319} \\ 
    \end{array}  
\\[1\baselineskip]
  & 
  \begin{array}{c|ccccc}
     0                                   & 0                                   &   &  &    \\
     \frac{3375509829940}{4525919076317} & \frac{3375509829940}{4525919076317} & 0 &  & \\
     \frac{272778623835}{1039454778728}  & 0                                   & \frac{272778623835}{1039454778728} & 0  & \\
     1                                   & 0                                   & \frac{673488652607}{2334033219546} & \frac{1660544566939}{2334033219546} & 0 \\ \hline 
                                         & 0                                   & \frac{673488652607}{2334033219546} & \frac{493801219040}{853653026979}   & \frac{184814777513}{1389668723319} 
    \end{array}                 
\end{align}
\end{subequations} 

Remark: The original publication contains a typo that causes $\aim_{4,3}$ and $\aex_{4,3}$ to be swapped.

\subsection{RK-CB3e \cite{TI_Cavaglieri2015a}}
 
\begin{equation}
  \label{eq:butcher:RK-CB3e}
  \begin{array}{c|ccccc}
      0            & 0  &              &              &      \\
      \frac{1}{3}  & 0  & \frac{1}{3}  &              &      \\
      1            & 0  & \frac{1}{2}  &  \frac{1}{2} &      \\
      1            & 0  & \frac{3}{4}  & -\frac{1}{4} & \frac{1}{2} \\ \hline 
                   & 0   &  \frac{3}{4} & -\frac{1}{4} &  \frac{1}{2}      
  \end{array}
   \qquad  
  \begin{array}{c|ccccc}
      0            &  0            &               &   &      \\
      \frac{1}{3}  &  \frac{1}{3}  &  0            &   &      \\
      1            &  0            &  1            & 0 &      \\
      1            &  0            &  \frac{3}{4}  &  \frac{1}{4} &  0 \\ \hline 
                   &  0            &  \frac{3}{4}  & -\frac{1}{4} &  \frac{1}{2}     
  \end{array} 
\end{equation}

\subsection{RK-CB4 \cite{TI_Cavaglieri2015a}}

\begin{subequations} 
\begin{align}
  \label{eq:butcher:RK-CB4:im}
  &
   \begin{array}{c|cccccc}
     0           & 0          &       &   &   &   &   \\
     \frac{1}{4} & \frac{1}{8} &  \frac{1}{8}  &  &   &   &    \\
     \frac{3}{4} & \frac{216145252607}{961230882893}  &  \frac{257479850128}{1143310606989} &  \frac{30481561667}{101628412017}  & &   &     \\
     \frac{3}{8} & \frac{232049084587}{1377130630063} & -\frac{381180097479}{1276440792700} & -\frac{54660926949}{461115766612} & \frac{344309628413}{552073727558}   &  &   \\
     \frac{1}{2} & \frac{232049084587}{1377130630063} &  \frac{322009889509}{2243393849156}    & -\frac{100836174740}{861952129159} & -\frac{250423827953}{1283875864443}  &  \frac{1}{2} &    \\
     1           & \frac{232049084587}{1377130630063} &  \frac{322009889509}{2243393849156}    & -\frac{195109672787}{1233165545817} & -\frac{340582416761}{705418832319} & \frac{463396075661}{409972144477} & \frac{323177943294}{1626646580633}  \\ \hline 
                 & \frac{232049084587}{1377130630063} &  \frac{322009889509}{2243393849156}    &  -\frac{195109672787}{1233165545817} & -\frac{340582416761}{705418832319}   & \frac{463396075661}{409972144477} & \frac{323177943294}{1626646580633}     
 \end{array}  
\\[\baselineskip]
& \begin{array}{c|cccccc}    
     0           & 0            &   &  &  &  & \\
     \frac{1}{4} & \frac{1}{4}  & 0 &  &  &  & \\
     \frac{3}{4} & \frac{153985248130}{1004999853329} & \frac{902825336800}{1512825644809} &  0 &  &  &  \\
     \frac{3}{8} & \frac{232049084587}{1377130630063} & \frac{99316866929}{820744730663}   &  \frac{82888780751}{969573940619}   &  0 &   &      \\
     \frac{1}{2} & \frac{232049084587}{1377130630063} & \frac{322009889509}{2243393849156} &  \frac{57501241309}{765040883867}   &  \frac{76345938311}{676824576433}    & 0 & \\
     1           & \frac{232049084587}{1377130630063} & \frac{322009889509}{2243393849156} & -\frac{195109672787}{1233165545817} & -\frac{4099309936455}{6310162971841} & \frac{1395992540491}{933264948679} & 0 \\\hline 
                 & \frac{232049084587}{1377130630063} & \frac{322009889509}{2243393849156} & -\frac{195109672787}{1233165545817} & -\frac{340582416761}{705418832319}   & \frac{463396075661}{409972144477}  & \frac{323177943294}{1626646580633}       
 \end{array}                 
\end{align}
\end{subequations} 

\bigskip

\subsection{RK-ARS3 \cite{TI_Ascher1997a}}
 
\begin{equation}
  \label{eq:butcher:RK-ARS3}
  \begin{array}{c|cccccc}
    0           & 0 &               &              &             &  \\
    \frac{1}{2} & 0 &  \frac{1}{2}  &              &             &  \\
    \frac{2}{3} & 0 &  \frac{1}{6}  &  \frac{1}{2} &             &  \\
    \frac{1}{2} & 0 & -\frac{1}{2}  &  \frac{1}{2} & \frac{1}{2} &  \\
    1           & 0 &  \frac{3}{2}  & -\frac{3}{2} & \frac{1}{2} &\frac{1}{2} \\\hline 
                & 0 &  \frac{3}{2}  & -\frac{3}{2} & \frac{1}{2} &\frac{1}{2}
  \end{array}
   \qquad  
  \begin{array}{c|cccccc}
    0           & 0             &               &             &   &  \\
    \frac{1}{2} & \frac{1}{2}   &  0            &             &   &  \\
    \frac{2}{3} & \frac{11}{18} &  \frac{1}{18} & 0           &   &  \\
    \frac{1}{2} & \frac{5}{6}   & -\frac{5}{6}  & \frac{1}{2} & 0 &  \\
    1           & \frac{1}{4}   &  \frac{7}{4}  & \frac{3}{4} & -\frac{7}{4} & 0 \\\hline 
                & \frac{1}{4}   &  \frac{7}{4}  & \frac{3}{4} & -\frac{7}{4} & 0
  \end{array}
\end{equation}

\renewcommand{\arraystretch}{1.0}


\printbibliography 

@article{SE_Janssen2011a,
	author = {B{\"a}rbel Janssen and Guido Kanschat},
	doi = {10.1137/090778523},
	journal = {{SIAM} Journal on Scientific Computing},
	month = {jan},
	number = {4},
	pages = {2095--2114},
	publisher = {Society for Industrial {\&} Applied Mathematics ({SIAM})},
	title = {Adaptive Multilevel Methods with Local Smoothing for $H^1$- and $H^{\mathrm{curl}}$-Conforming High Order Finite Element Methods},
	url = {https://doi.org/10.1137%2F090778523},
	volume = {33},
	year = 2011}

@article{SE_Lottes2005a,
	author = {James W. Lottes and Paul F. Fischer},
	category = {SE},
	doi = {10.1007/s10915-004-4787-3},
	journal = {Journal of Scientific Computing},
	month = {jul},
	number = {1},
	pages = {45--78},
	publisher = {Springer Science and Business Media {LLC}},
	title = {Hybrid Multigrid/Schwarz Algorithms for the Spectral Element Method},
	url = {https://doi.org/10.1007%2Fs10915-004-4787-3},
	volume = {24},
	year = 2005}

@article{SE_Winters2018a,
	author = {Andrew R. Winters and Rodrigo C. Moura and Gianmarco Mengaldo and Gregor J. Gassner and Stefanie Walch and Joaquim Peiro and Spencer J. Sherwin},
	category = {SE},
	doi = {10.1016/j.jcp.2018.06.016},
	journal = {Journal of Computational Physics},
	month = {nov},
	pages = {1--21},
	publisher = {Elsevier {BV}},
	title = {A comparative study on polynomial dealiasing and split form discontinuous Galerkin schemes for under-resolved turbulence computations},
	url = {https://doi.org/10.1016%2Fj.jcp.2018.06.016},
	volume = {372},
	year = 2018}

@article{SE_Bassi2020a,
	author = {Francesco Bassi and Alessandro Colombo and Andrea Crivellini and Krzysztof J. Fidkowski and Matteo Franciolini and Antonio Ghidoni and Gianmaria Noventa},
	category = {SE},
	doi = {10.2514/1.j058847},
	journal = {{AIAA} Journal},
	month = {sep},
	number = {9},
	pages = {3963--3977},
	publisher = {American Institute of Aeronautics and Astronautics ({AIAA})},
	title = {Entropy-Adjoint p-Adaptive Discontinuous Galerkin Method for the Under-Resolved Simulation of Turbulent Flows},
	url = {https://doi.org/10.2514%2F1.j058847},
	volume = {58},
	year = 2020}

@article{SE_Fehn2018a,
	author = {Niklas Fehn and Wolfgang A. Wall and Martin Kronbichler},
	category = {SE},
	doi = {10.1002/fld.4511},
	journal = {International Journal for Numerical Methods in Fluids},
	month = {may},
	number = {1},
	pages = {32--54},
	publisher = {Wiley},
	title = {Efficiency of high-performance discontinuous Galerkin spectral element methods for under-resolved turbulent incompressible flows},
	url = {https://doi.org/10.1002%2Ffld.4511},
	volume = {88},
	year = 2018}

@article{SE_Joshi2016a,
	author = {Sumedh M. Joshi and Peter J. Diamessis and Derek T. Steinmoeller and Marek Stastna and Greg N. Thomsen},
	category = {SE},
	doi = {10.1016/j.compfluid.2016.04.021},
	journal = {Computers {\&} Fluids},
	month = {nov},
	pages = {120--129},
	publisher = {Elsevier {BV}},
	title = {A post-processing technique for stabilizing the discontinuous pressure projection operator in marginally-resolved incompressible inviscid flow},
	volume = {139},
	year = {2016}}

@article{SE_Akbas2018a,
	author = {Mine Akbas and Alexander Linke and Leo G. Rebholz and Philipp W. Schroeder},
	category = {SE},
	doi = {10.1016/j.cma.2018.07.019},
	journal = {Computer Methods in Applied Mechanics and Engineering},
	month = {nov},
	pages = {917--938},
	publisher = {Elsevier {BV}},
	title = {The analogue of grad-div stabilization in {DG} methods for incompressible flows: Limiting behavior and extension to tensor-product meshes},
	volume = {341},
	year = {2018}}

@article{SE_Stiller2016a,
	author = {J\"{o}rg Stiller},
	category = {SE},
	doi = {10.1007/s10915-016-0345-z},
	journal = {Journal of Scientific Computing},
	month = {dec},
	number = {1},
	pages = {81--96},
	publisher = {Springer Science and Business Media {LLC}},
	title = {Nonuniformly Weighted {S}chwarz Smoothers for Spectral Element Multigrid},
	volume = {72},
	year = {2016}}

@article{SE_Stiller2016b,
	author = {J\"{o}rg Stiller},
	category = {SE},
	doi = {10.1016/j.jcp.2016.09.041},
	journal = {Journal of Computational Physics},
	month = {dec},
	pages = {317--336},
	publisher = {Elsevier {BV}},
	title = {Robust multigrid for high-order discontinuous {G}alerkin methods: A fast {P}oisson solver suitable for high-aspect ratio Cartesian grids},
	volume = {327},
	year = {2016}}

@incollection{SE_Stiller2017a,
	author = {J\"{o}rg Stiller},
	booktitle = {Lecture Notes in Computational Science and Engineering},
	category = {SE},
	doi = {10.1007/978-3-319-65870-4_12},
	pages = {189--201},
	publisher = {Springer International Publishing},
	title = {Robust Multigrid for Cartesian Interior Penalty {DG} Formulations of the {P}oisson Equation in 3D},
	year = {2017}}

@book{SE_Canuto2011a,
	author = {Canuto, Claudio and Hussaini, M. Yousuff and Quarteroni, Alfio and Zang, Thomas A.},
	category = {SE},
	date = {2011-09-07},
	ean = {9783540307259},
	isbn = {3540307257},
	pagetotal = {592},
	publisher = {Springer Berlin Heidelberg},
	title = {Spectral Methods. Fundamentals in Single Domains},
	url = {https://www.ebook.de/de/product/5194670/claudio_canuto_m_yousuff_hussaini_alfio_quarteroni_thomas_a_zang_spectral_methods.html},
	year = {2011}}

@book{SE_Karniadakis2005a,
	author = {George Karniadakis and Spencer Sherwin},
	category = {SE},
	doi = {10.1093/acprof:oso/9780198528692.001.0001},
	month = {jun},
	publisher = {Oxford University Press},
	title = {Spectral/hp Element Methods for Computational Fluid Dynamics},
	year = {2005}}

@article{SE_Beck2014a,
	author = {Andrea D. Beck and Thomas Bolemann and David Flad and Hannes Frank and Gregor J. Gassner and Florian Hindenlang and Claus-Dieter Munz},
	category = {SE},
	doi = {10.1002/fld.3943},
	journal = {International Journal for Numerical Methods in Fluids},
	month = {aug},
	number = {8},
	pages = {522--548},
	publisher = {Wiley},
	title = {High-order discontinuous {G}alerkin spectral element methods for transitional and turbulent flow simulations},
	volume = {76},
	year = {2014}}

@article{SE_Schaal2015a,
	author = {Kevin Schaal and Andreas Bauer and Praveen Chandrashekar and R{\"u}diger Pakmor and Christian Klingenberg and Volker Springel},
	category = {SE},
	doi = {10.1093/mnras/stv1859},
	journal = {Monthly Notices of the Royal Astronomical Society},
	month = {sep},
	number = {4},
	pages = {4279--4301},
	publisher = {Oxford University Press ({OUP})},
	title = {Astrophysical hydrodynamics with a high-order discontinuous {G}alerkin scheme and adaptive mesh refinement},
	volume = {453},
	year = {2015}}

@article{SE_Fehn2019a,
	author = {Niklas Fehn and Martin Kronbichler and Christoph Lehrenfeld and Gert Lube and Philipp W. Schroeder},
	category = {SE},
	doi = {10.1002/fld.4763},
	journal = {International Journal for Numerical Methods in Fluids},
	month = {aug},
	number = {11},
	pages = {533--556},
	publisher = {Wiley},
	title = {High-order {DG} solvers for underresolved turbulent incompressible flows: A comparison of {$L^2$} and {$H(\operatorname{div})$} methods},
	volume = {91},
	year = {2019}}

@article{TI_Le1991a,
	author = {Hung Le and Parviz Moin},
	category = {TI},
	doi = {10.1016/0021-9991(91)90215-7},
	journal = {Journal of Computational Physics},
	month = {feb},
	number = {2},
	pages = {369--379},
	publisher = {Elsevier {BV}},
	title = {An improvement of fractional step methods for the incompressible Navier-Stokes equations},
	url = {https://doi.org/10.1016%2F0021-9991%2891%2990215-7},
	volume = {92},
	year = 1991}

@article{TI_Pareschi2005a,
	author = {Lorenzo Pareschi and Giovanni Russo},
	category = {TI},
	doi = {10.1007/s10915-004-4636-4},
	journal = {Journal of Scientific Computing},
	month = {oct},
	number = {1},
	pages = {129--155},
	publisher = {Springer Science and Business Media {LLC}},
	title = {Implicit{\textendash}Explicit Runge{\textendash}Kutta Schemes and Applications to Hyperbolic Systems with Relaxation},
	url = {https://doi.org/10.1007%2Fs10915-004-4636-4},
	volume = {25},
	year = 2005}

@article{TI_Kennedy2019a,
	author = {Christopher A. Kennedy and Mark H. Carpenter},
	category = {TI},
	doi = {10.1016/j.apnum.2018.10.007},
	journal = {Applied Numerical Mathematics},
	month = {feb},
	pages = {183--205},
	publisher = {Elsevier {BV}},
	title = {Higher-order additive Runge{\textendash}Kutta schemes for ordinary differential equations},
	url = {https://doi.org/10.1016%2Fj.apnum.2018.10.007},
	volume = {136},
	year = 2019}

@article{TI_Boscarino2017a,
	author = {Sebastiano Boscarino and Lorenzo Pareschi and Giovanni Russo},
	category = {TI},
	doi = {10.1137/m1111449},
	journal = {{SIAM} Journal on Numerical Analysis},
	month = {jan},
	number = {4},
	pages = {2085--2109},
	publisher = {Society for Industrial {\&} Applied Mathematics ({SIAM})},
	title = {A Unified {IMEX} Runge--Kutta Approach for Hyperbolic Systems with Multiscale Relaxation},
	url = {https://doi.org/10.1137%2Fm1111449},
	volume = {55},
	year = 2017}

@article{TI_Boscarino2013a,
	author = {S. Boscarino and L. Pareschi and G. Russo},
	category = {TI},
	doi = {10.1137/110842855},
	journal = {{SIAM} Journal on Scientific Computing},
	month = {jan},
	number = {1},
	pages = {A22--A51},
	publisher = {Society for Industrial {\&} Applied Mathematics ({SIAM})},
	title = {Implicit-Explicit Runge--Kutta Schemes for Hyperbolic Systems and Kinetic Equations in the Diffusion Limit},
	url = {https://doi.org/10.1137%2F110842855},
	volume = {35},
	year = 2013}

@article{TI_Boscarino2007a,
	author = {Sebastiano Boscarino},
	category = {TI},
	doi = {10.1137/060656929},
	journal = {{SIAM} Journal on Numerical Analysis},
	month = {jan},
	number = {4},
	pages = {1600--1621},
	publisher = {Society for Industrial {\&} Applied Mathematics ({SIAM})},
	title = {Error Analysis of {IMEX} Runge{\textendash}Kutta Methods Derived from Differential-Algebraic Systems},
	url = {https://doi.org/10.1137%2F060656929},
	volume = {45},
	year = 2007}

@article{TI_Bermudez2020a,
	author = {A. Berm{\'{u}}dez and S. Busto and M. Dumbser and J.L. Ferr{\'{\i}}n and L. Saavedra and M.E. V{\'{a}}zquez-Cend{\'{o}}n},
	category = {TI},
	doi = {10.1016/j.jcp.2020.109743},
	journal = {Journal of Computational Physics},
	month = {nov},
	pages = {109743},
	publisher = {Elsevier {BV}},
	title = {A staggered semi-implicit hybrid {FV}/{FE} projection method for weakly compressible flows},
	url = {https://doi.org/10.1016%2Fj.jcp.2020.109743},
	volume = {421},
	year = 2020}

@article{TI_Busto2018a,
	author = {S. Busto and J.L. Ferr{\'{\i}}n and E.F. Toro and M.E. V{\'{a}}zquez-Cend{\'{o}}n},
	category = {TI},
	doi = {10.1016/j.jcp.2017.10.004},
	journal = {Journal of Computational Physics},
	month = {jan},
	pages = {169--192},
	publisher = {Elsevier {BV}},
	title = {A projection hybrid high order finite volume/finite element method for incompressible turbulent flows},
	url = {https://doi.org/10.1016%2Fj.jcp.2017.10.004},
	volume = {353},
	year = 2018}

@article{TI_Sanderse2013a,
	author = {B. Sanderse},
	category = {RK},
	doi = {10.1016/j.jcp.2012.07.039},
	journal = {Journal of Computational Physics},
	month = {jan},
	pages = {100--131},
	publisher = {Elsevier {BV}},
	title = {Energy-conserving Runge{\textendash}Kutta methods for the incompressible Navier{\textendash}Stokes equations},
	url = {https://doi.org/10.1016%2Fj.jcp.2012.07.039},
	volume = {233},
	year = 2013}

@article{TI_Boscarino2009a,
	author = {Sebastiano Boscarino},
	category = {TI},
	doi = {10.1016/j.apnum.2008.10.003},
	journal = {Applied Numerical Mathematics},
	month = {jul},
	number = {7},
	pages = {1515--1528},
	publisher = {Elsevier {BV}},
	title = {On an accurate third order implicit-explicit Runge{\textendash}Kutta method for stiff problems},
	url = {https://doi.org/10.1016%2Fj.apnum.2008.10.003},
	volume = {59},
	year = 2009}

@article{TI_Hagstrom2006a,
	author = {Thomas Hagstrom and Ruhai Zhou},
	category = {TI},
	doi = {10.2140/camcos.2006.1.169},
	journal = {Communications in Applied Mathematics and Computational Science},
	month = {dec},
	number = {1},
	pages = {169--205},
	publisher = {Mathematical Sciences Publishers},
	title = {On the spectral deferred correction of splitting methods for initial value problems},
	url = {https://doi.org/10.2140%2Fcamcos.2006.1.169},
	volume = {1},
	year = 2006}

@article{TI_Boscarino2018a,
	author = {Sebastiano Boscarino and Jing-Mei Qiu and Giovanni Russo},
	category = {TI},
	doi = {10.1137/16m1105232},
	journal = {{SIAM} Journal on Scientific Computing},
	month = {jan},
	number = {2},
	pages = {A787--A816},
	publisher = {Society for Industrial {\&} Applied Mathematics ({SIAM})},
	title = {Implicit-Explicit Integral Deferred Correction Methods for Stiff Problems},
	url = {https://doi.org/10.1137%2F16m1105232},
	volume = {40},
	year = 2018}

@article{TI_Frank1997a,
	author = {J. Frank and W. Hundsdorfer and J.G. Verwer},
	category = {TI},
	doi = {10.1016/s0168-9274(97)00059-7},
	journal = {Applied Numerical Mathematics},
	month = {nov},
	number = {2-3},
	pages = {193--205},
	publisher = {Elsevier {BV}},
	title = {On the stability of implicit-explicit linear multistep methods},
	url = {https://doi.org/10.1016%2Fs0168-9274%2897%2900059-7},
	volume = {25},
	year = 1997}

@article{TI_Kanevsky2007a,
	author = {Alex Kanevsky and Mark H. Carpenter and David Gottlieb and Jan S. Hesthaven},
	category = {TI},
	doi = {10.1016/j.jcp.2007.02.021},
	journal = {Journal of Computational Physics},
	month = {aug},
	number = {2},
	pages = {1753--1781},
	publisher = {Elsevier {BV}},
	title = {Application of implicit{\textendash}explicit high order Runge{\textendash}Kutta methods to disconti\-nuous{\textendash}Galerkin schemes},
	url = {https://doi.org/10.1016%2Fj.jcp.2007.02.021},
	volume = {225},
	year = 2007}

@article{TI_Nikitin2006a,
	author = {Nikolay Nikitin},
	category = {TI},
	doi = {10.1002/fld.1122},
	journal = {International Journal for Numerical Methods in Fluids},
	number = {2},
	pages = {221--233},
	publisher = {Wiley},
	title = {Third-order-accurate semi-implicit Runge{\textendash}Kutta scheme for incompressible Na\-vier{\textendash}Stokes equations},
	url = {https://doi.org/10.1002%2Ffld.1122},
	volume = {51},
	year = {2006}}

@article{TI_Jothiprasad2003a,
	author = {Giridhar Jothiprasad and Dimitri J. Mavriplis and David A. Caughey},
	category = {TI},
	doi = {10.1016/s0021-9991(03)00330-9},
	journal = {Journal of Computational Physics},
	month = {nov},
	number = {2},
	pages = {542--566},
	publisher = {Elsevier {BV}},
	title = {Higher-order time integration schemes for the unsteady Navier{\textendash}Stokes equations on unstructured meshes},
	url = {https://doi.org/10.1016%2Fs0021-9991%2803%2900330-9},
	volume = {191},
	year = {2003}}

@inproceedings{TI_Persson2010a,
	author = {Per-Olof Persson and David Willis and Jaime Peraire},
	booktitle = {48th {AIAA} Aerospace Sciences Meeting Including the New Horizons Forum and Aerospace Exposition},
	category = {TI},
	doi = {10.2514/6.2010-724},
	month = {jan},
	publisher = {American Institute of Aeronautics and Astronautics},
	title = {The Numerical Simulation of Flapping Wings at Low Reynolds Numbers},
	url = {https://doi.org/10.2514%2F6.2010-724},
	year = {2010}}

@article{TI_Pazner2017a,
	author = {Will Pazner and Per-Olof Persson},
	category = {TI},
	doi = {10.1016/j.jcp.2017.01.050},
	journal = {Journal of Computational Physics},
	month = {apr},
	pages = {700--717},
	publisher = {Elsevier {BV}},
	title = {Stage-parallel fully implicit Runge{\textendash}Kutta solvers for discontinuous Galerkin fluid simulations},
	url = {https://doi.org/10.1016%2Fj.jcp.2017.01.050},
	volume = {335},
	year = {2017}}

@article{TI_Pan2021a,
	author = {Yu Pan and Zhen-Guo Yan and Joaquim Peir{\'{o}} and Spencer J. Sherwin},
	category = {TI},
	doi = {10.1007/s42967-021-00138-1},
	journal = {Communications on Applied Mathematics and Computation},
	month = {aug},
	publisher = {Springer Science and Business Media {LLC}},
	title = {Development of a Balanced Adaptive Time-Stepping Strategy Based on an Implicit {JFNK}-{DG} Compressible Flow Solver},
	url = {https://doi.org/10.1007%2Fs42967-021-00138-1},
	year = {2021}}

@article{TI_Stiller2020a,
	author = {J{\"o}rg Stiller},
	category = {TI},
	doi = {10.1016/j.jcp.2020.109840},
	journal = {Journal of Computational Physics},
	month = {dec},
	pages = {109840},
	publisher = {Elsevier {BV}},
	title = {A spectral deferred correction method for incompressible flow with variable viscosity},
	url = {https://doi.org/10.1016%2Fj.jcp.2020.109840},
	volume = {423},
	year = {2020}}

@book{TI_Hundsdorfer2003a,
	author = {Willem Hundsdorfer and Jan Verwer},
	category = {TI},
	doi = {10.1007/978-3-662-09017-6},
	publisher = {Springer Berlin Heidelberg},
	title = {Numerical Solution of Time-Dependent Advection-Diffusion-Reaction Equations},
	url = {https://doi.org/10.1007%2F978-3-662-09017-6},
	year = {2003}}

@techreport{TI_Kennedy2016a,
	author = {Christopher A. Kennedy and Mark H. Carpenter},
	category = {TI},
	institution = {NASA},
	month = {dec},
	number = {2016-219173},
	title = {Diagonally implicit Runge{\textendash}Kutta methods for Ordinary Differential Equations: A Review},
	type = {Technical Report NASA/TM},
	year = {2016}}

@article{TI_Ong2020a,
	author = {Ong, Benjamin W. and Spiteri, Raymond J.},
	category = {TI},
	doi = {10.1007/s10915-020-01235-8},
	issn = {1573-7691},
	journal = {Journal of Scientific Computing},
	month = {Jun},
	number = {3},
	publisher = {Springer Science and Business Media LLC},
	title = {Deferred Correction Methods for Ordinary Differential Equations},
	url = {http://dx.doi.org/10.1007/s10915-020-01235-8},
	volume = {83},
	year = {2020}}

@article{TI_Layton2007a,
	author = {Layton, Anita and Minion, Michael},
	category = {TI},
	doi = {10.2140/camcos.2007.2.1},
	issn = {1559-3940},
	journal = {Communications in Applied Mathematics and Computational Science},
	month = {Aug},
	number = {1},
	pages = {1--34},
	publisher = {Mathematical Sciences Publishers},
	title = {Implications of the choice of predictors for semi-implicit Picard integral deferred correction methods},
	url = {http://dx.doi.org/10.2140/camcos.2007.2.1},
	volume = {2},
	year = {2007}}

@article{TI_Minion2004a,
	author = {Minion, Michael L.},
	category = {TI},
	doi = {10.1016/j.apnum.2003.11.005},
	issn = {0168-9274},
	journal = {Applied Numerical Mathematics},
	month = {Mar},
	number = {3-4},
	pages = {369--387},
	publisher = {Elsevier BV},
	title = {Semi-implicit projection methods for incompressible flow based on spectral deferred corrections},
	url = {http://dx.doi.org/10.1016/j.apnum.2003.11.005},
	volume = {48},
	year = {2004}}

@article{TI_Speck2015a,
	author = {Speck, Robert and Ruprecht, Daniel and Emmett, Matthew and Minion, Michael and Bolten, Matthias and Krause, Rolf},
	category = {TI},
	doi = {10.1007/s10543-014-0517-x},
	issn = {1572-9125},
	journal = {BIT Numerical Mathematics},
	month = {Aug},
	number = {3},
	pages = {843--867},
	publisher = {Springer Science and Business Media LLC},
	title = {A multi-level spectral deferred correction method},
	url = {http://dx.doi.org/10.1007/s10543-014-0517-x},
	volume = {55},
	year = {2015}}

@article{TI_Christlieb2011a,
	author = {Christlieb, Andrew and Morton, Maureen and Ong, Benjamin and Qiu, Jing-Mei},
	category = {TI},
	doi = {10.4310/cms.2011.v9.n3.a10},
	issn = {1945-0796},
	journal = {Communications in Mathematical Sciences},
	number = {3},
	pages = {879--902},
	publisher = {International Press of Boston},
	title = {Semi-implicit integral deferred correction constructed with additive {R}unge--{K}utta methods},
	url = {http://dx.doi.org/10.4310/CMS.2011.v9.n3.a10},
	volume = {9},
	year = {2011}}

@article{TI_Ascher1997a,
	author = {Ascher, Uri M. and Ruuth, Steven J. and Spiteri, Raymond J.},
	category = {TI},
	doi = {10.1016/s0168-9274(97)00056-1},
	issn = {0168-9274},
	journal = {Applied Numerical Mathematics},
	month = {Nov},
	number = {2-3},
	pages = {151--167},
	publisher = {Elsevier BV},
	title = {Implicit-explicit {R}unge-{K}utta methods for time-dependent partial differential equations},
	url = {http://dx.doi.org/10.1016/S0168-9274(97)00056-1},
	volume = {25},
	year = {1997}}

@article{TI_John2006a,
	author = {John, Volker and Matthies, Gunar and Rang, Joachim},
	category = {TI},
	doi = {10.1016/j.cma.2005.10.007},
	issn = {0045-7825},
	journal = {Computer Methods in Applied Mechanics and Engineering},
	month = {Sep},
	number = {44-47},
	pages = {5995--6010},
	publisher = {Elsevier BV},
	title = {A comparison of time-discretization/linearization approaches for the incompressible {N}avier--{S}tokes equations},
	url = {http://dx.doi.org/10.1016/j.cma.2005.10.007},
	volume = {195},
	year = {2006}}

@article{TI_Noventa2016a,
	author = {Noventa, G. and Massa, F. and Bassi, F. and Colombo, A. and Franchina, N. and Ghidoni, A.},
	category = {TI},
	doi = {10.1016/j.compfluid.2016.03.007},
	issn = {0045-7930},
	journal = {Computers {\&} Fluids},
	month = {Nov},
	pages = {248--260},
	publisher = {Elsevier BV},
	title = {A high-order Discontinuous {G}alerkin solver for unsteady incompressible turbulent flows},
	url = {http://dx.doi.org/10.1016/j.compfluid.2016.03.007},
	volume = {139},
	year = {2016}}

@article{TI_Uranga2010a,
	author = {Uranga, A. and Persson, P.-O. and Drela, M. and Peraire, J.},
	category = {TI},
	doi = {10.1002/nme.3036},
	issn = {0029-5981},
	journal = {International Journal for Numerical Methods in Engineering},
	month = {Oct},
	number = {1-5},
	pages = {232--261},
	publisher = {Wiley},
	title = {Implicit Large Eddy Simulation of transition to turbulence at low Reynolds numbers using a Discontinuous {G}alerkin method},
	url = {http://dx.doi.org/10.1002/nme.3036},
	volume = {87},
	year = {2010}}

@article{TI_Leriche2006a,
	author = {Leriche, E. and Perchat, E. and Labrosse, G. and Deville, M. O.},
	category = {TI},
	doi = {10.1007/s10915-004-4798-0},
	issn = {1573-7691},
	journal = {Journal of Scientific Computing},
	month = {Jan},
	number = {1},
	pages = {25--43},
	publisher = {Springer Science and Business Media LLC},
	title = {Numerical Evaluation of the Accuracy and Stability Properties of High-order Direct {S}tokes Solvers with or without Temporal Splitting},
	url = {http://dx.doi.org/10.1007/s10915-004-4798-0},
	volume = {26},
	year = {2006}}

@article{TI_Persson2011a,
	author = {Persson, Per-Olof},
	category = {TI},
	doi = {10.2514/6.2011-684},
	isbn = {9781600869501},
	journal = {49th AIAA Aerospace Sciences Meeting including the New Horizons Forum and Aerospace Exposition},
	month = {Jan},
	publisher = {American Institute of Aeronautics and Astronautics},
	title = {High-Order LES Simulations using Implicit-Explicit {R}unge-{K}utta Schemes},
	url = {http://dx.doi.org/10.2514/6.2011-684},
	year = {2011}}

@article{TI_Bassi2015a,
	author = {Bassi, F. and Botti, L. and Colombo, A. and Ghidoni, A. and Massa, F.},
	category = {TI},
	doi = {10.1016/j.compfluid.2015.06.007},
	issn = {0045-7930},
	journal = {Computers {\&} Fluids},
	month = {Sep},
	pages = {305--320},
	publisher = {Elsevier BV},
	title = {Linearly implicit {R}osenbrock-type {R}unge--{K}utta schemes applied to the Discontinuous {G}alerkin solution of compressible and incompressible unsteady flows},
	url = {http://dx.doi.org/10.1016/j.compfluid.2015.06.007},
	volume = {118},
	year = {2015}}

@article{TI_Guermond2003a,
	author = {J.L. Guermond and Jie Shen},
	category = {TI},
	doi = {10.1016/j.jcp.2003.07.009},
	journal = {Journal of Computational Physics},
	month = {nov},
	number = {1},
	pages = {262--276},
	publisher = {Elsevier {BV}},
	title = {A new class of truly consistent splitting schemes for incompressible flows},
	volume = {192},
	year = {2003}}

@article{TI_Guermond2003b,
	author = {J. L. Guermond and Jie Shen},
	category = {TI},
	doi = {10.1137/s0036142901395400},
	journal = {{SIAM} Journal on Numerical Analysis},
	month = {jan},
	number = {1},
	pages = {112--134},
	publisher = {Society for Industrial {\&} Applied Mathematics ({SIAM})},
	title = {Velocity-Correction Projection Methods for Incompressible Flows},
	volume = {41},
	year = {2003}}

@article{TI_Dutt2000a,
	author = {Alok Dutt and Leslie Greengard and Vladimir Rokhlin},
	category = {TI},
	doi = {10.1023/a:1022338906936},
	journal = {Bit Numerical Mathematics},
	number = {2},
	pages = {241--266},
	publisher = {Springer Nature},
	title = {Spectral Deferred Correction Methods for Ordinary Differential Equations},
	volume = {40},
	year = {2000}}

@article{TI_Minion2003b,
	author = {Michael L. Minion},
	category = {TI},
	doi = {10.4310/cms.2003.v1.n3.a6},
	journal = {Communications in Mathematical Sciences},
	number = {3},
	pages = {471--500},
	publisher = {International Press of Boston},
	title = {Semi-implicit spectral deferred correction methods for ordinary differential equations},
	volume = {1},
	year = {2003}}

@article{TI_Minion2018a,
	author = {M.L. Minion and R.I. Saye},
	category = {TI},
	doi = {10.1016/j.jcp.2018.08.054},
	journal = {Journal of Computational Physics},
	month = {dec},
	pages = {797--822},
	publisher = {Elsevier {BV}},
	title = {Higher-order temporal integration for the incompressible {Na}vier-{S}tokes equations in bounded domains},
	volume = {375},
	year = {2018}}

@article{TI_Cavaglieri2015a,
	author = {Daniele Cavaglieri and Thomas Bewley},
	category = {TI},
	doi = {10.1016/j.jcp.2015.01.031},
	journal = {Journal of Computational Physics},
	month = {apr},
	pages = {172--193},
	publisher = {Elsevier {BV}},
	title = {Low-storage implicit/explicit {R}unge-{K}utta schemes for the simulation of stiff high-dimensional {ODE} systems},
	volume = {286},
	year = {2015}}

@article{TI_Kennedy2003a,
	author = {Christopher A. Kennedy and Mark H. Carpenter},
	category = {TI},
	doi = {10.1016/s0168-9274(02)00138-1},
	journal = {Applied Numerical Mathematics},
	month = {jan},
	number = {1-2},
	pages = {139--181},
	publisher = {Elsevier {BV}},
	title = {Additive {R}unge-{K}utta schemes for convection-diffusion-reaction equations},
	volume = {44},
	year = {2003}}

@article{TI_Fehn2017a,
	author = {Niklas Fehn and Wolfgang A. Wall and Martin Kronbichler},
	category = {TI},
	doi = {10.1016/j.jcp.2017.09.031},
	journal = {Journal of Computational Physics},
	month = {dec},
	pages = {392--421},
	publisher = {Elsevier {BV}},
	title = {On the stability of projection methods for the incompressible {N}avier-{S}tokes equations based on high-order discontinuous {G}alerkin discretizations},
	volume = {351},
	year = {2017}}

@article{TI_Karniadakis1991a,
	author = {George Em Karniadakis and Moshe Israeli and Steven A Orszag},
	category = {TI},
	doi = {10.1016/0021-9991(91)90007-8},
	journal = {Journal of Computational Physics},
	month = {dec},
	number = {2},
	pages = {414--443},
	publisher = {Elsevier {BV}},
	title = {High-order splitting methods for the incompressible {N}avier-{S}tokes equations},
	volume = {97},
	year = {1991}}

@article{TI_Klein2015a,
	author = {B. Klein and F. Kummer and M. Keil and M. Oberlack},
	category = {TI},
	doi = {10.1002/fld.3994},
	journal = {International Journal for Numerical Methods in Fluids},
	month = {jan},
	number = {10},
	pages = {571--589},
	publisher = {Wiley},
	title = {An extension of the {SIMPLE} based discontinuous {G}alerkin solver to unsteady incompressible flows},
	volume = {77},
	year = {2015}}

@article{TI_Colomes2016a,
	author = {Oriol Colom{\'{e}}s and Santiago Badia},
	category = {TI},
	doi = {10.1002/nme.4987},
	journal = {International Journal for Numerical Methods in Engineering},
	month = {sep},
	number = {5},
	pages = {372--400},
	publisher = {Wiley},
	title = {Segregated {R}unge-{K}utta methods for the incompressible {N}avier-{S}tokes equations},
	volume = {105},
	year = {2016}}

@article{TI_Ahmed2017a,
	author = {Naveed Ahmed and Simon Becher and Gunar Matthies},
	category = {TI},
	doi = {10.1016/j.cma.2016.09.026},
	journal = {Computer Methods in Applied Mechanics and Engineering},
	month = {jan},
	pages = {28--52},
	publisher = {Elsevier {BV}},
	title = {Higher-order discontinuous {G}alerkin time stepping and local projection stabilization techniques for the transient {S}tokes problem},
	volume = {313},
	year = {2017}}

@book{TI_Hairer1993a,
	author = {Ernst Hairer and Syver P. N{\o}rsett and Gerhard Wanner},
	category = {TI},
	doi = {10.1007/978-3-540-78862-1},
	publisher = {Springer Berlin Heidelberg},
	title = {Solving Ordinary Differential Equations I},
	year = {1993}}

@book{TI_Hairer1996a,
	author = {Ernst Hairer and Gerhard Wanner},
	category = {TI},
	doi = {10.1007/978-3-642-05221-7},
	publisher = {Springer Berlin Heidelberg},
	title = {Solving Ordinary Differential Equations {II}},
	year = {1996}}

@article{TI_Boscarino2016a,
	author = {Sebastiano Boscarino and Francis Filbet and Giovanni Russo},
	category = {TI},
	doi = {10.1007/s10915-016-0168-y},
	journal = {Journal of Scientific Computing},
	month = {jan},
	number = {3},
	pages = {975--1001},
	publisher = {Springer Science and Business Media {LLC}},
	title = {High Order Semi-implicit Schemes for Time Dependent Partial Differential Equations},
	volume = {68},
	year = {2016}}

@article{TI_Tavelli2016a,
	author = {Maurizio Tavelli and Michael Dumbser},
	category = {TI},
	doi = {10.1016/j.jcp.2016.05.009},
	journal = {Journal of Computational Physics},
	month = {aug},
	pages = {294--323},
	publisher = {Elsevier {BV}},
	title = {A staggered space{\textendash}time discontinuous {G}alerkin method for the three-dimensional incompressible {N}avier{\textendash}{S}tokes equations on unstructured tetrahedral meshes},
	volume = {319},
	year = {2016}}

@article{TI_Southworth2021a,
	author = {Ben S. Southworth and Oliver Krzysik and Will Pazner},
	category = {TI},
	date = {2021-01-05},
	eprint = {2101.01776},
	eprintclass = {math.NA},
	eprinttype = {arXiv},
	title = {Fast solution of fully implicit Runge-Kutta and discontinuous Galerkin in time for numerical PDEs, Part II: nonlinearities and DAEs},
	year = {2021}}

@inproceedings{Heinecke2016a,
	author = {Alexander Heinecke and Greg Henry and Maxwell Hutchinson and Hans Pabst},
	booktitle = {SC '16: Proceedings of the International Conference for High Performance Computing, Networking, Storage and Analysis},
	doi = {10.1109/sc.2016.83},
	month = {nov},
	pages = {981--991},
	publisher = {{IEEE}},
	title = {{LIBXSMM}: Accelerating Small Matrix Multiplications by Runtime Code Generation},
	url = {https://doi.org/10.1109%2Fsc.2016.83},
	year = 2016}

@article{Piatkowski2018a,
	author = {Marian Piatkowski and Steffen M{\"u}thing and Peter Bastian},
	doi = {10.1016/j.jcp.2017.11.035},
	journal = {Journal of Computational Physics},
	month = {mar},
	pages = {220--239},
	publisher = {Elsevier {BV}},
	title = {A stable and high-order accurate discontinuous Galerkin based splitting method for the incompressible Navier{\textendash}Stokes equations},
	url = {https://doi.org/10.1016%2Fj.jcp.2017.11.035},
	volume = {356},
	year = 2018}

@article{CartonDeWiart2014a,
	author = {C. Carton de Wiart and K. Hillewaert and M. Duponcheel and G. Winckelmans},
	doi = {10.1002/fld.3859},
	journal = {International Journal for Numerical Methods in Fluids},
	month = {nov},
	number = {7},
	pages = {469--493},
	publisher = {Wiley},
	title = {Assessment of a discontinuous Galerkin method for the simulation of vortical flows at high Reynolds number},
	url = {https://doi.org/10.1002%2Ffld.3859},
	volume = {74},
	year = 2014}

@article{Gassner2013a,
	author = {Gregor J. Gassner and Andrea D. Beck},
	doi = {10.1007/s00162-011-0253-7},
	journal = {Theoretical and Computational Fluid Dynamics},
	month = {jan},
	number = {3-4},
	pages = {221--237},
	publisher = {Springer Science and Business Media {LLC}},
	title = {On the accuracy of high-order discretizations for underresolved turbulence simulations},
	url = {https://doi.org/10.1007%2Fs00162-011-0253-7},
	volume = {27},
	year = 2013}

@article{Brachet1991a,
	author = {M. E. Brachet},
	doi = {10.1016/0169-5983(91)90026-f},
	journal = {Fluid Dynamics Research},
	month = {oct},
	number = {1-4},
	pages = {1--8},
	publisher = {{IOP} Publishing},
	title = {Direct simulation of three-dimensional turbulence in the Taylor{\textendash}Green vortex},
	url = {https://doi.org/10.1016%2F0169-5983%2891%2990026-f},
	volume = {8},
	year = 1991}

@article{Taylor1937a,
	author = {Geoffrey Ingram Taylor and Albert Edward Green},
	doi = {10.1098/rspa.1937.0036},
	journal = {Proceedings of the Royal Society of London. Series A - Mathematical and Physical Sciences},
	month = {feb},
	number = {895},
	pages = {499--521},
	publisher = {The Royal Society},
	title = {Mechanism of the production of small eddies from large ones},
	url = {https://doi.org/10.1098%2Frspa.1937.0036},
	volume = {158},
	year = 1937}

@phdthesis{Niemann2018a,
	author = {Martin Niemann},
	date = {2018},
	institution = {Institute of Fluid Mechanics, TU Dresden},
	isbn = {978-3-95908-123-8},
	school = {Institute of Fluid Mechanics, TU Dresden},
	title = {Buoyancy Effects in Turbulent Liquid Metal Flow},
	titleaddon = {A Study in Square Channels and Behind Sudden Expansions},
	year = {2018}}


\end{document}